\documentclass[a4paper,11pt]{amsart}

\usepackage{amsmath,amsthm,amssymb}
\usepackage[mathscr]{eucal}
 \usepackage{cite}
\usepackage{upgreek}
\usepackage[bookmarks,bookmarksnumbered]{hyperref}

\numberwithin{equation}{section}

\theoremstyle{plain}
\newtheorem{main}{Theorem}
\newtheorem{mcor}[main]{Corollary}
\newtheorem{mlemma}[main]{Lemma}

\newtheorem{theorem}{Theorem}[section]

\newtheorem{lemma}[theorem]{Lemma}
\newtheorem{proposition}[theorem]{Proposition}
\newtheorem{corollary}[theorem]{Corollary}
\theoremstyle{definition}
\newtheorem{definition}[theorem]{Definition}
\newtheorem{example}[theorem]{Example}
\newtheorem{notation}[theorem]{Notation}
\newtheorem{remark}[theorem]{Remark}
\newtheorem{fact}[theorem]{Fact}

\begin{document}

\title[Profinite actions with spectral gap]
{Orbit equivalence and Borel reducibility rigidity for \ profinite actions with spectral gap}
\author[Adrian Ioana]{Adrian Ioana}
\thanks{The author\ was partially supported by  NSF  Grant DMS \#1161047, NSF Career Grant DMS \#1253402, and a Sloan Foundation Fellowship.}
\address{Mathematics Department; University of California, San Diego, CA 90095-1555 (United States).}
\email{aioana@ucsd.edu}

\begin{abstract} 
We study equivalence relations $\mathcal R(\Gamma\curvearrowright G)$ that arise from left translation actions of countable groups on their profinite completions. Under the assumption that the action $\Gamma\curvearrowright G$ is free and has spectral gap, we describe precisely when $\mathcal R(\Gamma\curvearrowright G)$ is orbit equivalent or Borel reducible to another such equivalence relation $\mathcal R(\Lambda\curvearrowright H)$.
As a consequence, we provide explicit uncountable families of free ergodic probability measure preserving (p.m.p.) profinite actions of $SL_2(\mathbb Z)$ and its non-amenable subgroups (e.g. $\mathbb F_n$, with $2\leqslant n\leqslant\infty$) whose orbit equivalence relations are mutually  not orbit equivalent and not Borel reducible. 
In particular,  we show that if $S$ and $T$ are distinct sets of primes, then the orbit equivalence relations associated to the actions $SL_2(\mathbb Z)\curvearrowright\prod_{p\in S}SL_2(\mathbb Z_p)$ and $SL_2(\mathbb Z)\curvearrowright\prod_{p\in T}SL_2(\mathbb Z_p)$ are neither orbit equivalent nor Borel reducible.   This settles a conjecture of S. Thomas \cite{Th01,Th06}. Other applications include the first  calculations of outer automorphism groups for concrete treeable  p.m.p. equivalence relations, and the first concrete examples of free ergodic p.m.p. actions of   $\mathbb F_{\infty}$ whose orbit equivalence relations have trivial fundamental group.  
\end{abstract}

\maketitle

\section{Introduction and statement of main results}

\subsection{Introduction}
The general goal of this paper is to establish new rigidity results for countable equivalence relations, in both the measure theoretic and Borel contexts. Our main technical result (Theorem \ref{mainthm}) gives necessary and sufficient conditions for equivalence relations, which are associated to ``translation profinite" actions $\Gamma\curvearrowright G$ with spectral gap, to be orbit equivalent or Borel reducible. The novelty of this theorem lies in that there are no assumptions on the groups, but instead, all the assumptions are imposed on their actions. 
In particular, our result applies to many natural families of translation profinite actions of $SL_2(\mathbb Z)$ and the free groups $\mathbb F_n$, that are known to have spectral gap as a consequence of Selberg's theorem and its recent generalizations  \cite{BG05, BV10}. 

As an application, we provide explicit uncountable families of  free ergodic p.m.p. actions of $SL_2(\mathbb Z)$ and the free groups whose orbit equivalence relations are pairwise not orbit equivalent and not Borel reducible. Additionally, we compute the outer automorphism and fundamental groups of the orbit equivalence relations arising from several of these actions. This improves on several results from \cite{GP03, OP07, PV08, Ga08, Hj08} which showed the existence of such families of actions.

In order to state our results in detail, we first need to review a few concepts concerning countable equivalence relations. Recall that a {\it standard probability space} $(X,\mu)$ is a Polish space $X$ equipped with its $\sigma$-algebra of Borel sets and a Borel probability measure $\mu$.

Countable p.m.p. and countable Borel equivalence relations are fundamental objects of study in orbit equivalence  theory and descriptive set theory. 
If $\Gamma\curvearrowright (X,\mu)$ is a p.m.p. (respectively, Borel) action of a countable group $\Gamma$ on a standard probability space $(X,\mu)$, then the  {\it orbit equivalence relation} 
$\mathcal R(\Gamma\curvearrowright X):=\{(x,y)\in X^2\;|\;\Gamma\cdot x=\Gamma\cdot y\}$ 
is a countable p.m.p. (respectively, Borel) equivalence relation. Conversely, it was shown in \cite[Theorem 1]{FM77} that any countable p.m.p. and any countable Borel  equivalence relation can be realized in 
this way. 

The study of countable equivalence relations is organized using the notions of orbit equivalence and Borel reducibility.
Firstly, let $\mathcal R$, $\mathcal S$  be countable p.m.p. equivalence relations on standard probability spaces $(X,\mu)$, $(Y,\nu)$. Then $\mathcal R$ is said to be {\it orbit equivalent} to $\mathcal S$ if there exists an isomorphism of probability spaces $\theta:(X,\mu)\rightarrow (Y,\nu)$ such that $(\theta\times\theta)(\mathcal R)=\mathcal S$. Moreover, we say that $\mathcal R$, $\mathcal S$ are {\it stably orbit equivalent} if there exist Borel subsets $A\subset X,B\subset Y$ of positive measure such that the restrictions $\mathcal R_{|A}$, $\mathcal S_{|B}$ are orbit equivalent.
Two p.m.p. actions $\Gamma\curvearrowright (X,\mu)$, $\Lambda\curvearrowright (Y,\nu)$ are (stably) orbit equivalent if their orbit equivalence relations are (stably) orbit equivalent.

Secondly, let $\mathcal R$, $\mathcal S$ be countable Borel equivalence relations on standard Borel spaces $X$, $Y$. Then $\mathcal R$ is  {\it Borel reducible} to 
$\mathcal S$ if there exists a Borel map $\theta:X\rightarrow Y$ such that $(x,y)\in\mathcal R$ iff $(\theta(x),\theta(y))\in\mathcal S$.
The condition that $\mathcal R$ is Borel reducible to $\mathcal S$ is usually interpreted to mean that the classification problem associated to $\mathcal R$ is at most as complicated as the classification problem associated to $\mathcal S$.

\subsection{Orbit equivalence and Borel reducibility rigidity} 
We are now ready to state the main technical result of this paper. 
 Recall that an ergodic p.m.p. action $\Gamma\curvearrowright (X,\mu)$ is said to have {\it spectral gap} if the  Koopman representation of $\Gamma$ on $L^2(X,\mu)\ominus\mathbb C1$ has no almost invariant vectors. 
 Let $\Gamma$ be a residually finite group and  $G=\varprojlim\Gamma/\Gamma_n$ be its profinite completion with respect to a descending chain $\{\Gamma_n\}_n$ of finite index, normal subgroups with trivial intersection, $\cap_n\Gamma_n=\{e\}$. Let $\rho:\Gamma\hookrightarrow G$ be the embedding given by $\rho(g)=(g\Gamma_n)_n$. 
Then the {\it left translation} action $\Gamma\curvearrowright G$ defined by $g\cdot x=\rho(g)x$ is free, ergodic, and preserves the Haar measure of $G$.
Moreover, this action is {\it profinite}, i.e. it is an inverse limit of actions of $\Gamma$ on finite probability spaces.
 
\begin{main}\label{mainthm}
 Let $\Gamma$, $\Lambda$ be residually finite groups.  Let $G=\varprojlim\Gamma/\Gamma_n$, $H=\varprojlim\Lambda/\Lambda_n$ be  profinite completions of  $\Gamma$, $\Lambda$ with respect to descending chains $\{\Gamma_n\}_n$, $\{\Lambda_n\}_n$ of finite index normal subgroups with trivial intersection. Denote by $m_G$ and $m_H$ the Haar measures of $G$ and $H$.
Assume  that the left translation action $\Gamma\curvearrowright (G,m_G)$ has spectral gap.
Then the following hold:
\begin{enumerate}
\item    $\Gamma\curvearrowright (G,m_G)$ is stably orbit equivalent (respectively, orbit equivalent) to $\Lambda\curvearrowright (H,m_H)$ if and only if we can find open subgroups $G_0<G$ and $H_0<H$ and a continuous isomorphism $\delta:G_0\rightarrow H_0$ such that $\delta(\Gamma\cap G_0)=\Lambda\cap H_0$ (and, respectively,  $[G:G_0]=[H:H_0]$).
\item  $\mathcal R(\Gamma\curvearrowright G)$ is Borel reducible to $\mathcal R(\Lambda\curvearrowright H)$ if and only if we can find an open subgroup $G_0<G$, a closed subgroup $H_0<H$, and a continuous isomorphism $\delta:G_0\rightarrow H_0$ such that $\delta(\Gamma\cap G_0)=\Lambda\cap H_0$.
\end{enumerate}
\end{main}

Next, let us make a few comments on the statement of Theorem \ref{mainthm} and discuss some classes of actions to which it applies.

Assume that $\Gamma$ has property (T) of Kazhdan (e.g. take $\Gamma=SL_n(\mathbb Z)$, for  $n\geqslant 3$). Then every ergodic p.m.p. action $\Gamma\curvearrowright (X,\mu)$ has spectral gap. In particular, the first part of Theorem \ref{mainthm} holds for any left translation action of $\Gamma$  on one of its profinite completions $G$. Let us point out that in this case much more can be said. Indeed,  if $\Gamma$ has property (T), we showed that the action $\Gamma\curvearrowright G$ is  {\it orbit equivalent superrigid}, in the sense that any free ergodic p.m.p. action $\Lambda\curvearrowright (Y,\nu)$ that is orbit equivalent to $\Gamma\curvearrowright G$ is virtually conjugate to it (see \cite[Theorem A]{Io08}). This means that we can find a finite index subgroup $\Lambda_0<\Lambda$ and an open subgroup $G_0<G$  such that  an ergodic component of the action $\Lambda_0\curvearrowright (Y,\nu)$ is conjugate to the action $\Gamma\cap G_0\curvearrowright G_0$. Moreover,  note that if $\Gamma$ has property (T), then the action $\Gamma\curvearrowright G$ satisfies a cocycle superrigidity theorem \cite[Theorem B]{Io08} which can be used to deduce the second part of Theorem \ref{mainthm}.

In this paper, however, we are interested in studying actions of
groups such as $SL_2(\mathbb Z)$ and the free groups $\mathbb F_n$, for which property (T) fails.
Nevertheless,  these groups 
still satisfy a weak form of property (T). Recall that a countable group $\Gamma$ has {\it property $(\uptau)$ with respect to a family of subgroups $\{\Gamma_n\}_n$} if the unitary representation of $\Gamma$ on $\bigoplus_n\ell^2_0(\Gamma/\Gamma_n)$ has no almost invariant vectors, where $\ell^2_0(\Gamma/\Gamma_n)=\ell^2(\Gamma/\Gamma_n)\ominus\mathbb C1$ (see \cite{LZ03,Lu12}).

Selberg's famous theorem implies that  $SL_2(\mathbb Z)$ has property ($\uptau$) with respect to the family of congruence subgroups $\Gamma(n)=\ker(SL_2(\mathbb Z)\rightarrow SL_2(\mathbb Z/n\mathbb Z))$. 
Recently, Selberg's theorem has been vastly generalized starting with the breakthrough work of J. Bourgain and A. Gamburd in \cite{BG05}. In particular, J. Bourgain and P. Varj\'{u} have shown that any non-amenable subgroup $\Gamma<SL_2(\mathbb Z)$ has property ($\uptau$) with respect to the family of subgroups $\{\Gamma\cap\Gamma(n)\}_{n\geqslant 1}$ (see \cite[Theorem 1]{BV10}).
Note that a translation action $\Gamma\curvearrowright  \varprojlim\Gamma/\Gamma_n$ has spectral gap if and only if $\Gamma$ has property ($\uptau$) with respect to $\{\Gamma_n\}_{n}$. Thus, the results mentioned above can be used to construct many interesting families of translation actions with spectral gap.
 We single out two such families of actions that we will use subsequently:

{\bf Examples.}
Consider the profinite groups $G_S=\prod_{p\in S}SL_2(\mathbb F_p)$ and $K_S=\prod_{p\in S}SL_2(\mathbb Z_p)$, for a set of primes $S$. 
Whenever $S$ is infinite, we view $SL_2(\mathbb Z)$ as a subgroup of $G_S$, via the diagonal embedding. We also embed $SL_2(\mathbb Z)$ diagonally into $K_S$, for any set $S$. 
By the Strong Approximation Theorem (see e.g. \cite{LS03}) both of these embeddings are dense.

Given a subgroup $\Gamma<SL_2(\mathbb Z)$, we denote by $G_{\Gamma,S}$ and $K_{\Gamma,S}$ its closures in $G_S$ and $K_S$, respectively.
If $\Gamma$ is  non-amenable, then the translation actions $\Gamma\curvearrowright G_{\Gamma,S}$, $\Gamma\curvearrowright K_{\Gamma,S}$ have spectral gap \cite{BV10}.
Moreover, the Strong Approximation Theorem  implies that 
 $G_{\Gamma,S}<G_S$ and $K_{\Gamma,S}<K_S$ are open.

Let us point out that the closures $G_{\Gamma,S}$ and $K_{\Gamma,S}$ of $\Gamma$ are, at least in principle, explicit. In general, if $G=\varprojlim\Gamma/\Gamma_n$ is a profinite completion of a countable group $\Gamma$ and $\Gamma_0<\Gamma$ is a subgroup, then the closure of $\Gamma_0$ in $G$ is isomorphic, as a topological group, to the profinite group $\varprojlim\Gamma_0/(\Gamma_0\cap\Gamma_n)$.

{\bf Remarks.}
\begin{itemize}
\item The conclusion of the first part of Theorem \ref{mainthm} is optimal, in the sense that the translation action $\Lambda\curvearrowright (H,m_H)$ cannot be replaced with an arbitrary free ergodic p.m.p. action $\Lambda\curvearrowright (Y,\nu)$. 
In other words, the action $\Gamma\curvearrowright (G,m_G)$ is not necessarily orbit equivalent superrigid.
This is easy to see if $\Gamma$ is a free group. Indeed, by \cite[Theorem 2.27]{MS02}, any free ergodic p.m.p. action of  a free group is orbit equivalent to actions of uncountably many non-isomorphic groups, and hence cannot be orbit equivalent superrigid. 
\item If $\Gamma$ is the product of two groups having property ($\uptau$) and the action $\Gamma\curvearrowright G$ satisfies a certain growth condition, then the conclusion of Theorem \ref{mainthm} can be deduced from the cocycle rigidity result from \cite[Theorem C]{OP08}.
\item 
An analogous result to the first part of Theorem \ref{mainthm}, where orbit equivalence is replaced with weak equivalence of actions, was obtained recently
by M. Ab\'{e}rt and G. Elek in \cite[Theorem 2]{AE10}. Their result in particular shows that, under the assumptions of Theorem \ref{mainthm}, the actions $\Gamma\curvearrowright (G,m_G)$, $\Lambda\curvearrowright (H,m_H)$ are weakly equivalent if and only if they are conjugate.
\item The notion of ``spectral gap rigidity" has been introduced by S. Popa in the context of von Neumann algebras and has been used to great effect starting with the work \cite{Po06a}. 
\item
Several rigidity results were recently  proven in \cite{Po09} for II$_1$ factors $M$ that are an inductive limit of a sequence $\{M_n\}$ of subfactors  with spectral gap.   If $\Gamma$  is not inner amenable and 
$\Gamma\curvearrowright X=\varprojlim X_n$ is a profinite action with spectral gap, then the II$_1$ factor $M=L^{\infty}(X)\rtimes\Gamma$ has this property, where 
$M_n=L^{\infty}(X_n)\rtimes\Gamma$. Note, however, that one cannot directly apply \cite[Theorem 3.5]{Po09} as the relative commutant condition $(M_n'\cap M)'\cap M=M_n$ fails.

 \item We do not know to what extent Theorem \ref{mainthm} can be extended to arbitrary ergodic compact (or profinite) actions. Recall that these are actions of the form $\Gamma\curvearrowright G/K$, where $G$ is a compact (respectively, compact profinite) group in which $\Gamma$ embeds densely, and $K<G$ is a closed subgroup.
However, the proof of Theorem \ref{mainthm} relies on a cocycle rigidity result (Theorem \ref{cocycle}) for translation profinite actions
 which admits an analogue in the case of translation actions $\Gamma\curvearrowright G$ on compact connected groups $G$ with finite fundamental group (Theorem \ref{compact}). As a consequence, we are able to prove an analogue of Theorem \ref{mainthm} in this case (see Corollary \ref{Bor} and Theorem \ref{conjugacy}).
 Moreover, Theorem \ref{conjugacy} establishes an analogue of Theorem \ref{mainthm} for a fairly general class of compact actions.
 \end{itemize}

\subsection{Orbit inequivalent and Borel incomparable actions}
In the rest of the introduction, we discuss several applications of  Theorem \ref{mainthm} and of its method of proof. Our first applications provide examples  of actions of the free groups that are orbit inequivalent and Borel incomparable.
We begin by giving some context and motivation.

In the last 15 years, remarkable progress has been made in the study of countable equivalence relations, in both  the measure theoretic and Borel contexts (see the surveys \cite{Po06b,Fu09,Ga10} and \cite{Th06,TS07}).  In particular, considerable effort has been devoted to the investigation of equivalence relations that arise from actions of free groups and, more generally, of {\it treeable} equivalence relations (i.e. equivalence relations whose classes are the connected components of a Borel acyclic graph).

On the orbit equivalence side,  D. Gaboriau  proved that free ergodic p.m.p. actions $\mathbb F_n\curvearrowright X$, $\mathbb F_m\curvearrowright Y$ of free groups of different ranks ($n\not=m$) are never orbit equivalent \cite{Ga99,Ga01}. This result, however, offered little insight on how to distinguish between actions of the same free group. It was not until the work of D. Gaboriau and S. Popa \cite{GP03} that every non-abelian free group $\mathbb F_n$ ($2\leqslant n\leqslant\infty$)  was shown to admit a continuum of non-orbit equivalent free ergodic p.m.p. actions. However,  \cite{GP03} only demonstrates the existence of such a continuum of actions. This motivated our work \cite{Io06}, where we found an explicit 
 list of uncountably many orbit inequivalent  free ergodic p.m.p. actions of $\mathbb F_n$. 
Note that finding natural classes of orbit inequivalent actions of $\mathbb F_n$ is a  difficult task, since some obvious candidates  turn out to be orbit equivalent.  Indeed, L. Bowen proved that any two  Bernoulli actions of $\mathbb F_n$ are  orbit equivalent, if $2\leqslant n\leqslant\infty$, and moreover that any two Bernoulli actions of $\mathbb F_n,\mathbb F_m$ are stably orbit equivalent, whenever $2\leqslant n,m<\infty$ \cite{Bo09a,Bo09b} (see also \cite{MRV11}).

In descriptive set theory, the investigation of treeable equivalence relations started in \cite{JKL01} where it was proved that the orbit equivalence relation $E_{\infty\mathcal T}$ of the free part of the Bernoulli action $\mathbb F_{\infty}\curvearrowright \{0,1\}^{\mathbb F_{\infty}}$ is maximal among treeable countable Borel  equivalence relations, with respect to Borel reducibility. In the same paper, the authors asked whether there exist infinitely many  treeable countable Borel equivalence, up to Borel reducibility. 
After providing a first example of a non-hyperfinite treeable equivalence relation that lies strictly below $E_{\infty\mathcal T}$ \cite{Hj03},  G. Hjorth answered this question in the affirmative in \cite{Hj08}. More precisely, he proved that there exist uncountably many treeable countable Borel equivalence relations that are mutually incomparable with respect to Borel reducibility.
Since the proof of \cite{Hj08} uses a separability argument, it only provides an existence result, leaving open the problem of finding specific treeable equivalence relations that are Borel incomparable. In fact, at the time of the writing, not a single example of a pair of treeable countable Borel equivalence relations such that neither is Borel reducible to the other was known.

As a first application of Theorem \ref{mainthm}, we exhibit the first concrete uncountable families of actions of $SL_2(\mathbb Z)$ and the free groups $\mathbb F_n$ that are neither orbit equivalent nor Borel reducible. 

\begin{mcor}\label{corB}  Let $S$, $T$ be infinite sets of primes, and $\Gamma,\Lambda<SL_2(\mathbb Z)$ be non-amenable subgroups.

\begin{enumerate}
\item If $S\not=T$, then the actions $SL_2(\mathbb Z)\curvearrowright G_S$ and $SL_2(\mathbb Z)\curvearrowright G_T$ are not  stably orbit equivalent,  and the equivalence relations $\mathcal R(SL_2(\mathbb Z)\curvearrowright G_{S})$ and $\mathcal R(SL_2(\mathbb Z)\curvearrowright G_{T})$ are not comparable with respect to Borel reducibility.

\item If  $|S\Delta T|=\infty$, then the actions $\Gamma\curvearrowright G_{\Gamma,S}$ and $\Lambda\curvearrowright G_{\Lambda,T}$ are not  stably orbit equivalent,
and the equivalence relations $\mathcal R(\Gamma\curvearrowright G_{\Gamma,S})$ and $\mathcal R(\Lambda\curvearrowright G_{\Lambda,T})$ are not comparable with respect to Borel reducibility.
\end{enumerate}
\end{mcor}

Corollary \ref{corB} improves on a result of N. Ozawa and S. Popa who were the first to show the existence of uncountably many orbit inequivalent free ergodic {profinite} actions of $\mathbb F_n$  \cite[Theorem 5.4]{OP07}. More precisely, in the context of Corollary \ref{corB}, they proved that 
 if $\{S_i\}_{i\in I}$ is an uncountable family of infinite sets of primes such that $|S_i\cap S_j|<\infty$, for all $i\not=j$, then among the actions $\{\Gamma\curvearrowright G_{\Gamma,S_i}\}_{i\in I}$ there are uncountably many orbit equivalence classes. 

In \cite[Conjecture 5.7]{Th01} and  \cite[Conjecture 2.14]{Th06}, S. Thomas proposed a scenario for providing an explicit uncountable family of treeable countable Borel equivalence relations that are pairwise incomparable with respect to Borel reducibility. Specifically, he conjectured  
 that if $S,T$ are distinct nonempty sets of primes, then the orbit equivalence relations of the actions of $SL_2(\mathbb Z)$ on $K_S$, $K_T$ are neither Borel reducible nor stably orbit equivalent. 
Note that the analogous results where $SL_2$  is replaced by $SL_n$, for some $n\geqslant 3$, were shown to hold in \cite{Th01} and \cite{GG88}, respectively.

As a consequence of Theorem \ref{mainthm}, we settle this conjecture and, more generally, show that:

\begin{mcor}\label{corC} Let $S$, $T$ be nonempty sets of primes, and $\Gamma,\Lambda<SL_2(\mathbb Z)$ non-amenable subgroups.

If  $S\not=T$, then the actions $\Gamma\curvearrowright K_{\Gamma,S}$ and $\Lambda\curvearrowright K_{\Lambda,T}$ are not stably orbit equivalent, and  the equivalence relations $\mathcal R(\Gamma\curvearrowright K_{\Gamma,S})$ and $\mathcal R(\Lambda\curvearrowright K_{\Lambda,T})$ are not comparable with respect to Borel reducibility.
 \end{mcor}

{\bf Remarks.}
 \begin{itemize}
 \item Since the matrices $a=\bigl(\begin{smallmatrix}1&2\\ 0&1\end{smallmatrix} \bigr)$ and $b=\bigl(\begin{smallmatrix}1&0\\ 2&1\end{smallmatrix} \bigr)$ generate an isomorphic copy of  $\mathbb F_2$ inside $SL_2(\mathbb Z)$, we have an embedding of $\mathbb F_n$ into $SL_2(\mathbb Z)$, for any $2\leqslant n\leqslant\infty$. Therefore, both  Corollaries \ref{corB} and \ref{corC} yield concrete uncountable families of pairwise orbit inequivalent and Borel incomparable free ergodic p.m.p. actions of $\mathbb F_n$.  
 \item The orbit inequivalent actions of $\mathbb F_n$ provided by Corollaries \ref{corB} and \ref{corC} are different from the ones found in \cite{GP03, Io06}. Indeed, the latter actions admit quotients that have the relative property (T) \cite{Po01} and hence are not orbit equivalent to profinite actions.
 \item By a well-known result from \cite{OP07}, if  $\Gamma\curvearrowright (X,\mu)$ is a free ergodic p.m.p. profinite action of a non-amenable subgroup of $SL_2(\mathbb Z)$, then the II$_1$ factor  $L^{\infty}(X)\rtimes\Gamma$ has a unique Cartan subalgebra, up to unitary conjugacy. This implies that the non-orbit equivalent actions given by Corollaries \ref{corB} and \ref{corC}  have in fact non-isomorphic II$_1$ factors.
 
 \end{itemize}

Our techniques also allow us to confirm a second conjecture of S. Thomas \cite[Conjecture 6.10]{Th01}. This conjecture asserts that the orbit equivalence relations arising from the natural actions of $SL_2(\mathbb Z)$ on the projective lines over the various $p$-adic fields are pairwise incomparable with respect to Borel reducibility. More generally, we prove:

\begin{mcor}\label{corD}
For a prime $p$, let $PG(1,\mathbb Q_p)=\mathbb Q_p\cup\{\infty\}$ be the projective line over the field $\mathbb Q_p$ of $p$-adic numbers. Consider the action of $GL_2(\mathbb Q_p)$ on $PG(1,\mathbb Q_p)$ by linear fractional transformations. Let $\Gamma<GL_2(\mathbb Q_p)$ and $\Lambda<GL_2(\mathbb Q_q)$ be countable subgroups, for some primes $p$ and $q$.

If $p\not=q$ and $\Gamma\cap SL_2(\mathbb Z)$ is non-amenable, then $\mathcal R(\Gamma\curvearrowright PG(1,\mathbb Q_p))$ is not Borel reducible to $\mathcal R(\Lambda\curvearrowright PG(1,\mathbb Q_q))$.
\end{mcor}

Corollary \ref{corD} strengthens and offers a different approach to a theorem of G. Hjorth and S. Thomas. 
Thus, the case $\Gamma=\Lambda=GL_2(\mathbb Q)$ of Corollary \ref{corD} recovers the main result of \cite{HT04}. The original motivation for this result came from the classification problem for torsion-free abelian groups of finite rank (for more on this, see the survey \cite{Th06}). More precisely, by \cite[Theorem 5.7]{Th00}, the main result of \cite{HT04} can be reformulated as follows:  if $p\not=q$ are distinct primes, then the classification problems for the $p$-local and $q$-local torsion-free abelian groups of rank $2$ are incomparable with respect to Borel reducibility.

\subsection{Calculations of outer automorphism and fundamental groups}
The second direction in which our techniques apply is the calculation of the {\it fundamental group} $\mathcal F(\mathcal R)$ and the {\it outer automorphism group} $Out(\mathcal R)$ of orbit equivalence relations $\mathcal R$ associated to actions of free groups  (see Section \ref{basic} for the definitions). 
Let us briefly recall known results along these lines.

First, it was shown in \cite{Ga99,Ga01}  that if $2\leqslant n<\infty$, then the fundamental group of $\mathcal R(\mathbb F_n\curvearrowright X)$ is trivial, for any free ergodic p.m.p. action of $\mathbb F_n$.
On the other hand, no calculations of fundamental groups of equivalence arising from actions of $\mathbb F_{\infty}$ or outer automorphism groups of equivalence relations arising from actions of $\mathbb F_n$ ($2\leqslant n\leqslant\infty$) were available until recently. S. Popa and S. Vaes then proved that the fundamental group of $\mathcal R(\mathbb F_{\infty}\curvearrowright X)$ can be equal to any countable subgroup as well as many uncountable subgroups of $\mathbb R_{+}^*$ \cite{PV08}. Moreover, they showed the existence of free ergodic p.m.p. actions $\mathbb F_{\infty}\curvearrowright X$ such that $\mathcal R(\mathbb F_{\infty}\curvearrowright X)$  has trivial outer automorphism group. Shortly after, D. Gaboriau proved that the outer automorphism of $\mathcal R(\mathbb F_n\curvearrowright X)$ can in fact be trivial for any $2\leqslant n\leqslant\infty$  \cite{Ga08}. Nevertheless, the proofs of \cite{PV08,Ga08} did not provide a single explicit example of an action of $\mathbb F_{\infty}$ (respectively, of $\mathbb F_n$, with $2\leqslant n\leqslant\infty$) for which the fundamental group (respectively, the outer automorphism group) of the orbit equivalence relation could be computed.

We address this problem here by exhibiting the first examples of free ergodic p.m.p. actions of $\mathbb F_{\infty}$ whose orbit equivalence relation has trivial fundamental group, and the first examples of actions of $\mathbb F_2$ for which the outer automorphism group of the orbit equivalence relation can be calculated.

\begin{mcor}\label{corE}
Let $S$ be an infinite set of primes and $\Gamma<SL_2(\mathbb Z)$  be a non-amenable subgroup. 

Then we have the following:
\begin{enumerate}
\item $Out(\mathcal R(SL_2(\mathbb Z)\curvearrowright G_S))\cong (G_S/Z)\rtimes\mathbb Z/2\mathbb Z$. Here, $Z=\{\pm I\}$ denotes the center of $G_S$ and the action of $\mathbb Z/2\mathbb Z$ on $G_S$ is given by conjugation with the matrix $\bigl(\begin{smallmatrix}1&0\\ 0&-1\end{smallmatrix} \bigr)\in GL_2(\mathbb Z).$
\item  The equivalence relation $\mathcal R(\Gamma\curvearrowright G_{\Gamma,S})$ and the associated II$_1$ factor $L^{\infty}(G_{\Gamma,S})\rtimes\Gamma$ have trivial fundamental groups, i.e. $\mathcal F(\mathcal R(\Gamma\curvearrowright G_{\Gamma,S}))=\mathcal F(L^{\infty}(G_{\Gamma,S})\rtimes\Gamma)=\{1\}$.
\end{enumerate}
\end{mcor}

{\bf Remarks.} 
\begin{itemize}
\item Consider a fixed embedding $\mathbb F_{\infty}\subset SL_2(\mathbb Z)$. Then the combination of Corollaries \ref{corC} and  \ref{corE} yields a continuum of pairwise non-stably orbit equivalent free ergodic p.m.p. actions of $\mathbb F_{\infty}$  whose equivalence relations and II$_1$ factors have trivial fundamental groups.
\item
Let us explain how Corollary \ref{corE} also leads to examples of orbit equivalence relations of free ergodic p.m.p. actions of $\mathbb F_2$ whose outer automorphism group can be explicitly computed. Let $S$ be an infinite set of primes containing $2$. Denote by $G_{S}(2)$ the kernel of the natural homomorphism $G_S\rightarrow SL_2(\mathbb Z/2\mathbb Z)$. Note that $Z\subset G_S(2)$ and that $SL_2(\mathbb Z)\cap G_S(2)$ is  the congruence subgroup $\Gamma(2)$. If $A\subset G_S(2)$ is a fundamental domain for the left translation action $Z\curvearrowright G_S(2)$, it is easy to see that $\mathcal R(SL_2(\mathbb Z)\curvearrowright G_S)_{|A}\cong\mathcal R(\Gamma(2)/Z\curvearrowright G_S(2)/Z)$. Since the quotient $\Gamma(2)/Z$ is isomorphic to $\mathbb F_2$ and the outer automorphism group  is insensitive to taking restrictions (i.e. $Out(\mathcal R_{|A})\cong Out(\mathcal R)$, for any ergodic countable p.m.p. equivalence relation $\mathcal R$),
we get examples of actions of $\mathbb F_2$ with the desired property.
\item Let $n$ be a natural number of the form $n=p_1p_2...p_k$, where $3\leqslant p_1<p_2<...<p_k$ are prime numbers.  Let $S$ be an infinite set of primes containing $p_1,p_2,...,p_k$. Denote by $G_S(n)$ the kernel of the natural homomorphism $G_S\rightarrow\prod_{i=1}^kSL_2(\mathbb Z/p_i\mathbb Z)$. Then we have that $SL_2(\mathbb Z)\cap G_S(n)=\Gamma(n)$ and $\mathcal R(SL_2(\mathbb Z)\curvearrowright G_S)_{|G_S(n)}=\mathcal R(\Gamma(n)\curvearrowright G_S(n))$. Since $\Gamma(n)$ is a free group (of rank $1+\frac{1}{12}\prod_{i=1}^kp_i(p_i^2-1)$), by Corollary \ref{corE} we get more examples of actions of free groups with the same property as above.
\item Let $\Gamma$ be a free group and $\Gamma\curvearrowright^{\sigma} (X,\mu)$ be a free ergodic p.m.p. profinite action. Denote $\mathcal R=\mathcal R(\Gamma\curvearrowright X)$ and $M=L^{\infty}(X)\rtimes\Gamma$. Then the main result of \cite{OP07} implies that $L^{\infty}(X)$ is the unique Cartan subalgebra of $M$, up to unitary conjugacy. As a consequence,  $Out(M)\cong H^1(\sigma)\rtimes Out(\mathcal R)$, where $H^1(\sigma)$ denotes the first cohomology group of $\sigma$ with values in $\mathbb T$. Thus, if $\sigma$ is any of the actions from the previous two remarks, then the outer automorphism of its II$_1$ factor $M$ can be computed explicitly in terms of $H^1(\sigma)$.

\end{itemize}
 
Our next result provides further computations of fundamental and outer autormorphism groups.

\begin{mcor}\label{corF}
Let $p$ be a prime number and $\Gamma<SL_2(\mathbb Z)$ be a non-amenable subgroup.

Then we have the following
\begin{enumerate}
\item $Out(\mathcal R(SL_2(\mathbb Z)\curvearrowright SL_2(\mathbb Z_p)))$ is  a $(\mathbb Z/2\mathbb Z)^2$-extension of $PSL_2(\mathbb Q_p)$, i.e. there is an (explicit) exact sequence $1\rightarrow PSL_2(\mathbb Q_p)\rightarrow Out(\mathcal R(SL_2(\mathbb Z)\curvearrowright SL_2(\mathbb Z_p)))\rightarrow (\mathbb Z/2\mathbb Z)^2\rightarrow 1$.
\item $\mathcal F(\mathcal R(\Gamma\curvearrowright K_{\Gamma}))=\mathcal F(L^{\infty}(K_{\Gamma})\rtimes\Gamma)=\{1\}$, where $K_{\Gamma}$ denotes the closure of $\Gamma$ in $SL_2(\mathbb Z_p)$.
\end{enumerate}
\end{mcor}

Note that a statement similar to the first part of Corollary \ref{corF}, where $SL_2$ is replaced with $SL_n$, for some $n\geqslant 3$,  was obtained by A. Furman in \cite[Theorem 1.6]{Fu03}.

\subsection{Treeable equivalence relations with trivial outer automorphism group}
Corollaries \ref{corE} and \ref{corF} leave open the problem of finding concrete actions of $\mathbb F_n$ for which $\mathcal R(\mathbb F_n\curvearrowright X)$ has trivial outer automorphism group. Furthermore, they do not provide any examples of treeable equivalence relations with trivial outer automorphism groups. 
To approach these problems, we show that, under fairly general assumptions on a countable subgroup $\Gamma<SO(3)$, 
 the outer automorphism group of the equivalence relations associated to the natural actions of $\Gamma$  on  $SO(3),S^2$ and $P^2(\mathbb R)$ can be explicitly expressed in terms of the normalizer of $\Gamma$ in $SO(3)$.

\begin{mcor}\label{corG} Let $\Gamma$ be a countable icc  subgroup of the special orthogonal group $G=SO(3)$.
Assume that $\Gamma$ contains matrices $g_1,...,g_k$  which have algebraic entries and generate a dense subgroup of $G$.

Consider the free ergodic  p.m.p. actions $\Gamma\curvearrowright (G,m)$, $\Gamma\curvearrowright (S^2,\lambda)$ and  $\Gamma\curvearrowright (P^2(\mathbb R),\mu)$, where $m$ denotes the Haar measure of $G$, while $\lambda$ and $\mu$ denote the Lebesgue measures of the $2$-dimensional sphere, $S^2$, and the $2$-dimensional real projective space, $P^2(\mathbb R)$. 

Then we have that $Out(\mathcal R(\Gamma\curvearrowright G))=N_G(\Gamma)/\Gamma\times G$,
$Out(\mathcal R(\Gamma\curvearrowright S^2))\cong N_G(\Gamma)/\Gamma\times(\mathbb Z/2\mathbb Z)$ and
 $Out(\mathcal R(\Gamma\curvearrowright P^2(\mathbb R)))\cong N_G(\Gamma)/\Gamma$, where $N_G(\Gamma)$ denotes the normalizer of $\Gamma$ in $G$.

Moreover, $\mathcal F(\mathcal R(\Gamma\curvearrowright G))=\mathcal F(\mathcal R(\Gamma\curvearrowright S^2))=\mathcal F(\mathcal R(\Gamma\curvearrowright P^2(\mathbb R)))=\{1\}$.
\end{mcor}

The work of J. Bourgain and A. Gamburd \cite{BG06} implies that, under the assumptions imposed on $\Gamma$, the left translation action $\Gamma\curvearrowright (G,m_G)$ has spectral gap. This fact will be a key ingredient in the proof of Corollary \ref{corG}.

In view of Corollary \ref{corG}, in order to give examples of actions of $\Gamma=\mathbb F_n$  whose equivalence relations have trivial outer automorphism group, 
it suffices to a find a copy of $\Gamma$ inside $G$ which is generated by matrices with algebraic entries and  has trivial normalizer. We were, however, unable to compute the normalizer of $\Gamma$, for any of the known ``algebraic" embeddings of $\Gamma$ into $G$ (see Remark \ref{Swi}).

Nevertheless, we managed to calculate the normalizer of certain countable subgroups $\Gamma<G$ (see Corollary \ref{cor}). These groups $\Gamma$, although are not free, are treeable, in the sense that any 
equivalence relation  $\mathcal R(\Gamma\curvearrowright X)$ arising from a free p.m.p. action of $\Gamma$ is treeable.  
Moreover, for some of these groups we showed that $N_G(\Gamma)=\Gamma$. This leads to the first concrete examples of treeable countable p.m.p. equivalence 
relations which have trivial outer automorphism group.

\begin{mcor}\label{corH}
Let $p\geqslant 4$, $q\geqslant 6$ be even integers such that $p\not=q$ and $q=2s$, with $s$ odd.

Denote by $\Gamma$ the subgroup of $G=SO(3)$ generated by  the rotation about the $x$-axis  by angle $\frac{2\pi}{p}$ and  the rotation about 
the $z$-axis by angle $\frac{2\pi}{q}$.

 Then  $Out(\mathcal R(\Gamma\curvearrowright G))\cong G$,
  $Out(\mathcal R(\Gamma\curvearrowright S^2))\cong \mathbb Z/2\mathbb Z$ and $Out(\mathcal R(\Gamma\curvearrowright P^2(\mathbb R)))\cong \{e\}$.
  
  Moreover, $\mathcal F(\mathcal R(\Gamma\curvearrowright G))=\mathcal F(\mathcal R(\Gamma\curvearrowright S^2))=\mathcal F(\mathcal R(\Gamma\curvearrowright P^2(\mathbb R)))=\{1\}$.
\end{mcor}

 C. Radin and L. Sadun \cite{RS98} showed that $\Gamma$ is isomorphic to an amalgamated free product of the form $D_p*_{D_2}D_q$, where $D_n=\mathbb Z/n\mathbb Z\rtimes\mathbb Z/2\mathbb Z$ denotes the dihedral group of order $2n$. By a result of D. Gaboriau \cite{Ga99} it follows that 
 $\Gamma$ is treeable, showing that the equivalence relations from Corollary \ref{corH} are indeed  treeable.

\subsection{Comments on the proofs} Since all the results stated above are derived  from either Theorem \ref{mainthm} or variations of it, let us give an outline of its proof.
There are two main ingredients in the proof of Theorem \ref{mainthm}.

The first is a criterion for untwisting cocycles of translation profinite actions  (Theorem \ref{cocycle}). As we explain in the next paragraphs,  this criterion shows that any cocycle satisfying a certain local condition is essentially  cohomologous to a homomorphism.

 In  \cite[Theorem B]{Io08}, we proved a cocycle superrigidity theorem for profinite actions. 
 In \cite[Theorem 5.21]{Fu09}, A. Furman provided an alternative proof of \cite{Io08}. His proof applies to the wider the class of compact actions and  has partially inspired our approach.
  
 The main result of \cite{Io08}, in the formulation given in \cite{Fu09},  shows that if $\Gamma$ has property (T) and $\Gamma\curvearrowright G$ is a translation profinite action,  then any cocycle $w:\Gamma\times G\rightarrow\Lambda$ taking values into a countable group $\Lambda$ is ``virtually cohomologous to a homomorphism". More precisely, we can find an open subgroup $G_0<G$ and a homomorphism $\delta:\Gamma\cap G_0\rightarrow\Lambda$ such that the restriction of $w$ to $(\Gamma\cap G_0)\times G_0$ is cohomologous to $\delta$.

In both \cite{Io08} and \cite{Fu09}, one combines property (T)  with S. Popa's deformation/rigidity approach to deduce that $w$ satisfies a ``local uniformity" condition, which is then exploited to conclude that $w$ is virtually cohomogolous to a homomorphism.
This condition, in the form which appears in \cite{Fu09}, amounts to the existence of a neighborhood $V$ of the identity in $G$ and of a constant $C\in (\frac{31}{32},1)$ such that
  $m_G(\{x\in G\;|\;w(g,xt)=w(g,x)\})\geqslant C$, for all $g\in\Gamma$ and every $t$ in $V$.

The second ingredient in the proof of Theorem \ref{mainthm} is an elementary lemma which, roughly speaking, asserts that if $\Gamma\curvearrowright (X,\mu)$ is a p.m.p. action with spectral gap, then 
any ``almost $\Gamma$-invariant" Borel map $\rho:X\rightarrow Y$ into a Polish space $Y$ is ``almost constant" (Lemma \ref{spectral}).
In precise terms, this means that for any $\varepsilon>0$, we can find a finite subset $F\subset\Gamma$ and $\delta>0$ such that whenever a Borel map $\rho:X\rightarrow Y$ satisfies
$\mu(\{x\in X\;|\;\rho(gx)=\rho(x)\})\geqslant 1-\delta$, for all $g\in F$, there exists $y\in Y$ with $\mu(\{x\in X\;|\;\rho(x)=y\})\geqslant 1-\varepsilon.$
In particular, for any such map $\rho$ we have that $\mu(\{x\in X|\rho(gx)=\rho(x)\})\geqslant 1-2\varepsilon$, for all $g\in\Gamma$.

 To explain how  the above ingredients are combined to prove Theorem \ref{mainthm}, assume that the translation action $\Gamma\curvearrowright (G,m_G)$ has spectral gap and let $\theta:G\rightarrow H$ be a Borel map such that $\theta(\Gamma x)\subset\Lambda\theta(x)$, for almost every $x\in G$. Consider the cocycle $w:\Gamma\times G\rightarrow\Lambda$ defined by $\theta(gx)=w(g,x)\theta(x)$.
 
 The key idea of the proof is to show that $w$ satisfies the local uniformity condition defined above. Once this is achieved, the first ingredient of proof implies that $w$ is virtually cohomologous to a homomorphism, and the conclusion of Theorem \ref{mainthm} follows easily by using standard arguments.

To prove the local uniformity condition, for every $t\in G$, we introduce a Borel map $\rho_t:G\rightarrow H$ given by $\rho_t(x)=\theta(x)^{-1}\theta(xt)$. It is clear that if  we fix $g\in\Gamma$, then $\theta_t(gx)=\theta_t(x)\Longleftrightarrow w(g,xt)=w(g,x)$. Since $w$ takes values into a countable group, it follows that $m_G(\{x\in G\;|\;\theta_t(gx)=\theta_t(x)\})\rightarrow 1$, as $t$ approaches the identity in $G$.
Since  $g\in\Gamma$ is arbitrary, the second ingredient of proof implies that we must have that $\inf_{g\in\Gamma}m_G(\{x\in G\;|\;\theta_t(gx)=\theta_t(x)\})\rightarrow 1$, as $t$ approaches the identity in $G$.
This implies that $w$ indeed satisfies the local uniformity condition.

 \subsection{Approximately trivial cocycles} As explained in the previous subsection, the proofs of our main results ultimately rely on studying cocycles. As a byproduct of this study, we obtain two results which give instances of when ``approximately trivial" cocycles are cohomologous to the trivial homomorphism. I am grateful to one of the referees for pointing out that these results are of independent interest, and suggesting that they be included in the introduction.

 \begin{mlemma}\label{I} Let $\Gamma\curvearrowright (X,\mu)$ be a strongly ergodic p.m.p. action,  $H$ be a Polish group and $w:\Gamma\times X\rightarrow H$  a cocycle.
Assume that there exists a sequence of Borel maps $\{\phi_n:X\rightarrow H\}_{n\geqslant 1}$ such that for all $g\in\Gamma$ we have that $\mu(\{x\in X|w(g,x)=\phi_n(gx)\phi_n(x)^{-1}\})\rightarrow 1$, as $n\rightarrow\infty$.

Then $w$ is cohomologous to the trivial homomorphism from $\Gamma$ to $H$, i.e. there exists a Borel map $\psi:X\rightarrow H$ such that $w(g,x)=\psi(gx)\psi(x)^{-1}$, for all $g\in\Gamma$ and almost every $x\in X$.
\end{mlemma}

For the definition of strong ergodicity, see Subsection \ref{stronger}.
If $H$ is a compact Polish group, then the proof of \cite[Proposition 2.3]{Sc81} implies the following stronger version of Lemma \ref{I}.

\begin{mlemma}\emph{\cite{Sc81}}\label{J}
 Let $\Gamma\curvearrowright (X,\mu)$ be a strongly ergodic p.m.p. action,  $H$  a compact Polish group and $w:\Gamma\times X\rightarrow H$  a cocycle.
Assume that there exists a sequence of Borel maps $\{\phi_n:X\rightarrow H\}_{n\geqslant 1}$, such that  $\lim_{n\rightarrow\infty}\phi_n(gx)\phi_n(x)^{-1}=w(g,x)$, for all $g\in\Gamma$ and almost every $x\in X$.

Then $w$ is cohomologous to the trivial homomorphism from $\Gamma$ to $H$, i.e. there exists a Borel map $\psi:X\rightarrow H$ such that $w(g,x)=\psi(gx)\psi(x)^{-1}$, for all $g\in\Gamma$ and almost every $x\in X$.

\end{mlemma}

 \subsection{Organization of the paper} Besides the introduction, this paper has nine other sections. In Sections 2 and 3, we collect a number of results that are needed in the rest of the paper, including the two main ingredients of the proof of Theorem \ref{mainthm}. In Sections 4-6, we prove several rigidity results for homomorphisms between equivalence relations associated to compact actions. In particular, Theorem \ref{mainthm} is proven in Section 4. Finally, the last four sections are devoted to the proofs of the corollaries presented in the introduction.
\subsection{Acknowledgements} I would like to thank Sorin Popa and Stefaan Vaes for many helpful comments, and the two referees for suggestions that helped improve the exposition.

\section {Preliminaries}

In this section we collect several basic notions and results that we will use throughout the paper.

Recall that a {\it standard Borel space}  is a Polish space $X$  (i.e.  a separable complete metrizable topological space) endowed with its $\sigma$-algebra of Borel subsets. 
A {\it standard probability space} $(X,\mu)$ is a Polish space $X$ endowed with a Borel probability measure $\mu$.

 Given a compact group $G$ and a closed subgroup $K<G$, we denote by $m_G$ the Haar measure of $G$ and by $m_{G/K}$ the unique $G$-invariant Borel probability measure on $G/K$.

\subsection {Countable equivalence relations}\label{basic}
We continue by recalling several notions about countable equivalence relations.

 An equivalence relation $\mathcal R$ on  a standard Borel space $X$ is  called {\it countable Borel} if 
 \begin{itemize}
 \item the equivalence class $[x]_{\mathcal R}:=\{y\in X|(x,y)\in\mathcal R\}$ is countable, for every $x\in X$, and
 \item $\mathcal R$ is a Borel subset of $X\times X$.
  \end{itemize}

If $\Gamma\curvearrowright X$ is a Borel action of a countable group $\Gamma$, then the {\it orbit equivalence relation} $$\mathcal R(\Gamma\curvearrowright X):=\{(x,y)\in X\times X|\;\Gamma x=\Gamma y\}$$ is countable Borel. Conversely, J. Feldman and C.C. Moore proved that every countable Borel equivalence relation arises in this way \cite[Theorem 1]{FM77}.

Let $\mathcal R$ and $\mathcal S$ be countable Borel equivalence relations on standard Borel spaces $X$ and $Y$. We say that  $\mathcal R$ is  {\it Borel reducible to} $\mathcal S$ and write $\mathcal R\leqslant_{B}\mathcal S$, if there exists a Borel map $\theta:X\rightarrow Y$ such that $(x,y)\in\mathcal R$ if and only if $(\theta(x),\theta(y))\in\mathcal S.$ In this case, $\theta$ is called a {\it reduction from $\mathcal R$ to $\mathcal S$}.

A Borel map $\theta:X\rightarrow Y$ is called a {\it homomorphism from $\mathcal R$ to $\mathcal S$} if $(\theta(x),\theta(y))\in\mathcal S$, for all $(x,y)\in\mathcal R$. If $X$ is endowed with a Borel probability measure $\mu$, then we say that a homomorphism $\theta:X\rightarrow Y$ from $\mathcal R$ to $\mathcal S$ is {\it trivial} if there is some $y\in Y$ such that $\theta(x)\in [y]_{\mathcal S}$, for $\mu$-almost every $x\in X$.

 A countable Borel equivalence relation $\mathcal R$ on a standard probability space $(X,\mu)$ is  called {\it probability measure preserving} (abbreviated {\it p.m.p.}) if every Borel automorphism $\theta:X\rightarrow X$ that satisfies $\theta(x)\in [x]_{\mathcal R}$, for all $x\in X$, preserves $\mu$. 
 
  The {\it full group} of $\mathcal R$, denoted $[\mathcal R]$, is the group of  automorphisms of $(X,\mu)$ such that $\theta(x)\in [x]_{\mathcal R}$, for $\mu$-almost every $x\in X$. The {\it full pseudogroup} of $\mathcal R$, denoted $[[\mathcal R]]$, is the set of $\mu$-preserving bijections $\theta:A\rightarrow B$ between Borel subsets of $X$ satisfying $\theta(x)\in [x]_{\mathcal R}$, for $\mu$-almost every $x\in A$.

 The {\it automorphism group} of $\mathcal R$, denoted Aut$(\mathcal R)$, is the group of automorphisms $\theta$ of $(X,\mu)$  such that $(x,y)\in\mathcal R$ if and only if $(\theta(x),\theta(y))\in\mathcal R$. Then $[\mathcal R]$ is a normal subgroup of Aut$(\mathcal R)$ and the {\it outer automorphism group} of $\mathcal R$ is defined as Out$(\mathcal R)=$ Aut$(\mathcal R)/[\mathcal R]$.

Let $\mathcal R$ and $\mathcal S$ be countable p.m.p. equivalence relations on  standard probability space $(X,\mu)$ and $(Y,\nu)$. We say that $\mathcal R$ and $\mathcal S$ are {\it orbit equivalent} and write $\mathcal R\sim_{\text{OE}} S$, if there exists an isomorphism of probability spaces $\theta:X\rightarrow Y$ such that $(x,y)\in\mathcal R$ if and only if $(\theta(x),\theta(y))\in\mathcal S$, almost everywhere. Moreover, we say that $\mathcal R$ and $\mathcal S$ are {\it stably orbit equivalent} if $\mathcal R_{|X_0}\sim_{\text{OE}}\mathcal S_{|Y_0}$, for some Borel subsets $X_0\subset X$ and $Y_0\subset Y$ of positive measure.  Here, we denote by   $\mathcal R_{|X_0}:=\mathcal R\cap (X_0\times X_0)$ the restriction of $\mathcal R$ to a Borel subset $X_0$ of $X$.

 Two p.m.p. actions $\Gamma\curvearrowright (X,\mu)$ and $\Lambda\curvearrowright (Y,\nu)$ are said to be {\it (stably) orbit equivalent} if the associated equivalence relations $\mathcal R(\Gamma\curvearrowright X)$ and $\mathcal R(\Lambda\curvearrowright Y)$ are (stably) orbit equivalent.

Let $\mathcal R$ be a countable ergodic p.m.p. equivalence relation on a standard probability space $(X,\mu)$.
The {\it fundamental group} of  $\mathcal R$, denoted $\mathcal F(\mathcal R)$,  is  the multiplicative group of $t>0$ for which there exist Borel sets $X_1,X_2\subset X$ of positive measure such that $\mathcal R_{|X_1}\sim_{\text{OE}}\mathcal R_{|X_2}$ and $\mu(X_1)=t\mu(X_2)$.

\subsection{Profinite and compact actions} 
A p.m.p. action $\Gamma\curvearrowright (X,\mu)$  of a countable group $\Gamma$ on a standard probability space $(X,\mu)$ is  called {\it profinite} if it is the inverse limit of a sequence of  p.m.p. actions of $\Gamma$ on finite probability spaces (see \cite[Definition 1.1]{Io08}). This   means that  we can write  $(X,\mu)=\varprojlim\; (X_n,\mu_n)$, where $X_n$ is a probability space of finite cardinality such that the subalgebra $L^{\infty}(X_n)\subset L^{\infty}(X)$ is  $\Gamma$-invariant, for all $n$.

 If such an action is ergodic, then it is isomorphic to an action of the form  $\Gamma\curvearrowright \varprojlim(\Gamma/\Gamma_n,\mu_n)$, where $\{\Gamma_n\}_n$ is a descending chain of finite index subgroups of $\Gamma$ and $\mu_n$ denotes the normalized counting measure on $\Gamma/\Gamma_n$ (see \cite[Remark 1.3]{Io08}).

If $\Gamma_n$ is a normal subgroup of $\Gamma$, for all $n$, then the  action $\Gamma\curvearrowright \varprojlim(\Gamma/\Gamma_n,\mu_n)$  can be identified with the left translation action $\Gamma\curvearrowright (G,m_G)$, where $G:=\varprojlim{\Gamma/\Gamma_n}$ is  the {\it profinite completion} of  $\Gamma$ with respect to $\{\Gamma_n\}$.  Note that the action $\Gamma\curvearrowright (G,m_G)$ is free if and only  $\cap_n\Gamma_n=\{e\}$.

Now, recall that a p.m.p. action $\Gamma\curvearrowright (X,\mu)$ is {\it compact} if the image of $\Gamma$ in the automorphism group of $(X,\mu)$ is compact. Any compact action is isomorphic to an action of the form $\Gamma\curvearrowright (G/K,m_{G/K})$, where $G$ is some compact group in which $\Gamma$ embeds densely and  $K<G$ is a closed subgroup. Note that any profinite action is compact.
More precisely, any ergodic profinite p.m.p. action $\Gamma\curvearrowright (X,\mu)$ is isomorphic to a compact action  $\Gamma\curvearrowright (G/K,m_{G/K})$, where $G$ is some profinite completion of $\Gamma$. 

{\bf Convention.} In order to differentiate the actions $\Gamma\curvearrowright G$ from the general profinite/compact actions $\Gamma\curvearrowright G/K$, we will refer to the former as {\it translation} profinite/compact actions.

Next, we give some examples of translation profinite actions that we will use later in the paper.

\begin {example}\label{profinite} Let $S$ be a set of primes. We define the profinite groups $$G_S=\prod_{p\in S}SL_2(\mathbb F_p)\;\;\;\text{and}\;\;\;K_S=\prod_{p\in S}SL_2(\mathbb Z_p),$$
where $\mathbb F_p$ is the field with $p$ elements and $\mathbb Z_p$ is the ring of $p$-adic integers. 
The Strong Approximation Theorem  gives that the diagonal embeddings of $SL_2(\mathbb Z)$ into $G_S$ and into $K_S$ are dense. 
Therefore, the left translation actions $SL_2(\mathbb Z)\curvearrowright G_S$ and $SL_2(\mathbb Z)\curvearrowright K_S$ are profinite and ergodic.

For a subgroup $\Gamma<SL_2(\mathbb Z)$, we denote by $G_{\Gamma,S}$ and $K_{\Gamma,S}$ the closure of $\Gamma$ in $G_S$ and $K_S$.
If $\Gamma$ is non-amenable, then an extension of the Strong Approximation Theorem to linear groups established in the 1980s by Nori, Weisfeiler and others guarantees that $G_{\Gamma,S}<G_S$ and $K_{\Gamma,S}<K_S$ are open subgroups (see \cite[Theorem 16.4.1]{LS03}). 
Thus, the left translation actions $\Gamma\curvearrowright G_{\Gamma,S}$ and $\Gamma\curvearrowright K_{\Gamma,S}$ are also profinite and ergodic. 

\end{example}

\subsection{Spectral gap and strong ergodicity}\label{stronger} An ergodic p.m.p. action $\Gamma\curvearrowright (X,\mu)$ is  called {\it strongly ergodic} if for any sequence  $\{A_n\}_n$ of Borel subsets of $X$ satisfying $\mu(gA_n\Delta A_n)\rightarrow 0$, for all $g\in \Gamma$, we must have that $\mu(A_n)(1-\mu(A_n))\rightarrow 0$.
The action $\Gamma\curvearrowright (X,\mu)$ is said to have {\it spectral gap} if the associated unitary representation $\pi$ of $\Gamma$ on $L^2_0(X)=L^2(X)\ominus\mathbb C1$ has spectral gap, i.e. there is no sequence $\{\xi_n\}_n$ of unit vectors in $L^2_0(X)$ satisfying $\|\pi(g)(\xi_n)-\xi_n\|_2\rightarrow 0$, for all $g\in\Gamma$.

If an action has spectral gap then it is strongly ergodic. The converse is false in general (see \cite[Theorem A.3.2]{HK05}), even within the class of profinite actions (see \cite[Theorem 5]{AE10}). 
Nevertheless,  the converse is true for translation profinite actions  \cite{AE10}. Moreover, the following holds:

\begin{proposition}\label{tau}
Let $\Gamma$ be a residually finite group and $\{\Gamma_n\}_n$ a descending chain  of finite index, normal subgroups with $\cap_n\Gamma_n=\{e\}$. Let   $G=\varprojlim\Gamma/\Gamma_n$ be the profinite completion of $\Gamma$ with respect to  $\{\Gamma_n\}_n$. Consider the left translation action $\Gamma\curvearrowright (G,m_G)$, where $m_G$ is the Haar measure of $G$.
 
 Then the following conditions are equivalent:

\begin{enumerate}
\item The action $\Gamma\curvearrowright (G,m_G)$ is strongly ergodic.

\item The action $\Gamma\curvearrowright (G,m_G)$ has spectral gap.

\item  $\Gamma$ 
has \text{\bf property $(\tau)$} with respect to the family of subgroups $\{\Gamma_n\}_n$, i.e. the representation of $\Gamma$ on $\oplus_{n}\ell^2_0(\Gamma/\Gamma_n)$ has spectral gap, where  $\ell^2_0(\Gamma/\Gamma_n)$ denotes the Hilbert space of functions $f\in\ell^2(\Gamma/\Gamma_n)$ with $\sum_{x\in \Gamma/\Gamma_n}f(x)=0$.
\end{enumerate}
Moreover, if $\Gamma$ is generated by a finite set $S$, then  conditions $(1)$-$(3)$ are also equivalent to
\begin{enumerate}
\setcounter{enumi}{3}
\item The Cayley graphs Cay$(\Gamma/\Gamma_n,S)$ form a sequence of \text{\bf expanders}.
\end{enumerate}
\end{proposition}

{\it Proof.}
The implication (1) $\Rightarrow$ (2) is clear, while (2) $\Rightarrow$ (1) was proved in \cite[Theorem 4]{AE10}. The proof that (2) $\Leftrightarrow$ (3) $\Leftrightarrow$ (4) is well-known and straightforward (see e.g. \cite{LZ03} or \cite{Lu12}). 
\hfill$\blacksquare$

Next, we collect from the literature several examples of profinite actions with spectral gap.

\begin{example}\label{selberg}   If $\Gamma$ has Kazhdan's property (T) (e.g. if $\Gamma=\text{SL}_{n}(\mathbb Z)$, for $n\geqslant 3$), then any ergodic action of $\Gamma$ has spectral gap. 
In this paper, however, we are interested in profinite actions with spectral gap of groups such as $SL_2(\mathbb Z)$ and $\mathbb F_n$ ($n\geqslant 2$), which do not have property (T).  

 One main source of examples will be  the famous theorem of Selberg asserting that $SL_2(\mathbb Z)$ has property ($\tau$) with respect to the congruence subgroups $\Gamma(m):=\ker(SL_2(\mathbb Z)\rightarrow SL_2(\mathbb Z/m\mathbb Z))$  (see\cite[Chapter 4]{LZ03}). This implies that the actions $SL_2(\mathbb Z)\curvearrowright G_S$ and $SL_2(\mathbb Z)\curvearrowright K_S$ defined in Example \ref{profinite} have spectral gap. Moreover, it gives that whenever $\Gamma<SL_2(\mathbb Z)$ is a finite index subgroup, then  the actions $\Gamma\curvearrowright G_{\Gamma,S}$ and $\Gamma\curvearrowright K_{\Gamma,S}$ have spectral gap, where  $G_{\Gamma,S}$ and $K_{\Gamma,S}$ denote the closure of $\Gamma$ in $G_S$ and $K_S$, respectively.

Recently, Selberg's theorem has been vastly generalized in a series of papers starting with J. Bourgain and A. Gamburd's breakthrough work \cite{BG05}. In particular, J. Bourgain and P. Varj\'{u} have  shown that any non-amenable subgroup  $\Gamma<SL_2(\mathbb Z)$ has property ($\tau$) with respect to the family $\{\Gamma\cap\Gamma(m)\}_{m\in\mathbb Z}$  (see  \cite[Theorem 1]{BV10}) 
Hence, the actions  $\Gamma\curvearrowright G_{\Gamma,S}$ and $\Gamma\curvearrowright K_{\Gamma,S}$ have spectral gap, whenever $\Gamma<SL_2(\mathbb Z)$ is a non-amenable subgroup. 
\end{example}

\begin{remark}
Note that if $\Gamma$ is a {\it co-amenable} subgroup of $SL_2(\mathbb Z)$, then property ($\tau$) of $\Gamma$ relative to the family  $\{\Gamma\cap\Gamma(m)\}_{m\in\mathbb Z}$ can be deduced from Selberg's theorem (see \cite{Sh99}). If we consider a finite index embedding of $\mathbb F_n$, $2\leqslant n<\infty$, into $SL_2(\mathbb Z)$, then the commutator subgroup $\Gamma=[\mathbb F_n,\mathbb F_n]\cong\mathbb F_{\infty}$ is co-amenable inside $SL_2(\mathbb Z)$. We thus have an alternative way of seeing that the actions  $\Gamma\curvearrowright G_{\Gamma,S}$ and $\Gamma\curvearrowright K_{\Gamma,S}$ have spectral gap, for this specific group $\Gamma$.
\end{remark}

We continue this section with a lemma that will be a key ingredient in the proof of Theorem \ref{mainthm}.

\begin {lemma} \label{spectral}Let $\Gamma\curvearrowright (X,\mu)$ be a strongly ergodic p.m.p. action and let $\varepsilon>0$. 

Then we can find $\delta>0$ and $F\subset\Gamma$ finite such that if a Borel map $\rho:X\rightarrow Y$ with values into a standard Borel space $Y$  satisfies $\mu(\{x\in X|\;\rho(gx)=\rho(x)\})\geqslant 1-\delta$, for all $g\in F$, then there exists $y\in Y$ such that $\mu(\{x\in X|\;\rho(x)=y\})\geqslant 1-\varepsilon$.

\end{lemma}

{\it Proof.} We may clearly assume that $\varepsilon\in (0,1)$. Then we can find $\delta>0$ and $F\subset\Gamma$ finite such that there is no Borel set $A\subset X$  satisfying $\mu(A)\in (\frac{\varepsilon}{2},1-\frac{\varepsilon}{2})$ and $\mu(g^{-1} A\Delta A)\leqslant\delta$, for all $g\in F$.

Let $\rho:X\rightarrow Y$ be as in the hypothesis and denote by $\nu=\rho_*\mu$ the push-forward of $\mu$ through $\rho$. Let $B\subset Y$ be a Borel set and put $A=\rho^{-1}(B)$. Then $g^{-1}A\Delta A\subset\{x\in X|\rho(gx)\not=\rho(x)\}$ and therefore $\mu(g^{-1}A\Delta A)\leqslant \delta$, for all $g\in F$. It follows that $\nu(B)=\mu(A)\not\in (\frac{\varepsilon}{2},1-\frac{\varepsilon}{2})$.

Thus, $\nu$ is Borel probability measure on $Y$ such that $\nu(B)\not\in(\frac{\varepsilon}{2},1-\frac{\varepsilon}{2})$, for every Borel set $B\subset Y$. This implies that
  there is $y\in Y$ such that $\nu(\{y\})\geqslant 1-\varepsilon$.
\hfill$\blacksquare$

\subsection{Smooth actions} 
Next, we recall the notion of smooth actions and prove an elementary result that we will need in Section \ref{hrig}.

\begin{definition} A Borel space $X$ is called {\it countably separated} if there exists a sequence of Borel sets which separate points. A  Borel action $H\curvearrowright X$ of a topological group $H$ on a standard Borel space $X$ is called {\it smooth} if the quotient Borel structure on $X/H$ is countably separated.  
\end{definition}

\begin{lemma}\label{smooth}

Let $H$ be a locally compact Polish group, $\Omega$ be a Polish space, and $H\curvearrowright\Omega$ be a  smooth continuous action. Denote by $\Omega_0$  the set of $x\in\Omega$ such that $hx\not=x$, for all $h\in H\setminus\{e\}$.

Then $\Omega_0$ is a $G_{\delta}$ subset of $\Omega$ and there exists a Borel map $p:\Omega_0\rightarrow H$ such that $p(hx)=hp(x)$, for all $h\in H$ and $x\in\Omega_0$. 
\end{lemma}

{\it Proof.}
We denote by $d_H$ and $d_{\Omega}$ the distance functions on $H$ and $\Omega$ which give the respective Polish topologies. Since $H$ is locally compact, the distance $d_H$ can be chosen proper, i.e. such that the closed ball $\{h\in H|d_H(h,e)\leqslant M\}$ is compact, for all $M>0$.

Let  $m,n\geqslant 1$. We denote by $U_{m,n}$ the set of $x\in\Omega$ such that $d_{\Omega}(hx,x)>\frac{1}{m}$, for all $h\in H$ satisfying $\frac{1}{n}\leqslant d_{H}(h,e)\leqslant n$. We claim that $U_{m,n}$ is an open set. Indeed, let $x_k$ be a sequence in $\Omega\setminus U_{m,n}$ which converges to a point $x\in\Omega$. Then we can find a sequence $h_k\in H$ such that $d_{\Omega}(h_kx_k,x_k)\leqslant\frac{1}{m}$ and $\frac{1}{n}\leqslant d_H(h_k,e)\leqslant n$, for all $k\geqslant 1$. Since the set of $h\in H$ which satisfy $\frac{1}{n}\leqslant d_H(h,e)\leqslant n$ is compact, after passing to a subsequence we may assume that $h_k$ converges to some $h\in H$. But then we have that  $\frac{1}{n}\leqslant d_H(h,e)\leqslant n$ and $d_{\Omega}(hx,x)=\lim\limits_{k\rightarrow\infty}d_{\Omega}(h_kx_k,x_k)\leqslant\frac{1}{m}$. Thus, $x\in\Omega\setminus U_{m,n}$. This shows that $\Omega\setminus U_{m,n}$ is closed and hence that  $U_{m,n}$ is open.
Since $\Omega_0=\cap_{n\geqslant 1}(\cup_{m\geqslant 1}U_{m,n})$, we deduce that $\Omega_0$ is a $G_{\delta}$ set.

Next, since the action $H\curvearrowright\Omega$ is smooth, it admits a {\it Borel selector}. More precisely,  we can find a Borel map $s:\Omega\rightarrow\Omega$ such that $s(x)\in Hx$ and $s(x)=s(hx)$, for all $x\in\Omega$ and every $h\in H$  (see \cite[Exercise 18.20 iii)]{Ke95}). 
Since the action $H\curvearrowright \Omega_0$ is free, for every $x\in\Omega_0$, there is a unique $p(x)\in H$ such that $x=p(x)s(x)$.
The map $p:\Omega_0\rightarrow H$ clearly satisfies $p(hx)=hp(x)$, for all $h\in H$ and $x\in \Omega_0$. 

Let us show that $p$ is Borel. To this end, let $F\subset H$ be a closed subset.  Then the map $f:\Omega\rightarrow [0,\infty)$ given by $f(x)=\inf_{h\in F}d_{\Omega}(x,hs(x))$ is Borel. 
Note that if $x\in\Omega_0$ and $p(x)\in F$, then $f(x)=0$. Conversely, let $x\in\Omega_0$ such that $f(x)=0$. We claim that $p(x)\in F$.
Indeed, there is a sequence $\{h_{n}\}_{n\geqslant 1}$ in $F$ such that $h_{n}s(x)\rightarrow x$, as $n\rightarrow\infty$. Since the action $H\curvearrowright\Omega$ is smooth and the stabilizer of $x$ is trivial, the map $H\ni h\rightarrow hx\in Hx$ is a homeomorphism \cite[Theorem 2.1.14]{Zi84}. Thus, we can find $h\in H$ such that $h_{n}\rightarrow h$, as $n\rightarrow\infty$. Since $F$ is closed, we get that $h\in F$. This implies that $hs(x)=x$ and hence $p(x)=h\in F$.

 Altogether, it follows that $\{x\in\Omega_0|p(x)\in F\}=\Omega_0\cap\{x\in\Omega|f(x)=0\}$ is a Borel set. Since this holds for any closed set $F\subset H$, we get that  $p$ is Borel. \hfill$\blacksquare$

\subsection{Extensions of homomorphisms}  We finish this section with an elementary result about extending homomorphisms from a 
dense subgroup of a Polish group to the whole group.

\begin{lemma}\label{ext}  Let $G$ be a locally compact Polish group, $m_G$ a Haar measure of $G$, and $H$ a Polish group. Let $\Gamma<G$ be  a dense subgroup and $\delta:\Gamma\rightarrow H$ a homomorphism. Assume that  $\theta:G\rightarrow H$ is a Borel map such that for all $g\in\Gamma$ we have $\theta(gx)=\delta(g)\theta(x)$, for $m_G$-almost every $x\in G$.

Then $\delta$ extends to a continuous homomorphism $\delta:G\rightarrow H$ and we can find $h\in H$ such that $\theta(g)=\delta(g)h$, for $m_G$-almost every $g\in G$. 
\end{lemma}

{\it Proof.} 
Let $f:G\times G\rightarrow H$ be given by $f(x,y)=\theta(x)^{-1}\theta(y)$. Define $S$ to be the set of $g\in G$ such that $f(gx,gy)=f(x,y)$, for $(m_G\times m_G)$-almost every $(x,y)\in G\times G$. Since $H$ is Polish, $S$ is a closed subgroup of $G$.  
Since $\Gamma<G$ is dense and $\Gamma\subset S$ by the hypothesis, we deduce that $S=G$. 

This implies that given $g\in G$, we have that $\theta(gx)\theta(x)^{-1}=\theta(gy)\theta(y)^{-1}$, for $(m_G\times m_G)$-almost every $(x,y)\in G\times G$. Thus, we can find $\delta(g)\in H$ such that $\theta(gx)\theta(x)^{-1}=\delta(g)$, for almost every $x\in G$. It is easy to see that
 $\delta:G\rightarrow H$ must be a continuous homomorphism.

Finally, since $\theta(gx)\theta(x)^{-1}=\delta(g)$, for $(m_G\times m_G)$-almost every $(g,x)\in G\times G$, we can find $x_0\in G$ such that $\theta(gx_0)\theta(x_0)^{-1}=\delta(g)$, for $m_G$-almost every $g\in G$. Hence, for almost every $g\in G$ we have that $\theta(g)=\delta(gx_0^{-1})\theta(x_0)=\delta(g)(\delta(x_0^{-1})\theta(x_0))$. Thus, the conclusion holds for $h=\delta(x_0^{-1})\theta(x_0)$.
\hfill$\blacksquare$

\section{Cocycle rigidity} The purpose of this section is twofold. Firstly, we extract from \cite{Io08} and \cite{Fu09} two rigidity results for cocycles associated with profinite and compact actions. There results will be essential in the proofs of our main results. Secondly, we prove that for strongly ergodic actions, any cocycle that is approximately cohomologous to the trivial cocycle must be cohomologous to the trivial cocycle. We start by recalling some notation.

  Let $\Gamma\curvearrowright (X,\mu)$ be a p.m.p. action of a countable group $\Gamma$ on a standard probability space $(X,\mu)$. 
A measurable map $w:\Gamma\times X\rightarrow\Lambda$ into a group $\Lambda$ is called a {\it measurable cocycle} (or just a {\it cocycle}) if it satisfies the identity $w(gh,x)=w(g,hx)w(h,x)$, for all $g,h\in\Gamma$ and almost every $x\in X$. Two cocycles $w_1,w_2:\Gamma\times X\rightarrow\Lambda$ are said to be {\it cohomologous} if there exists a Borel map $\phi:X\rightarrow\Lambda$ such that $w_1(g,x)=\phi(gx)w_2(g,x)\phi(x)^{-1}$, for all $g\in\Gamma$ and almost every $x\in X$.

\subsection{Cocycle rigidity for profinite actions}
Now, assume that $\Gamma$ has property (T) and let $\Gamma\curvearrowright (X,\mu)=\varprojlim\; (X_n,\mu_n)$ be a free ergodic profinite p.m.p. action.  We proved that any cocycle $w:\Gamma\times X\rightarrow\Lambda$ with values into a countable group $\Lambda$ is cohomologous to a cocycle which factors through the map $\Gamma\times X\rightarrow\Gamma\times X_n$, for some $n$ (see  \cite[Theorem B]{Io08}). 
In \cite{Fu09}, A. Furman  extended this result  from profinite to compact actions (see \cite[Theorem 5.21]{Fu09}). 

A main ingredient in the proof of Theorem \ref{mainthm} is the following criterion for untwisting cocycles for profinite actions $\Gamma\curvearrowright (X,\mu)$. This criterion is an easy consequence of the proof of \cite[Theorem B]{Io08} and is implicitly proved in the proof of  \cite[Theorem 5.21]{Fu09}. 
A crucial aspect of this criterion is that it applies to arbitrary residually finite groups $\Gamma$, that are not assumed to have property (T).

\begin{theorem} \emph{\cite{Io08} \cite {Fu09}}\label{cocycle}
Let $\Gamma$ be a residually finite group and $\{\Gamma_n\}_n$ be a descending chain  of finite index, normal subgroups with trivial intersection, $\cap_n\Gamma_n=\{e\}$. Let   $G=\varprojlim\Gamma/\Gamma_n$ be the profinite completion of $\Gamma$ with respect to  $\{\Gamma_n\}_n$ and consider the left translation action $\Gamma\curvearrowright (G,m_G)$.
 
 Let $\Lambda$ be a countable group and $w:\Gamma\times G\rightarrow \Lambda$ be a cocycle. Assume that for some constant $C\in (\frac{31}{32},1)$ we can find a neighborhood $V$ of the identity $e$ in $G$ such that  
\begin{equation}\label{uniform} m_G(\{x\in G|\;w(g,xt)=w(g,x)\})\geqslant C, \;\;\;\text{for all} \;\;\; g\in\Gamma\;\;\;\text{and every}\;\;\;t\in V.
\end{equation}

Then we can find an open subgroup $G_0<G$ such that  the restriction of $w$ to $(\Gamma\cap G_0)\times G_0$ is cohomologous to  a homomorphism $\delta:\Gamma\cap G_0\rightarrow\Lambda$.
\end{theorem}

 For the reader's convenience, we give two proofs of Theorem \ref{cocycle}, following \cite{Io08} and \cite{Fu09}, respectively. 

\vskip 0.05in
{\it First proof of Theorem \ref{cocycle}.}  Denote by $c$ the counting measure on $\Lambda$.
Following \cite{Io08},   consider the infinite measure preserving action of $\Gamma$ on $(Z,\rho):=(G\times G\times\Lambda,m_G\times m_G\times c)$ given by $$g\cdot(x,y,\lambda)=(gx,gy,w(g,x)\lambda w(g,y)^{-1}).$$ Denote by $\pi:\Gamma\rightarrow \mathcal U(L^2(Z,\rho))$ the associated Koopman representation.
For $n\geqslant 0$, we let $r_n:G\rightarrow\Gamma/\Gamma_n$ be the quotient homomorphism and $\zeta_n$ be the characteristic function of the subset $\{(x,y)\in G\times G|r_n(x)=r_n(y)\}\times\{e\}$ of  $Z$.

Since $V$ is a neighborhood of the identity, we can find $n\geqslant 0$ such that $V$ contains $G_n:=\ker(r_n)$.
We claim that  $\|\pi(g)(\zeta_n)-\zeta_n\|_2\leqslant \sqrt{2-2C}\;\|\zeta_n\|_2$, for all $g\in\Gamma$. Towards this, let $g\in\Gamma$.
Since $r_n(x)=r_n(y)$ if and only if $x^{-1}y\in G_n$, a simple computation shows that  $$\langle\pi(g)(\zeta_n),\zeta_n\rangle=$$ $$(m_G\times m_G)(\{(x,y)\in G\times G|\;r_n(x)=r_n(y)\;\text{and}\; w(g,x)=w(g,y)\}=$$ $$(m_G\times m_G)(\{x,t)\in G\times G|\;t\in G_n\;\text{and}\;w(g,x)=w(g,xt)\}.$$

Since $m_G(\{x\in G|w(g,xt)=w(g,x)\})\geqslant C$, for all $t\in G_n$,  and $m_G(G_n)=\|\zeta_n\|_2^2=|\Gamma/\Gamma_n|^{-1}$ we deduce that $\langle \pi(g)\zeta_n,\zeta_n\rangle\geqslant C\mu(G_n)=C\|\zeta_n\|_2^2$. This  implies the claim.

Now, define $\xi_n:=\sqrt{|\Gamma/\Gamma_n|}\;\zeta_n\in L^2(Z,\rho)$. Then $\xi_n$ is a unit vector and the above claim says that $\|\pi(g)(\xi_n)-\xi_n\|_2\leqslant\sqrt{2-2C}$, for all $g\in\Gamma$.
By using a standard averaging argument, we can find a $\pi(\Gamma)$-invariant vector $\eta\in L^2(Z,\rho)$ such that $\|\eta-\xi_n\|_2\leqslant\sqrt{2-2C}$. 

Since $\sqrt{2-2C}<\frac{1}{4},$ Part 2 of the proof of \cite[Theorem B]{Io08} is satisfied.  Continuing as in the proof of \cite[Theorem B]{Io08}  gives that $w$ is cohomologous  to a cocycle which factors through the map $\Gamma\times G\rightarrow \Gamma\times\Gamma/\Gamma_N$, for some $N\geqslant n$. Then $G_0:=\ker(r_N)$ clearly satisfies the conclusion. \hfill$\blacksquare$

\vskip 0.05in
{\it Second proof of Theorem \ref{cocycle}.} Let us  now explain how the proof of  \cite[Theorem 5.21]{Fu09} also implies  Theorem \ref{cocycle}. 
Following \cite{Fu09}, fix $t\in V$ and   define a new cocycle $w_t:\Gamma\times G\rightarrow\Lambda$ by letting $w_t(g,x)=w(g,xt)$.
Since   $\Gamma\curvearrowright (G,m_G)$ is ergodic and $C>\frac{7}{8}$,  \cite [Lemma 2.1]{Io08} implies that $w_t$ is cohomologous to $w$. 
Hence, we can find a Borel map $\phi_t:G\rightarrow\Lambda$ such that $w_t(g,x)=\phi_t(gx)w(g,x)\phi_t(x)^{-1}$, for all $g\in\Gamma$ and almost every $x\in G$. 

Moreover, a close inspection of the proof of \cite [Lemma 2.1]{Io08} shows that $\phi_t$ verifies the following: there exists $\eta_t\in L^2(G\times\Lambda,m_G\times c)$
such that  $\phi_t(x)$ is the unique  $\lambda\in\Lambda$ satisfying 
$|\eta_t(x,\lambda)|>\frac{1}{2}$, for almost every $x\in G$, and $\|\eta_t-1_{G\times\{e\}}\|_2\leqslant\sqrt{2-2C}<\frac{1}{4}$. Since $$m_G(\{x\in G|\;|\eta_t(x,e)|\leqslant\frac{1}{2}\})\leqslant 4\int_{G}|\eta_t(x,e)-1|^2\;\text{d}m_G(x)\leqslant 4\|\eta_t-1_{G\times\{e\}}\|_2^2<\frac{1}{4},$$
we conclude that $m_G(\{x\in G|\phi_t(x)=e\})>\frac{3}{4}$. The proof of \cite[Theorem 5.21]{Fu09} now applies verbatim to give the conclusion. \hfill$\blacksquare$

\subsection{Cocycle rigidity for compact actions}
Later on, we will also need the following variant of Theorem \ref{cocycle} for compact actions. This results is an immediate consequence of \cite{Fu09}.

\begin{theorem}  \emph{\cite{Fu09}}\label{compact} Let $\Gamma$ be a countable group together with a dense embedding $\tau:\Gamma\hookrightarrow G$ into  a connected, compact group $G$. Consider the left translation action $\Gamma\curvearrowright (G,m_G)$, where $m_G$ is the Haar measure of $G$. Assume that  $\pi_1(G)$, the fundamental group of $G$, is finite.

Let $\Lambda$ be a countable group and $w:\Gamma\times G\rightarrow\Lambda$ be a  cocycle.  Assume that for some constant $C\in (\frac{31}{32},1)$ we can find a neighborhood $V$ of the identity $e$ in $G$ such that condition \ref{uniform} holds.

We have the following:
\begin{enumerate}
\item If any homomorphism $\pi_1(G)\rightarrow\Lambda$ is trivial (e.g. if $\pi_1(G)=\{e\}$ or if $\Lambda$ is torsion free), then $w$ is cohomologous to a homomorphism $\delta:\Gamma\rightarrow \Lambda$.
\item In general, we can find a subgroup $\Lambda_0<\Lambda$, a finite subgroup $\Lambda_1<Z(\Lambda_0)$, a Borel map $\phi:G\rightarrow\Lambda$, and a homomorphism $\delta:\Gamma\rightarrow\Lambda_0/\Lambda_1$ such that if $p:\Lambda_0\rightarrow\Lambda_0/\Lambda_1$ denotes the quotient homomorphism, then $w'(g,x)=\phi(gx)^{-1}w(g,x)\phi(x)\in\Lambda_0$ and $p(w'(g,x))=\delta(g)$, for all $g\in\Gamma$ and almost every $x\in G$.
\end{enumerate}
\end{theorem}

Here and after we denote by $Z(\Lambda)$ the {\it center} of a group $\Lambda$.

{\it Proof.} (1) Assume that the only homomorphism $\pi_1(G)\rightarrow\Lambda$ is the trivial one.
As in the second proof of Theorem \ref{cocycle}, for every $t\in V$ we  can find a Borel map $\phi_t:G\rightarrow\Lambda$  which satisfies $m_G(\{x\in G|\phi_t(x)=e\})>\frac{3}{4}$ and $w_t(g,x)=\phi_t(gx)w(g,x)\phi_t(x)^{-1}$, for all $g\in\Gamma$ and almost every $x\in G$. The conclusion then follows from the proof of \cite[Theorem 5.21]{Fu09}.

(2) Let $\tilde G$ be the universal covering group of $G$, and denote by $\pi:\tilde G\rightarrow G$ the covering map. Let $\tilde\Gamma=\pi^{-1}(\Gamma)$ and define a cocycle $\tilde w:\tilde\Gamma\times\tilde G\rightarrow \Lambda$ by letting $\tilde w(g,x)=w(\pi(g),\pi(x))$. Since $\tilde G$ is connected, it follows that $\tilde\Gamma<\tilde G$ is dense. Since $\pi_1(\tilde G)=\{e\}$,  the first part of the proof gives that $\tilde w$ is cohomologous to a homomorphism $\delta:\tilde\Gamma\rightarrow\Lambda$. Let $\psi:\tilde G\rightarrow\Lambda$ be a Borel map satisfying $\tilde w(g,x)=\psi(gx)\delta(g)\psi(x)^{-1}$, for all $g\in\tilde\Gamma$ and almost every $x\in\tilde G$.

Fix $k\in\ker{\pi}$ and define $\rho_k:\tilde G\rightarrow\Lambda$ by letting $\rho_k(x)=\psi(x)^{-1}\psi(xk)$. We claim that $\rho_k$ is  constant. Note  that $\rho_k(gx)=\delta(g)\rho_k(x)\delta(g)^{-1}$, for all 
$g\in\tilde\Gamma$ and almost every $x\in\tilde G$.
Let $h\in\Lambda$  such that $A_h=\{x\in \tilde G|\rho_k(x)=h\}$ has positive measure. Then $gA_h=A_{\delta(g)h\delta(g)^{-1}}$, for all $g\in\tilde\Gamma$. This implies that for any $g_1,g_2\in\tilde\Gamma$, the sets $g_1A_h$ and $g_2A_h$ are either disjoint or equal. Hence, $A_h$  is invariant under some finite index subgroup $\Gamma_0<\tilde\Gamma$. 
Since $\tilde G$ is connected and $\tilde\Gamma<\tilde G$ is dense, we get that $\Gamma_0<\tilde G$ is also dense. As a consequence,  $A_h=\tilde G$, almost everywhere. 
This proves the claim.

Thus, there is a homomorphism $\rho:\ker{\pi}\rightarrow\Lambda$ such that $\psi(x)^{-1}\psi(xk)=\rho(k)$,  for all $k\in\ker{\pi}$ and almost every $x\in\tilde G$. Moreover,  $\rho(k)$ commutes with $\delta(\tilde\Gamma)$, for all $k\in\ker{\pi}$.
Since $\tilde w(kg,x)=\tilde w(g,x)$, we also get that $\psi(x)^{-1}\psi(kx)=\delta(k)^{-1}$, for all $k\in\ker{\pi}$ and almost every $x\in G$. By combining these two facts, we derive that $\delta(k)^{-1}=\rho(x^{-1}kx)$,  which implies that $\delta(k)^{-1}=\rho(k)$, for all $k\in\ker{\pi}$. Thus, if $\Lambda_0:=\delta(\tilde\Gamma)$ and $\Lambda_1:=\delta(\ker{\pi})=\rho(\ker{\pi})$, then $\Lambda_1$ is finite and $\Lambda_1<Z(\Lambda_0)$.
In particular, we have a homomorphism $\bar{\delta}:\Gamma=\tilde\Gamma/\ker{\pi}\rightarrow\Lambda_0/\Lambda_1$ given by $\bar{\delta}(x \ker{\pi})=\delta(x)\Lambda_1$.

Finally, choose a Borel map $\phi:G\rightarrow\Lambda$ such that $\phi\circ\pi=\psi$.  Note that if $y,y'\in\tilde G$ are such that $\pi(y)=\pi(y')$, then $y^{-1}y'\in\ker{\pi}$, hence $\psi(y')=\psi(y)\rho(y^{-1}y')\in\psi(y)\Lambda_1$. It follows  that if $g\in\Gamma$ and $\tilde g\in\delta^{-1}(\{g\})$, then $w'(g,x):=\phi(gx)^{-1}w(g,x)\phi(x)\in\delta(\tilde{g})\Lambda_1\subset\Lambda_0$, for almost every $x\in G$.  This shows that $p(w'(g,x))=\delta(\tilde g)\Lambda_1=\bar{\delta}(g)$ and
finishes the proof.
\hfill$\blacksquare$

\subsection{Approximately trivial cocycles}
We end this section by proving Lemmas \ref{I} and \ref{J}.

\subsection{Proof of Lemma \ref{I}} Let $n\geqslant 1$.
Since the action $\Gamma\curvearrowright (X,\mu)$ is strongly ergodic, by applying Lemma \ref{spectral} we can find $\delta_n>0$ and $F_n\subset\Gamma$ finite such that if a Borel map $\rho:X\rightarrow H$ satisfies $\mu(\{x\in X|\rho(gx)=\rho(x)\})\geqslant 1-\delta_n$, for every $g\in F_n$, then there exists $y\in H$ such that we have $\mu(\{x\in X|\rho(x)=y\})\geqslant 1-\frac{1}{2^n}$. Moreover, we may clearly assume that $\delta_{n+1}\leqslant\delta_n\leqslant\frac{1}{2^n}$ and $F_n\subset F_{n+1}$, for every $n\geqslant 1$.

Now, after replacing the sequence $\{\phi_n:X\rightarrow H\}_{n\geqslant 1}$ with a subsequence, we may assume that 
\begin{equation}\label{alho}\mu(\{x\in X|w(g,x)=\phi_n(gx)\phi_n(x)^{-1}\})\geqslant 1-\frac{\delta_n}{2}\geqslant 1-\frac{1}{2^{n+1}},\;\;\text{for all $g\in F_n$ and  $n\geqslant 1$}.\end{equation}

Next, we claim that there exists a sequence of Borel maps  $\{\psi_n:X\rightarrow H\}_{n\geqslant 1}$ such that

\begin{enumerate}
\item $\mu(\{x\in X|w(g,x)=\psi_n(gx)\psi_n(x)^{-1}\})\geqslant 1-\frac{\delta_n}{2},$ for all $g\in F_n$ and every $n\geqslant 1$, and

\item $\mu(\{x\in X|\psi_n(x)=\psi_{n+1}(x)\})\geqslant 1-\frac{1}{2^n},$ for every $n\geqslant 1$.
\end{enumerate}

We proceed by induction. Thus, let $\psi_1:=\phi_1$. Assume that we have constructed $\psi_1,...,\psi_n$ and let us construct $\psi_{n+1}$. 
To this end, we define $\rho:X\rightarrow H$ by letting $\rho(x)=\psi_n(x)^{-1}\phi_{n+1}(x)$. 

If $g\in F_n$, then  by combining equation \ref{alho} with condition (1), we derive that the set of $x\in X$ such that $w(g,x)=\psi_n(gx)\psi_n(x)^{-1}=\phi_{n+1}(gx)\phi_{n+1}(x)^{-1}$ has measure $\geqslant 1-\delta_n$.
From this it follows that $\mu(\{x\in X|\rho(gx)=\rho(x)\})\geqslant 1-\delta_n$, for all $g\in F_n$.  The above implies that we can find $h\in H$ such that $\mu(\{x\in X|\rho(x)=h\})\geqslant 1-\frac{1}{2^n}$. 

Therefore, if we define $\psi_{n+1}(x):=\phi_{n+1}(x)h^{-1}$, then $\mu(\{x\in X|\psi_n(x)=\psi_{n+1}(x)\})\geqslant 1-\frac{1}{2^n}$ and also $\mu(\{x\in X|w(g,x)=\psi_{n+1}(gx)\psi_{n+1}(x)^{-1}\})=\mu(\{x\in X|w(g,x)=\phi_{n+1}(gx)\phi_{n+1}(x)^{-1}\})\geqslant 1-\delta_n$, for all $g\in F_n$. This finishes the proof of the claim.

For  $N\geqslant 1$, let $X_N:=\{x\in X|\psi_n(x)=\psi_N(x),\;\;\text{for all $n\geqslant N$}\}$. Then condition (2) implies that $\mu(X_N)\geqslant 1-\frac{1}{2^{N-1}}$.
Thus, if $X'=\cup_{N\geqslant 1}X_N$, then $\mu(X')=1$.
We define $\psi:X\rightarrow H$ by letting $$\psi(x)=\begin{cases}\psi_N(x),\;\;\text{if $x\in X_N$, for some $N\geqslant 1$, and}\\ e,\;\;\text{if $x\in X\setminus X'$}. \end{cases}$$

We claim that $\psi$ satisfies the conclusion of the lemma. To see this, fix $g\in\Gamma$. For every $N\geqslant 1$, we define  $Y_N:=\{x\in X|w(g,x)=\psi_n(gx)\psi_n(x),\;\text{for all $n\geqslant N$}\}$. Then equation \ref{alho} gives that $\mu(Y_N)\geqslant 1-\frac{1}{2^{N+1}}$. Hence,  $Y'=\cup_{N\geqslant 1}Y_N$ verifies $\mu(Y')=1$. It is now clear that for every $x\in g^{-1}X'\cap X'\cap Y'$ we have  $w(g,x)=\psi(gx)\psi(x)^{-1}$. Since $\mu(X')=\mu(Y')=1$, we are done.
\hfill$\blacksquare$

\subsection{Proof of Lemma \ref{J}} Let $d:H\times H\rightarrow [0,+\infty)$ be a left-right invariant metric on $H$. 
Define $\phi_{m,n}:X\rightarrow H$ by  $\phi_{m,n}(x)=\phi_m^{-1}(x)\phi_n(x)$, for $m,n\geqslant 1$. Then $\lim_{m,n\rightarrow\infty}d(\phi_{m,n}(gx),\phi_{m,n}(x))=0$, for all $g\in\Gamma$ and almost every $x\in X$. By using the strong ergodicity of the action $\Gamma\curvearrowright (X,\mu)$ as in the proof of \cite[Proposition 2.3]{Sc81} we can find $h_{m,n}\in H$ such that $\lim_{m,n\rightarrow\infty}d(\phi_{m,n}(x),h_{m,n})=0$, for almost every $x\in X$. This further implies that we can find $h_n\in H$ such that the limit $\psi(x)=\lim_{n\rightarrow\infty}\phi_n(x)h_n$ exists, for almost every $x\in X$.
But then $\psi:X\rightarrow H$ is a Borel map such that $w(g,x)=\psi(gx)\psi(x)^{-1}$, for all $g\in\Gamma$ and almost every $x\in X$.
\hfill$\blacksquare$

\section {Homomorphism rigidity and proof of Theorem \ref{mainthm}}\label{section3}

\subsection{Homomorphism rigidity for profinite actions}  The main goal of this section is to prove Theorem \ref{mainthm}.  The proof of Theorem \ref{mainthm} relies on the following rigidity result for homomorphisms between equivalence relations arising from translation profinite actions with spectral gap.

\begin{theorem}\label{maintech} 
Let $\Gamma$ be a residually finite group. Let   $G=\varprojlim\Gamma/\Gamma_n$ be the profinite completion of $\Gamma$ with respect to  a chain  of finite index normal subgroups with trivial intersection. Assume that the left translation action $\Gamma\curvearrowright (G,m_G)$ has spectral gap.

Let $\Lambda$ be a countable subgroup of a Polish group $H$ and consider the left translation  action $\Lambda\curvearrowright H$. 
Let $\theta:G\rightarrow H$ be a Borel map such that $\theta(\Gamma x)\subset\Lambda\theta(x)$, for almost every $x\in G$. Let $w:\Gamma\times G\rightarrow\Lambda$ be the cocycle defined by the formula $\theta(gx)=w(g,x)\theta(x)$, for all $g\in\Gamma$ and almost every $x\in G$.

Then we can find an open subgroup $G_0<G$,   a continuous homomorphism $\delta:G_0\rightarrow H$, a Borel map $\phi:G_0\rightarrow\Lambda$ and $h\in H$ such that

\begin{itemize}
\item $\delta(\Gamma\cap G_0)\subset\Lambda$,
\item $w(g,x)=\phi(gx)\delta(g)\phi(x)^{-1}$, for all $g\in\Gamma\cap G_0$ and almost every $x\in G_0$, and
\item $\theta(x)=\phi(x)\delta(x)h$, for almost every $x\in G_0$.
 \end{itemize}

\end{theorem}

\begin{remark} There are two useful ways of interpreting the conclusion of Theorem \ref{maintech}.

Firstly, Theorem \ref{maintech} describes all homomorphisms between $\mathcal R(\Gamma\curvearrowright G)$ and $\mathcal R(\Lambda\curvearrowright H)$.
Thus, it shows that any Borel map $\theta:G\rightarrow H$ which satisfies $\theta(\Gamma x)\subset\Lambda\theta(x)$, for almost every $x\in G$,  arises from a homomorphism $\delta:G_0\rightarrow H$, for some open subgroup $G_0<G$.

Secondly, Theorem \ref{maintech} can be viewed as a rigidity result for cocycles $w:\Gamma\times G\rightarrow\Lambda$.  
More precisely, assume that if we view $w$ as a cocycle with values in $H$, then $w$ is cohomologous to the trivial cocycle. This means that there exists a Borel map $\theta:G\rightarrow H$ such that $w(g,x)=\theta(gx)\theta(x)^{-1}$, for all $g\in\Gamma$ and almost every $x\in G$. Theorem \ref{maintech} then shows that there exists an open subgroup $G_0<G$ such that the restriction of $w$ to $(\Gamma\cap G_0)\times G_0$ is cohomologous to a homomorphism $\delta:\Gamma\cap G_0\rightarrow\Lambda$.
\end{remark}

{\it Proof.}
 We claim that can find an open subgroup $G_0<G$ and a Borel map $\phi:G_0\rightarrow\Lambda$ such that  $w(g,x)=\phi(gx)\delta(g)\phi(x)^{-1},$ for all $g\in\Gamma\cap G_0$ and almost every $x\in G_0$. 
To this end, let $\varepsilon\in (0,\frac{1}{64})$.

 Since the action $\Gamma\curvearrowright (G,m_G)$ has spectral gap, it is strongly ergodic. Lemma \ref{spectral} yields $\delta>0$ and a finite set $F\subset\Gamma$ such that if a Borel map $\rho:G\rightarrow H$   satisfies $m_G(\{x\in X|\rho(gx)=\rho(x)\})\geqslant 1-\delta$, for all $g\in F$, then  $m_G(\{x\in X|\rho(x)=y\})\geqslant 1-\varepsilon$, for some $y\in H$.

Now, if $A\subset G$ is a Borel set, then $\lim\limits_{t\rightarrow e}m_G(At\Delta A)=0$. Since $\Lambda$ is countable, it follows that
 $\lim\limits_{t\rightarrow e}m_G(\{x\in G|\;w(g,xt)=w(g,x)\})=1$. Thus,
we can find a neighborhood $V$ of $e$ in $G$ such that 
\begin{equation}\label{const}m_G(\{x\in G|\;w(g,xt)=w(g,x)\})\geqslant 1-\delta, \;\;\text{for all}\;\; g\in F\;\;\text{and}\;\;t\in V.
\end{equation}

Fix $t\in V$ and define $\rho_t:G\rightarrow H$ by $\rho_t(x)=\theta(x)^{-1}\theta(xt)$. Then for almost every $x\in G$ we have that $\rho_t(gx)=\theta(gx)^{-1}\theta(gxt)=\theta(x)^{-1}w(g,x)^{-1}w(g,xt)\theta(xt).$ Thus, $\rho_t(gx)=\rho_t(x)$ if and only if $w(g,xt)=w(g,x)$. Equation \ref{const}  gives that $m_G(\{x\in G|\;\rho_t(gx)=\rho_t(x)\})\geqslant 1-\delta$, for all $g\in F$. 

We deduce that there is $y_t\in H$ such that $m_G(\{x\in G|\;\rho_t(x)=y_t\})\geqslant 1-\varepsilon$. This implies that  $m_G(\{x\in G|\;\rho_t(gx)=\rho_t(x)\})\geqslant 1-2\varepsilon$, for all $g\in\Gamma$. Equivalently,  $$m_G(\{x\in G|\;w(g,xt)=w(g,x)\})\geqslant 1-2\varepsilon,\;\; \text{for all}\;\; g\in\Gamma.$$

 Since $1-2\varepsilon>\frac{31}{32}$ and $t\in V$ is arbitrary,  Theorem \ref{cocycle}  implies the claim. 

Denoting $\tilde\theta(x)=\phi(x)^{-1}\theta(x)$ and using the claim we have that \begin{equation}\label{hom}\tilde\theta(gx)=\delta(g)\tilde\theta(x),\;\;\text{for all}\;\; g\in\Gamma\cap G_0\;\;\text{and}\;\; \text{almost every}\;\; x\in G_0.\end{equation}

By applying Lemma \ref{ext} we get that $\delta$ extends to a continuous homomorphism $\delta:G_0\rightarrow H$ and that we can find $h\in H$ such that $\tilde\theta(g)=\delta(g)h$, for almost every $g\in G_0$.
\hfill$\blacksquare$

\subsection{Proof of part (1) of Theorem A}  In order to prove Theorem \ref{mainthm}, we handle parts (1) and (2) separately.
First, we use Theorem \ref{maintech}  to describe the stable orbit equivalences between translation profinite actions with spectral gap. More generally, we have:

\begin{corollary}\label{cor1} Let $\Gamma,\Lambda$ be residually finite groups.  Let $G=\varprojlim\Gamma/\Gamma_n$ and $H=\varprojlim\Lambda/\Lambda_n$ be the profinite completions of  $\Gamma$ and $\Lambda$ with respect to chains of finite index normal subgroups with trivial intersection. Assume that the left translation action $\Gamma\curvearrowright (G,m_G)$ has spectral gap.

Let $A\subset G$ and $B\subset H$ be Borel sets of positive  measure endowed with the probability measures obtained by restricting and rescaling $m_G$ and $m_H$. 
Let $\theta:A\rightarrow B$ be an isomorphism of probability spaces such that  $\theta(\Gamma x\cap A)\subset\Lambda\theta(x)\cap B$, for almost every $x\in A$. 

Then we can find $\tau\in [\mathcal R(\Gamma\curvearrowright G)]$, $\rho\in [\mathcal R(\Lambda\curvearrowright H)]$, open subgroups $G_0<G$ and $H_0<H$, a continuous isomorphism $\delta:G_0\rightarrow H_0$ and  $h\in H$  such that $\tau(G_0)\subset A$, $\frac{\mu(A)}{\nu(B)}=\frac{[H:H_0]}{[G:G_0]}$, 
$$\delta(\Gamma\cap G_0)\subset \Lambda\cap H_0\;\;\;\text{and}\;\;\;(\rho\circ\theta\circ\tau)(x)=\delta(x)h,\;\;\;\text{for almost every}\;\;\;x\in G_0.$$

If $\theta$ moreover satisfies $\theta(\Gamma x\cap A)=\Lambda\theta(x)\cap B$, for almost every $x\in A$, then  $\delta(\Gamma\cap G_0)=\Lambda\cap H_0$.
\end{corollary}
 
Corollary \ref{cor1} clearly  implies the {\it only if} assertion from part (1) of Theorem \ref{mainthm}. Note that the {\it if} assertion from part (1) of Theorem \ref{mainthm} is obvious because $\mathcal R(\Gamma\curvearrowright G)_{|G_0}=\mathcal R((\Gamma\cap G_0)\curvearrowright G_0)$, whenever $G_0<G$ is an open subgroup.
 
{\it Proof.} After replacing $\theta$ with $\theta\circ\tau$, for some $\tau\in[\mathcal R(\Gamma\curvearrowright G)]$, we may assume that $A$ contains an open subgroup $G_1<G$.  We will prove that the conclusion holds to $\tau=\text{id}_G$. 

Note that  the action $\Gamma\cap G_1\curvearrowright (G_1,m_{G_1})$  has spectral gap. Indeed, since the action $\Gamma\curvearrowright (G,m_G)$ has spectral gap, it follows that $\mathcal R(\Gamma\curvearrowright G)$ and hence $\mathcal R(\Gamma\curvearrowright G)_{|G_1}=\mathcal R((\Gamma\cap G_1)\curvearrowright G_1)$ are strongly ergodic. Thus, the action $\Gamma\cap G_1\curvearrowright (G_1,m_{G_1})$ is strongly ergodic and so it must have spectral gap by Proposition \ref{tau}.

 By applying Theorem \ref{maintech} to $\theta_{|G_1}$ we can find an open subgroup $G_0<G_1$,   a continuous homomorphism $\delta:G_0\rightarrow H$, a Borel map $\phi:G_0\rightarrow\Lambda$ and $h\in H$ such that $\delta(\Gamma\cap G_0)\subset\Lambda\cap H_0$ and $\theta(x)=\phi(x)\delta(x)h$, for almost every $x\in G_0.$

 Next, we claim that $K:=\ker(\delta)$ is finite. Assume that is not the case. Towards  a contradiction, let $\lambda\in\Lambda$ such that $C:=\{x\in G_0|\phi(x)=\lambda\}$ satisfies $m_G(C)>0$. Since $K$ is assumed infinite, we can find $g\in K\setminus\{e\}$ such that $m_G(g^{-1}C\cap C)>0$. Since $\theta(gx)=\theta(x)$,  for almost every $x\in g^{-1}C\cap C$, we get a contradiction with the fact that $\theta$ is 1-1.

Thus, after replacing $G_0$ with a smaller open subgroup of $G$, we may assume that $\delta$ is 1-1.
Hence, if we let $H_0:=\delta(G_0)$, then $\delta:G_0\rightarrow H_0$ is a continuous isomorphism. Note that $H_0$ is an open subgroup of $H$. 
Indeed, since $m_H(\theta(G_0))>0$ and $\theta(G_0)\subset\cup_{\lambda\in\Lambda}\lambda H_0h$, we get that $m_H(H_0)>0$. Therefore,
$[H:H_0]<\infty$ and since $H_0<H$ is a closed subgroup (being the continuous image of a compact group) it must be open.

Now, since the map $G_0\ni x\rightarrow\delta(x)h\in H_0h$ is 1-1 and $\delta(x)h\in\Lambda\theta(x)$, for almost every $x\in G_0$, we can find $\rho\in [\mathcal R(\Lambda\curvearrowright H)]$ such that \begin{equation}\label{rho}\rho(\theta(x))=\delta(x)h,\;\;\text{for almost every $x\in G_0$}.\end{equation} In particular, we have that $\nu(\theta(G_0))=\nu(H_0h)=\nu(H_0)=[H:H_0]^{-1}$. Since $\theta:A\rightarrow B$ is an isomorphism of probability spaces, we have that $\nu(\theta(G_0))=\frac{\nu(B)}{\mu(A)}\mu(G_0)=\frac{\nu(B)}{\mu(A)}[G:G_0]^{-1}$. By combining the last two identities we get that $\frac{\mu(A)}{\nu(B)}=\frac{[H:H_0]}{[G:G_0]}$. This finishes the proof of the first assertion.

For the moreover assertion, assume that $\theta(\Gamma x\cap A)=\Lambda\theta(x)\cap B$, for almost every $x\in A$. 
This implies that $\theta(\Gamma x\cap G_0)=\Lambda\theta(x)\cap\theta(G_0)$, for almost every $x\in G_0$. By applying $\rho$ to this identity and using equation \ref{rho} we get that  
$\delta(\Gamma x\cap G_0)=\Lambda\delta(x)\cap H_0,$ for almost every $x\in G_0$. 
In particular,  if $\lambda\in\Lambda\cap H_0$, then $\lambda\delta(x)\in \delta(\Gamma x\cap G_0)=\delta((\Gamma\cap G_0)x),$ for some $x\in G_0$.
Thus, $\lambda\in\delta(\Gamma\cap G_0)$ and therefore $\Lambda\cap H_0\subset\delta(\Gamma\cap G_0)$. Since the other inclusion also holds, we are done.
\hfill$\blacksquare$

\subsection{Proof of part (2) of Theorem \ref{mainthm}}
Next, we use Theorem \ref{maintech} to prove the second part of Theorem \ref{mainthm}. Moreover, we give an explicit characterization of when there exists a non-trivial homomorphism between two given equivalence relations arising from translation profinite actions with spectral gap. More generally, we have:

\begin{corollary}\label{bborel}
 Let $\Gamma$ be a residually finite group.  Let $G=\varprojlim\Gamma/\Gamma_n$ be the profinite completion of  $\Gamma$ with respect to a  chain  of finite index normal subgroups with trivial intersection. Assume  that the left translation action $\Gamma\curvearrowright (G,m_G)$ has spectral gap.

Let $\Lambda$ be a countable subgroup of a Polish group $H$ and consider the left translation action $\Lambda\curvearrowright H$.

Then the following hold:
\begin{enumerate}

\item  $\mathcal R(\Gamma\curvearrowright G)\leqslant_{B}\mathcal R(\Lambda\curvearrowright H)$ if and only if we can find an open subgroup $G_0<G$, a closed subgroup $H_0<H$, and a continuous isomorphism $\delta:G_0\rightarrow H_0$ such that $\delta(\Gamma\cap G_0)=\Lambda\cap H_0$.

\item There exists a non-trivial  homomorphism from $\mathcal R(\Gamma\curvearrowright G)$  to $\mathcal R(\Lambda\curvearrowright H)$
if and only if there exist an open subgroup $G_0<G$ and a continuous homomorphism $\delta:G_0\rightarrow H$ such that $\delta(\Gamma\cap G_0)\subset\Lambda$ and $\delta(G_0)\not\subset\Lambda$.

\end{enumerate}
\end{corollary}

\begin{remark}
(1) If any homomorphism from $\mathcal R(\Gamma\curvearrowright G)$  to $\mathcal R(\Lambda\curvearrowright H)$ is trivial, then one says that $\mathcal R(\Gamma\curvearrowright G)$  is  {\it $\mathcal R(\Lambda\curvearrowright H)$-ergodic} (see \cite[Appendix A]{HK05}). 

(2) If a continuous homomorphism $\delta:G_0\rightarrow H$ satisfies $\delta(G_0)\subset\Lambda$, then $\delta(G_0)$ must be finite. This implies that $\delta(G_0)=\delta(\Gamma\cap G_0)$ and the restriction of $\delta$ to some open subgroup $G_1<G_0$ is trivial. 
\end{remark}

{\it Proof.}
(1) To see the {\it if} assertion, assume that there exist an open subgroup $G_0<G$, a closed subgroup $H_0<H$, and a continuous isomorphism $\delta:G_0\rightarrow H_0$ such that $\delta(\Gamma\cap G_0)=\Lambda\cap H_0$. Then $\delta$ witnesses the fact that $\mathcal R((\Gamma\cap G_0)\curvearrowright G_0)\leqslant_{B}\mathcal R(\Lambda\curvearrowright H)$.

Since $\Gamma<G$ is dense, we have that $\Gamma G_0=G$. This implies that we can find $F\subset\Gamma$ finite such that $G=\sqcup_{g\in F}gG_0$. Define a Borel map $\alpha:G\rightarrow G_0$ by 
$\alpha(x)=g^{-1}x$, where $g\in F$ is such that $x\in gG_0$. Then $y\in \Gamma x$ if and only if $\alpha(y)\in (\Gamma\cap G_0)\alpha(x)$, showing that  $\mathcal R(\Gamma\curvearrowright G)\leqslant_{B}\mathcal R((\Gamma\cap G_0)\curvearrowright G_0)$. Altogether, we deduce that $\mathcal R(\Gamma\curvearrowright G)\leqslant_{B}\mathcal R(\Lambda\curvearrowright H)$.

 For the {\it only if} assertion, let
$\theta: G\rightarrow H$ be a Borel  map such that $\Gamma x=\Gamma y$ if and only if $\Lambda\theta(x)=\Lambda\theta (y)$.
By applying Theorem \ref{maintech} we can find an open subgroup $G_0<G$, a continuous homomorphism $\delta:G_0\rightarrow H$,  a Borel map $\phi:G_0\rightarrow \Lambda$ and $h\in H$ such that $\delta(\Gamma\cap G_0)\subset\Lambda$  and \begin{equation}\label{eq}\theta(x)=\phi(x)\delta(x)h,\;\;\;\text{for almost every $x\in G_0$}.\end{equation}
 
We claim that $K:=\ker(\delta)$ is finite.  If $g\in K$, then \ref{eq} implies that $\Lambda\theta(gx)=\Lambda\theta(x)$, for almost every $x\in G_0$. From this we deduce that $\Gamma gx=\Gamma x$, for almost every (thus, for some) $x\in G_0$ and hence $g\in \Gamma$. This shows that $K\subset\Gamma$, so $K$ is countable.
It follows that $K$ is a countable compact group. Since the Haar measure of $K$ has finite mass, it must be the case that $K$ is finite.

Next, since $\ker(\delta)$ is finite, after replacing $G_0$ with a smaller open subgroup of $G$, we may assume that $\delta$ is 1-1. Define $H_0:=\delta(G_0)$, then $H_0$ is a closed subgroup of $H$ and $\delta:G_0\rightarrow H_0$ is a continuous isomorphism.

 We  have that $\delta(\Gamma\cap G_0)\subset \Lambda\cap H_0$ and to finish the proof it remains to show the reverse inclusion. Thus, let $\lambda\in \Lambda\cap H_0$ and let $g\in G_0$ such that $\delta(g)=\lambda$. Then by using \ref{eq} we have that $\Lambda\theta(gx)=\Lambda\delta(gx)h=\Lambda\delta(x)h=\Lambda\theta(x)$, for almost every $x\in G_0$. We deduce that $\Gamma gx=\Gamma x$, for almost every $x\in G_0$. This implies that $g\in\Gamma$ and therefore $\lambda\in\delta(\Gamma\cap G_0)$.

(2) Assume that there is a homomorphism $\theta:G\rightarrow H$  from  $\mathcal R(\Gamma\curvearrowright G)$  to $\mathcal R(\Lambda\curvearrowright H)$  which is non-trivial, in the sense  $\theta(x)$ does not belong to the same $\Lambda$-orbit, for almost every $x\in G$. 
By Theorem \ref{maintech} we can find an open subgroup $G_0<G$, a continuous homomorphism $\delta:G_0\rightarrow H$, a Borel map $\phi:G_0\rightarrow\Lambda$ and $h\in H$ such that $\delta(\Gamma\cap G_0)\subset\Lambda$ and $\theta(x)=\phi(x)\delta(x)h$, for almost every $x\in G_0$. If $\delta(G_0)\subset\Lambda$, then $\theta(x)\in\Lambda h$, for almost every $x\in G_0$. Since $\Gamma G_0=G$, we would get that $\theta(x)\in\Lambda h$, for almost every $x\in G$, contradicting the assumption that $\theta$ is non-trivial.

Conversely, assume that we have a continuous homomorphism $\delta:G_0\rightarrow H$  such that $\delta(\Gamma\cap G_0)\subset\Lambda$ and $\delta(G_0)\not\subset\Lambda$, where $G_0<G$ is an open subgroup. Let $\alpha:G\rightarrow G_0$ be the map defined as in the proof of part (1). Then $\beta:=\delta\circ\alpha:G\rightarrow H$ is a homomorphism from $\mathcal R(\Gamma\curvearrowright G)$  to $\mathcal R(\Lambda\curvearrowright H)$. If $\beta$ is trivial, then we can find $h\in H$ such that $\delta(x)\in\Lambda h$, for almost every $x\in G_0$. This  implies that $\delta(G_0)\subset\Lambda$ and provides a contradiction.
\hfill$\blacksquare$

\subsection{Homomorphism rigidity for compact actions} We end this section by proving analogues of Theorem \ref{maintech} and Corollary \ref{bborel} for translation compact actions.

\begin{theorem}\label{homrig}
 Let $\Gamma$ be a countable group together with  a dense embedding $\tau:\Gamma\hookrightarrow G$ into a connected compact Polish group $G$. 
Assume that $\pi_1(G)$ is finite and 
 the left translation action $\Gamma\curvearrowright (G,m_G)$ has spectral gap.

Let $\Lambda$ be a countable subgroup of a Polish group $H$ and consider the left translation action $\Lambda\curvearrowright H$. Let $\theta:G\rightarrow H$ be a Borel map such that $\theta(\Gamma x)\subset\Lambda\theta(x)$, for almost  every $x\in G$.  Let $w:\Gamma\times G\rightarrow\Lambda$ be the cocycle defined by the formula $\theta(gx)=w(g,x)\theta(x)$, for all $g\in\Gamma$ and almost every $x\in G$.

Then we can find a finite abelian subgroup $\Delta<\Lambda$, such that if $H_1$ denotes the centralizer of $\Delta$ in $H$ and $\pi:H\rightarrow {\Delta}\backslash H$ the quotient, then we can find a continuous homomorphism 
$\delta:G\rightarrow H_1/\Delta$, a Borel map $\phi:G\rightarrow\Lambda$ and $h\in {\Delta}\backslash H$  satisfying

\begin{itemize}
\item $\delta(\Gamma)\subset (\Lambda\cap H_1)/\Delta$,
\item $w'(g,x)=\phi(gx)^{-1}w(g,x)\phi(x)\in\Lambda\cap H_1$ and $\pi(w'(g,x))=\delta(g)$, for all $g\in\Gamma$ and almost every $x\in G$, and
\item $\pi(\phi(x)^{-1}\theta(x))=\delta(x)h$, for almost every $x\in G$.
\end{itemize}
\end{theorem}

{\it Proof.} By repeating verbatim the first part of the proof of Theorem \ref{maintech} we can find $C\in (\frac{31}{32},1)$ and a neighborhood $V$ of $e$ in $G$ such that $$m_G(\{x\in X|w(g,xt)=w(g,x)\})\geqslant C,\;\;\;\text{for all $g\in\Gamma$ and $t\in V$}.$$

By applying Theorem \ref{compact} we find a subgroup $\Lambda_0<\Lambda$, a finite subgroup $\Delta<Z(\Lambda_0)$, a Borel map $\phi:G\rightarrow\Lambda$, and a homomorphism $\delta:\Gamma\rightarrow\Lambda_0/\Delta$ such that $w'(g,x)=\phi(gx)^{-1}w(g,x)\phi(x)\in\Lambda_0$, and $p(w'(g,x))=\delta(g)$, for all $g\in\Gamma$ and almost every $x\in G$. Here $p:\Lambda_0\rightarrow\Lambda_0/\Delta$ is the quotient homomorphism. 
Let $\pi:H\rightarrow {\Delta}\backslash H$ denote the quotient. Since $\Delta<\Lambda_0$ is central, we can identify $\Lambda_0/\Delta={\Delta}\backslash\Lambda_0$. Under this identification we have that $p=\pi_{|\Lambda_0}$, hence $\pi(w'(g,x))=\delta(g)$, for all $g\in\Gamma$ and almost every $x\in G$.

Define $\tilde\theta:G\rightarrow\Delta\backslash H$ by letting $\tilde\theta(x)=\pi(\phi(x)^{-1}\theta(x))$. Let $H_1$ be the centralizer of $\Delta$ in $H$.
Notice that $H_1/\Delta$ and hence $\Lambda_0/\Delta$ acts on $\Delta\backslash H$ by left multiplication. Then we have that \begin{equation}\label{q}\tilde{\theta}(gx)=\delta(g)\tilde{\theta}(x),\;\;\;\text{for all}\;\;\;g\in\Gamma\;\;\;\text{and almost every}\;\;\;x\in G.\end{equation}

Now, define the quotient topological space  $\Sigma:=(H_1/\Delta)\backslash (\Delta\backslash H)$.
Since $\Sigma$ is homeomorphic to $H_1\backslash H$ and  $H_1<H$ is a closed subgroup, \cite[Proposition 1.2.3]{BK96} implies that $\Sigma$ is a Polish space. If $q:\Delta\backslash H\rightarrow\Sigma$ denotes the quotient map, then $q(\tilde\theta(gx))=q(\tilde\theta(x))$, for all $g\in\Gamma$ and almost every $x\in G$. Since $\Gamma<G$ is dense, we get that the map $G\ni x\rightarrow q(\tilde\theta(x))\in\Sigma$ is constant. Thus, we can find $h_0\in \Delta\backslash H$ and a Borel map $\alpha:G\rightarrow H_1/\Delta$ such that $\tilde\theta(x)=\alpha(x)h_0$, for almost every $x\in G$. 

In combination with \ref{q} this implies that $\alpha(gx)h_0=\delta(g)\alpha(x)h_0$, for all $g\in\Gamma$ and almost every $x\in G$. Since the left multiplication action $H_1/\Delta\curvearrowright\Delta\backslash H$ is free, we get that \begin{equation}\label{q2}\alpha(gx)=\delta(g)\alpha(x),\;\;\text{for all}\;\;g\in\Gamma\;\;\text{and}\;\;\text{almost every}\;\;x\in G.\end{equation}

By applying Lemma \ref{ext}, we get that $\delta:\Gamma\rightarrow \Lambda_0/\Delta<H_1/\Delta$ extends to a continuous homomorphism $\delta:G\rightarrow H_1/\Delta$ and that we can find $h_1\in H_1/\Delta$ such that $\alpha(x)=\delta(x)h_1$, for almost every $x\in G$.
Thus, if we let $h=h_1h_0\in\Delta\backslash H$, then $\tilde\theta(x)=\delta(x)h$, for  almost every $x\in G$. This finishes the proof of the theorem.
\hfill$\blacksquare$

As a consequence of Theorem \ref{homrig} we obtain the following analogue of Corollary \ref{bborel} for translation compact actions.

\begin{corollary}\label{Bor}
 Assume the setting from Theorem \ref{homrig}.
Then we have the following:
\begin{enumerate}
\item $\mathcal R(\Gamma\curvearrowright G)\leqslant_{B}\mathcal R(\Lambda\curvearrowright H)$ if and only if  there exist a finite subgroup $\Sigma<\Gamma$ such that $\Sigma<Z(G)$,
a finite subgroup $\Delta<\Lambda$ and a closed subgroup $H_0<H$ such that $\Delta<Z(H_0)$, and a continuous isomorphism $\delta:G/\Sigma\rightarrow H_0/\Delta$ such that $\delta(\Gamma/\Sigma)=(\Lambda\cap H_0)/\Delta$. 
\item There exists a non-trivial homomorphism from $\mathcal R(\Gamma\curvearrowright G)$ to $\mathcal R(\Lambda\curvearrowright H)$ if and only  there exist a finite subgroup $\Delta<\Lambda$ and a closed subgroup $H_0<H$ such that $\Delta<Z(H_0)$, and a non-trivial continuous homomorphism $\delta:G\rightarrow H_0/\Delta$ such that $\delta(\Gamma)\subset(\Lambda\cap H_0)/\Delta$. 
\end{enumerate} 
\end{corollary}

{\it Proof.}
Since the {\it if} assertions of both (1) and (2) are immediate, let us prove the {\it only if} assertions. To this end, 
let $\theta:G\rightarrow H$ be a Borel map such that we have $\theta(x)\in\Lambda\theta(y)$ whenever $x\in\Gamma y$. By Theorem \ref{homrig}, we can find a finite subgroup $\Delta<\Gamma$ such that if $H_1$ is the centralizer of $\Delta$ in $H$ and $\pi:H\rightarrow {\Delta}\backslash H$ the quotient map, then there exists a continuous homomorphism 
$\delta:G\rightarrow H_1/\Delta$, a Borel map $\phi:G\rightarrow\Lambda$ and $h\in {\Delta}\backslash H$  satisfying $$\delta(\Gamma)\subset (\Lambda\cap H_1)/H_1\;\;\;\text{and}\;\;\;\pi(\phi(x)^{-1}\theta(x))=\delta(x)h,\;\;\;\text{for $m_G$-almost every $x\in G$}$$

Define $\bar{\theta}:G\rightarrow H$ by letting $\bar{\theta}(x)=\phi(x)^{-1}\theta(x)$ and $\tilde{\theta}:G\rightarrow\Delta\backslash H$ by putting $\tilde\theta=\pi\circ\bar{\theta}$.

 Let $H_0$ be the unique subgroup of $H_1$ which contains $\Delta$ and satisfies $\delta(G)=H_0/\Delta$.   Then we have that $\delta(\Gamma)\subset (\Lambda\cap H_0)/\Delta$.

In the rest of the proof we derive the {\it only if} assertions by treating separately cases (1) and (2).

(1) Assume that $\theta$ is a reduction, that is $\theta(x)\in\Lambda\theta(y)$ if and only if $x\in\Gamma y$.
 By arguing as in the proof of part (1) of Corollary \ref{bborel}, it follows that $\Sigma:=\ker(\delta)$ is contained in $\Gamma$ and is finite.  Since $G$ is connected, we get that $\Sigma\subset Z(G)$. We denote still by $\delta$ the resulting homomorphism $\delta:G/\Sigma\rightarrow H_0/\Delta$.

Then we have that $\delta(\Gamma/\Sigma)\subset (\Lambda\cap H_0)/\Delta$.
To show the reverse inclusion, let $\lambda\in\Lambda\cap H_0$. Then we can find $g\in G$ such that $\delta(g\Sigma)=\lambda\Delta$. We further have that $\tilde{\theta}(gx)=\delta(g\Sigma)\tilde{\theta}(x)=(\lambda\Delta)\tilde{\theta}(x)$, for almost every $x\in G$.  This implies that $\bar{\theta}(gx)\in\Lambda\bar{\theta}(x)$ and hence $\theta(gx)\in\Lambda\theta(x)$, for almost every $x\in G$. Thus, $gx\in\Gamma x$, for almost every $x\in G$. This implies that $g\in\Gamma$ and so $g\Sigma\in\Gamma/\Sigma$, as claimed. This finishes the proof of (1).

(2) Assume that $\theta$ is not trivial, i.e. $\theta(x)$ is not contained is a single $\Lambda$-orbit, for almost every $x\in G$. We claim that $\delta$ is non-trivial. Suppose by contradiction that $\delta(g)=e$, for all $g\in G$. Then  for every $g\in G$ we have that $\tilde{\theta}(gx)=\tilde{\theta}(x)$, for almost every $x\in G$.
By Fubini's theorem we can find $x_0\in G$ such that $\tilde{\theta}(gx_0)=\tilde{\theta}(x_0)$, for almost every $g\in G$. This clearly implies that $\theta(y)\in\Lambda\theta(x_0)$, for almost every $y\in G$, leading to a contradiction.
\hfill$\blacksquare$

\section{Homomorphism  rigidity for general compact actions}\label{hrig}

In Section \ref{section3}, we proved several rigidity results for translation compact actions with spectral gap. The purpose of this and the next section is to  establish analogous rigidity results for  more general compact actions.  More precisely, in this section we will use Theorem \ref{maintech} and Theorem \ref{homrig} to prove the following result which generalizes both of these theorems.

\begin{theorem}\label{homrig2}
Let $\Gamma$ be a countable group together with a dense embedding $\tau:\Gamma\hookrightarrow G$ into a compact Polish group $G$. Assume that the left translation action $\Gamma\curvearrowright (G,m_G)$ has spectral gap.
Let $\Lambda$ be a countable group together with an embedding  $\rho:\Lambda\hookrightarrow H$ into a locally compact Polish group $H$. Let $H\curvearrowright Y$ be a continuous action on a Polish space $Y$.

Let $\theta:G\rightarrow Y$ be a Borel map such that $\theta(\Gamma x)\subset\Lambda\theta(x)$, for almost every $x\in G$. 
Assume that there exist $k\geqslant 1$ and an $H$-invariant $G_{\delta}$ subset $\Omega\subset Y^k$ such that
\begin{itemize}
\item  the action $H\curvearrowright\Omega$ is smooth and free.
\item $(\theta(x_i))_{i=1}^k\in\Omega$, for $m_G^{\otimes_k}$-almost every $x=(x_i)_{i=1}^k\in G^k$.
\end{itemize}

Then the following hold:

\begin{enumerate}
\item Assume that $G=\varprojlim\Gamma/\Gamma_n$ is a profinite completion of $\Gamma$. 
Then we can find an open subgroup $G_0<G$, a continuous homomorphism $\delta:G_0\rightarrow H$, a Borel map $\phi:G_0\rightarrow\Lambda$, and $y\in Y$  such that $\delta(\Gamma\cap G_0)\subset\Lambda$ and $\theta(x)=\phi(x)\delta(x)y$, for almost every $x\in G_0$.

\item
Assume that $G$ is connected and $\pi_1(G)$ is finite.
Then we can find a finite abelian subgroup $\Delta<\Lambda$, such that if $H_1$ denotes the centralizer of $\Delta$ in $H$ and $\pi:Y\rightarrow\Delta\backslash Y$ the quotient, then we can find a continuous homomorphism $\delta:G\rightarrow H_1/\Delta$, a Borel map $\phi:G\rightarrow\Lambda$, and $y\in\Delta\backslash Y$ satisfying $\delta(\Gamma)\subset (\Lambda\cap H_1)/\Delta$
and $\pi(\phi(x)^{-1}\theta(x))=\delta(x)y$, for almost every $x\in G$. 
(Here, we are using the action of $H_1/\Delta$  on $\Delta\backslash Y$ by left multiplication).\end{enumerate}
\end{theorem}

Since the action $H\curvearrowright\Omega=H$ is free and smooth, Theorem \ref{homrig} generalizes Theorems \ref{maintech} and \ref{homrig}. 
Before proceeding to the proof, let us mention that part (1) of Theorem \ref{homrig2} will be employed  in Section \ref{padic} to prove Corollary \ref{corD}.

{\it Proof.} 
Let $w:\Gamma\times G\rightarrow\Lambda$ be a Borel map such that $\theta(gx)=w(g,x)\theta(x)$, for all $g\in\Gamma$ and almost every $x\in G$. Since the action $\Lambda\curvearrowright Y$ is not assumed free, we cannot conclude that $w$ is a cocycle at this point.
However, as a consequence of the claim below, it will follow that $w$ is a cocycle. 

Since $\Omega\subset Y^k$ is a $G_{\delta}$ set and $Y^k$ is a Polish space,  $\Omega$ is a Polish space (see \cite[Theorem 3.11]{Ke95}). 
Since the action $H\curvearrowright\Omega$ is continuous, free and smooth, Lemma \ref{smooth} implies that there exists a Borel map $p:\Omega\rightarrow H$ such that $p(hx)=hp(x)$, for all $h\in H$ and $x\in\Omega$.

The rest of the proof relies on the following:

{\bf Claim.} There exists a sequence of Borel maps $\phi_n:G\rightarrow H$ such that for all $g\in\Gamma$ we have that $m_G(\{x\in G|\;w(g,x)=\phi_n(gx)\phi_n(x)^{-1}\})\rightarrow 1$, as $n\rightarrow\infty$.

{\it Proof of the claim}. It suffices to show that for every $\varepsilon>0$ and any finite subset $F\subset\Gamma$, we can find a Borel map $\phi:G\rightarrow H$ such that $m_G(\{x\in G|\;w(g,x)=\phi(gx)\phi(x)^{-1}\})\geqslant 1-\varepsilon$, for all $g\in F$.

Since $\Lambda$ is countable, we can find a neighborhood $V$ of the identity $e$ in $G$ such that
\begin{equation}\label{constant}
m_G(\{x\in G|\;w(g,xt)=w(g,x)\})\geqslant 1-\frac{\varepsilon}{k|F|},\;\;\;\text{for all}\;\;\; g\in F\;\;\;\text{and}\;\;\; t\in V.
\end{equation}

Recall that $(\theta(x_i))_{i=1}^k\in\Omega$, for almost every $(x_i)_{i=1}^k\in G^k$. It follows that for all $x\in G$ we have that $(\theta(xt_i))_{i=1}^k\in\Omega$, for almost every $(t_i)_{i=1}^k\in G^k$. Since $m_G^{\otimes_k}(V^k)=m_G(V)^k>0$, by  Fubini's theorem we derive that there exist $t_1,...,t_k\in V$ such that $(\theta(xt_i))_{i=1}^k\in\Omega$, for almost every $x\in G$.
 
Further, we define $\psi:G\rightarrow Y^k$ by letting $\psi(x):=(\theta(xt_1),...,\theta(xt_k))$. Then $\psi(x)\in \Omega$, for almost every $x\in G$. 
We denote by $C$ the set of all $x\in G$ which satisfy the following three conditions:

\begin{itemize}
\item  $\psi(x)\in\Omega$ and $\psi(gx)\in\Omega$, for all $g\in F$,
\item  $\theta(gxt_i)=w(g,xt_i)\theta(xt_i)$, for all $g\in F$ and $1\leqslant i\leqslant k$, and
\item $w(g,x)=w(g,xt_i)$, for all $g\in F$ and $1\leqslant i\leqslant k$. 
\end{itemize}

Since $\psi(x)\in \Omega$ and $\theta(gx)=w(g,x)\theta(x)$, for all $g\in \Gamma$ and almost every $x\in G$, equation \ref{constant} implies that $m_G(C)\geqslant 1-\varepsilon$.

Finally, we define $\phi:G\rightarrow H$ by letting $\phi(x)=\begin{cases}\pi(\psi(x)),\;\;\text{if}\;\; x\in \psi^{-1}(\Omega) \\ e,\;\;\text{if}\;\; x\not\in\psi^{-1}(\Omega)\end{cases}$. 

Then for every $x\in C$ and $g\in F$ we have that $$\phi(gx)=\pi(\psi(gx))=\pi(\theta(gxt_1),...,\theta(gxt_k))=\pi(w(g,xt_1)\theta(xt_1),...,w(g,xt_k)\theta(xt_k))=$$ $$\pi(w(g,x)\theta(xt_1),...,w(g,x)\theta(xt_k))=w(g,x)\pi(\theta(xt_1),...,\theta(xt_k))=w(g,x)\phi(x).$$ 
 
This finishes the proof of the claim. \hfill$\square$

Now, since the action $\Gamma\curvearrowright G$ is strongly ergodic, the claim and Lemma \ref{I} imply that there exists a Borel map $\psi:G\rightarrow H$ such that $w(g,x)=\psi(gx)\psi(x)^{-1}$, for all $g\in\Gamma$ and almost every $x\in G$.  
We are therefore in position to apply Theorems \ref{maintech} and \ref{homrig}, in the cases (1) and (2), respectively.

 Assume that we are in case (1). By applying Theorem \ref{maintech} we can find an open subgroup $G_0<G$, a continuous homomorphism $\delta:G_0\rightarrow H$ and a Borel map $\phi:G_0\rightarrow\Lambda$ such that we have $w(g,x)=\phi(gx)\delta(g)\phi(x)^{-1}$, for all $g\in\Gamma\cap G_0$ and $m_G$-almost every $x\in G_0$.

Thus, if we define $\tilde\theta:G_0\rightarrow Y$ by $\tilde\theta(x)=\phi(x)^{-1}\theta(x)$, then $\tilde\theta(gx)=\delta(g)\tilde\theta(x)$, for all $g\in\Gamma\cap G_0$ and almost every $x\in G_0$. Since $\delta$ is continuous and $\Gamma\cap G_0<G_0$ is dense, it follows that for every $g\in G_0$ we have  that $\tilde\theta(gx)=\delta(g)\tilde\theta(x)$, for almost every $x\in G_0$. By Fubini's theorem, we can find $x_0\in G_0$ such that 
$\tilde\theta(gx_0)=\delta(g)\tilde\theta(x_0)$, for almost every $g\in G_0$. Denoting $y=\delta(x_0^{-1})\tilde\theta(x_0)\in Y$, we get that $\tilde\theta(g)=\delta(g)y$, for almost every $g\in G_0$. This clearly implies the conclusion.

Finalizing the proof in case (2) is  similar to the above and we leave the details to the reader.
\hfill$\blacksquare$

\section{Orbit equivalence rigidity for compact actions}

In this section, we  study orbit equivalences between general compact actions with spectral gap. 
More precisely, under the assumption that $G$ and $H$ are connected, we provide conditions which ensure that any orbit equivalence between the actions $\Gamma\curvearrowright G/K$ and $\Lambda\curvearrowright H/L$  comes from a conjugacy between them.
This result will be used in the Section \ref{tree} to give examples of ergodic treeable equivalence relations without outer automorphisms.

\begin{theorem}\label{conjugacy} Let $\Gamma$ and $\Lambda$ be two countable groups together with dense embeddings $\tau:\Gamma\hookrightarrow G$ and  $\sigma:\Lambda\hookrightarrow H$ into compact Polish groups $G$ and $H$. Let $K<G$ and $L<H$ be closed subgroups.  Assume that the following conditions are satisfied:

\begin{enumerate} 

\item The actions  $\Gamma\curvearrowright (G/K,m_{G/K})$ and $\Lambda\curvearrowright (H/L,m_{H/L})$ are  free.
\item The left translation action $\Gamma\curvearrowright (G,m_G)$ has spectral gap.
\item $G$ and $H$ are connected, $\pi_1(G)$ is finite, $\cap_{g\in G}gKg^{-1}=\{e\}$ and $\cap_{h\in H}hLh^{-1}=\{e\}$.
\item $\Gamma$ has no non-trivial finite normal subgroups and $\Lambda$ has infinite conjugacy classes (icc).
\item For any $h\in\Lambda\setminus\{e\}$, we have $$(m_H\times m_H)(\{(x,y)\in H\times H\;|\;h\in xLx^{-1}yLy^{-1}\})=0.$$
\end{enumerate}
Let $B\subset H/L$ be a Borel set with  $m_{H/L}(B)>0$. Endow $B$ with the probability measure obtained by restricting and rescaling $m_{H/L}$.
Let $\theta:G/K\rightarrow B$ be an isomorphism of probability spaces such that $\theta(\Gamma x)=\Lambda\theta(x)\cap B$, for $m_{G/K}$-almost every $x\in G/K$.

Then $B=H/L$, almost everywhere, and we can find a continuous isomorphism $\delta:G\rightarrow H$, $y\in H$,  and $\alpha\in [\mathcal R(\Lambda\curvearrowright H/L)]$ such that
\begin{itemize}
\item $\delta(\Gamma)=\Lambda$ and $\delta(K)=yLy^{-1}$.
\item $\tilde\theta:=\alpha\circ\theta:G/K\rightarrow H/L$ is given by $\tilde\theta(gK)=\delta(g)yL$, for $m_G$-almost every $g\in G$.
\end{itemize}
\end{theorem}

\begin{remark}
The technical condition (5) is imposed so that we can adapt the stategy of proof of Theorem \ref{maintech} to this new context. Condition (5) is trivially satisfied if $L=\{e\}$. More interestingly, it also holds if $L=SO(n)<H=SO(n+1)$, for $n\geqslant 2$ (see Claim 2 in the proof of Theorem \ref{OUT}).
\end{remark}

{\it Proof.} Let $\pi:G\rightarrow G/K$ denote the quotient map. Define the cocycle $w:\Gamma\times G/K\rightarrow \Lambda$ given by the formula $\theta(gx)=w(g,x)\theta(x)$. Denote $W=w\circ(\text{id}_{\Gamma}\times\pi):\Gamma\times G\rightarrow\Lambda$ and $\Theta=\theta\circ\pi:G\rightarrow H/L$. 
The proof of Theorem \ref{conjugacy} relies on four claims.

{\bf Claim 1.} There exist a subgroup $\Lambda_0<\Lambda$, a finite central subgroup $\Lambda_1<Z(\Lambda_0)$, a Borel map $\phi:G\rightarrow\Lambda$  and a homomorphism $\delta:\Gamma\rightarrow\Lambda_0/\Lambda_1$ such that 
$W'(g,x):=\phi(gx)W(g,x)\phi(x)^{-1}\in\Lambda_0$ and $p(W'(g,x))=\delta(g)$, for all $g\in\Gamma$ and almost every $x\in G$, where $p:\Lambda_0\rightarrow\Lambda_0/\Lambda_1$ denotes the quotient homomorphism.

{\it Proof of Claim 1}. We will derive the conclusion by applying Theorem \ref{compact}.
To this end, note that $W:\Gamma\times G\rightarrow\Lambda$ is a cocycle and we have \begin{equation}\label{eq1}\Theta(gx)=W(g,x)\Theta(x),\;\;\;\text{for all}\;\;\;g\in\Gamma\;\;\;\text{and}\;\;\; \text{almost every $x\in G$}.\end{equation}

 We start by adapting the beginning of the proof of Theorem \ref{maintech}.
Let $\varepsilon\in (0,\frac{1}{64})$. Since the action $\Gamma\curvearrowright (G,m_G)$ has spectral gap, Lemma \ref{spectral} provides $\delta>0$ and $F\subset\Gamma$ finite such that whenever $Y$ is a standard Borel space and $\rho:G\rightarrow Y$ is a Borel map satisfying $m_G(\{x\in G|\rho(gx)=\rho(x)\})\geqslant 1-\delta$, for all $g\in F$, we can find $y\in Y$ such that $m_G(\{x\in G|\rho(x)=y\})\geqslant 1-\varepsilon$. Thus, we have that $\mu(\{x\in X|\rho(gx)=\rho(x)\})\geqslant 1-2\varepsilon$, for all $g\in\Gamma$.

 Since $\Lambda$ is countable,
we can find a neighborhood $V$ of $e\in G$ such that \begin{equation}\label{eq2}m_G(\{x\in G|W(g,xt)=W(g,x)\})\geqslant 1-\delta,\;\;\;\text{for all}\;\;\; g\in F\;\;\;\text{and every}\;\;\; t\in V.\end{equation}

Let $Y$ denote the double coset space $L\backslash H/L$. Fix $t\in V$ and define $\rho_t:G\rightarrow Y$ by letting $\rho_t(x)=L\Theta(x)^{-1}\Theta(xt)L.$
 By equation \ref{eq1}, for almost every $x\in G$, we have that \begin{equation}\label{eq3}\rho_t(gx)=L\Theta(gx)^{-1}\Theta(gxt)L=L\Theta(x)^{-1}W(g,x)^{-1}W(g,xt)\Theta(xt)L\end{equation} 
Combining \ref{eq2} and \ref{eq3} implies that  $m_G(\{x\in G|\rho_t(gx)=\rho_t(x)\})\geqslant 1-\delta$, for all $g\in F$.
Since $Y$ is a standard Borel space, by using the above consequence of the spectral gap property, we get that \begin{equation}\label{eq4}m_G(\{x\in G|\rho_t(gx)=\rho_t(x)\})\geqslant 1-2\varepsilon,\;\;\;\text{for all $g\in\Gamma$\;\;\;\text{and every $t\in V$.}}\end{equation}

Now,  let $Z_t$ be the set of $x\in G$ for which there is $h\in\Lambda\setminus\{e\}$ such that $h\in\Theta(x)L\Theta(x)^{-1}\Theta(xt)L\Theta(xt)^{-1}$. We claim that $m_G(Z_t)=0$, for almost every $t\in G$. By Fubini's theorem, it suffices to show that the set $Z$ of $(x,y)\in G\times G$ for which there is $h\in\Lambda\setminus\{e\}$ such that $h\in\Theta(x)L\Theta(x)^{-1}\Theta(y)L\Theta(y)^{-1}$ has measure zero. 
To see this, note that  assumption (5) implies that $(m_{H/L}\times m_{H/L})((\Theta\times\Theta)(Z))=0$. Since $\theta$ is  measure preserving, it follows that $(m_G\times m_G)(Z)=0$.

Next, by \ref{eq3} we have that $\{x\in G|\rho_t(gx)=\rho_t(x)$ and $W(g,x)\not=W(g,xt)\}\subset Z_t$. Since $m_G(Z_t)=0$, by \ref{eq4} we deduce that $m_G(\{x\in G|W(g,x)=W(g,xt)\})\geqslant 1-2\varepsilon,$ for all $g\in\Gamma$ and every $t\in V$. 
 Since $1-2\varepsilon>\frac{31}{32}$, $G$ is connected and $\pi_1(G)$ is finite, part (2) of Theorem \ref{compact}  gives the claim. \hfill$\square$
\vskip 0.05in

We continue with the following:

{\bf Claim 2}. $\delta$ is injective.

{\it Proof of Claim 2}. Since $\Gamma$ has no non-trivial finite normal subgroups, in order to show that $\delta$ is injective, it suffices to argue that $\ker(\delta)$ is finite. Assume by contradiction that $\ker(\delta)$ is infinite and let $F\subset\Lambda$ be a finite set such that $m_{G}(\{x\in G|\phi(x)\in F\})>\frac{2}{3}$.  
Claim 1  implies that
 $$m_{G/K}(\{x\in G/K|w(g,x)\in F^{-1}\Lambda_1F\})>\frac{1}{3},\;\;\;\text{for all $g\in\ker{\delta}$}.$$
 
  Since $\ker(\delta)$ is assumed infinite, we can find $h\in F^{-1}\Lambda_1F$ and a sequence $\{g_n\}_n$ of distinct elements from $\ker(\delta)$ such that if  $X_n=\{x\in G/K|w(g_{n},x)=h\}$, then $m_{G/K}(X_n)>\frac{1}{3|F^{-1}\Lambda_1F|}$, for all $n$. Since the sets $\{X_n\}_n$ are mutually disjoint, this provides a contradiction. \hfill$\square$

\vskip 0.05in

{\bf Claim 3.} $\Lambda_1=\{e\}$, and $\delta(\Gamma)<\Lambda$ has finite index.

{\it Proof of Claim 3.} Let $\Sigma:=p^{-1}(\delta(\Gamma))<\Lambda_0$. We claim that $\Sigma<\Lambda$ has finite index. Assume by contradiction that $\Sigma<\Lambda$ has infinite index.  We denote $\mathcal R=\mathcal R(\Lambda\curvearrowright H/L)$ and $\mathcal S=\mathcal R(\Sigma\curvearrowright H/L)$. 
We define $\varphi_{\mathcal S}:[[\mathcal R]]\rightarrow [0,1]$ by letting $$\varphi_{\mathcal S}(\beta)=m_{H/L}(\{x\in H/L|\; \text{$\beta(x)$ is defined and}\;\; \beta(x)\in [x]_{\mathcal S}\}),\;\;\text{for every $\beta\in [[\mathcal R]]$}.$$

Since $m_{G}(\{x\in G|\phi(x)\in F\})>\frac{2}{3}$, Claim 1 implies that \begin{equation}\label{ineq1} m_{G/K}(\{x\in G/K|\; \theta(gx)\in F^{-1}\Sigma F\;\theta(x)\})>\frac{1}{3},\;\;\text{for all $g\in\Gamma$}.\end{equation} 

Now, for every $g\in\Gamma$, we define $\alpha_g=\theta g\theta^{-1}:B\rightarrow B$. Then $\alpha_g\in [[\mathcal R]]$ and equation \ref{ineq1} gives that \begin{equation}\label{ineq2}\sum_{h,k\in  F} \varphi_{\mathcal S}(h\alpha_gk^{-1})>\frac{1}{3},\;\;\;\text{for all $g\in\Gamma$.} \end{equation}

Note that the equivalence relation associated to the action $(\alpha_g)_{g\in\Gamma}$ of $\Gamma$ on $B$ is equal to $\mathcal R_{|B}$. By using \ref{ineq2} and  a straightforward modification of the proof of \cite[Theorem 2.5]{IKT08}, it follows that we can find a Borel set $A\subset H/L$ of positive measure and $\kappa\geqslant 1$ such that every ${\mathcal R}_{|A}$-class contains at most $\kappa$ ${\mathcal S}_{|A}$-classes. In other words, the inclusion ${\mathcal S}_{|A}\subset{\mathcal R}_{|A}$ has index at most $\kappa$.

We will show that this contradicts our assumption that $\Sigma<\Lambda$ has infinite index.
Let us first show that if $g_1,g_2,...,g_{\kappa+1}\in\Lambda$ and $g_i\Sigma\not=g_j\Sigma$, for all $i\not=j$, then $m_{H/L}(\cap_{i=1}^{\kappa+1}g_iA)=0$. Indeed, if $x\in \cap_{i=1}^{\kappa+1}g_iA$, let $y=g_1^{-1}x$. Then $y\in A$ and $(g_i^{-1}g_1)y\in [y]_{\mathcal R}\cap A=[y]_{\mathcal R_{|A}}$, for all $2\leqslant i\leqslant \kappa+1$. Since $\Sigma\not=\Sigma(g_i^{-1}g_1)\not=\Sigma(g_j^{-1}g_1)$, for all $i,j\in\{2,...,\kappa+1\}$ with $i\not=j$, we deduce that the $\mathcal S_{|A}$-classes of $y,(g_2^{-1}g_1)y,...,(g_{\kappa+1}^{-1}g_1)y$ are disjoint. This would imply that $[y]_{\mathcal R_{|A}}$ contains more than $\kappa$ $\mathcal S_{|A}$-classes, which is a contradiction.

Next, let $p$ be the smallest natural number such that whenever $g_1,g_2,...,g_p\in\Lambda$ are such that $g_i\Sigma\not=g_j\Sigma$, for all $1\leqslant i<j\leqslant p$, then $m_{H/L}(\cap_{i=1}^{p}g_iA)=0$. Then we clearly have that $1<p\leqslant \kappa+1$. Let $q$ be a natural number such that $q<p$ and $2q\geqslant p$. Since $q<p$ we can find $g_1,g_2,...,g_q\in\Lambda$ such that $g_i\Sigma\not=g_j\Sigma$, for all $i\not=j$, and $\tilde A=\cap_{i=1}^qg_iA$ satisfies $m_{H/L}(\tilde A)>0$.

Since $\Sigma<\Lambda$ has infinite index, for any finite set $S\subset\Lambda$, we have that $S\Sigma S\not=\Lambda$. Using this fact, we can find a sequence $\{h_a\}_{a=1}^{\infty}$ of elements in $\Lambda$ such that $h_a(\cup_{i=1}^qg_i\Sigma)\cap h_b(\cup_{i=1}^qg_i\Sigma)=\emptyset$, for all $a\not= b$. Since $h_a\tilde A\cap h_b\tilde A=(\cap_{i=1}^qh_ag_iA)\cap(\cap_{i=1}^qh_bg_iA)$ and $2q\geqslant p$, we conclude that $m_{H/L}(h_a\tilde A\cap h_b\tilde A)=0$, for all $a\not=b$. Since $m_{H/L}(\tilde A)>0$, this gives the desired contradiction.

We have therefore proved that $\Sigma<\Lambda$ has finite index. Since $\Sigma<\Lambda_0$ and $\Lambda_1<Z(\Lambda_0)$, we have that $\Sigma$ commutes with $\Lambda_1$. Since $\Sigma<\Lambda$ has finite index, we get that $\{ghg^{-1}|g\in\Lambda\}$ is finite, for every $h\in\Lambda_1$. As $\Lambda$ is assumed icc, we get that $\Lambda_1=\{e\}$. Hence, $\delta(\Gamma)=\Sigma$ has finite index in $\Lambda$.
\hfill$\square$

Next, we show that: 

{\bf Claim 4.} $\phi:G\rightarrow\Lambda$ factors through $\pi:G\rightarrow G/K$.

{\it Proof of Claim 4.} Let $k\in K$.  Claim 1 implies that $$\phi(gxk)\delta(g)\phi(xk)^{-1}=w(g,\pi(x))=\phi(gx)\delta(g)\phi(x)^{-1},$$ for all $g\in\Gamma$ and almost every $x\in G$. Hence, if we denote $\lambda(x)=\phi(x)^{-1}\phi(xk)$, then we have that $\lambda(gx)=\delta(g)\lambda(x)\delta(g)^{-1}$, for all $g\in\Gamma$ and almost every $x\in G$. 

For $h\in\Lambda$, let $A_h=\{x\in G|\lambda(x)=h\}$. Then $gA_h=A_{\delta(g)h\delta(g)^{-1}}$ and  $m_{H/L}(A_h)=m_{H/L}(A_{\delta(g)h\delta(g)^{-1}})$, for all $g\in\Gamma$.
Since  $\delta(\Gamma)<\Lambda$ has finite index by Claim 2 and $\Lambda$ is icc,  $\{\delta(g)h\delta(g)^{-1}|g\in\Gamma\}$ is infinite, for any $h\in\Lambda\setminus\{e\}$. Since the sets $\{A_h\}_{h\in\Lambda}$ are mutually disjoint, we conclude that $m_{H/L}(A_h)=0$, for all $h\in\Lambda\setminus\{e\}$. This implies that $\phi(x)=\phi(xk)$, for almost every $x\in G$. Therefore, $\phi$ factors through $\pi:G\rightarrow G/K$, as claimed. \hfill$\square$

We are now ready to finish the proof of Theorem \ref{conjugacy}.
By Claim 4, we can define $\tilde\theta:G/K\rightarrow H/L$ by letting $\tilde\theta(x)=\phi(x)^{-1}\theta(x)$. Also, let $\tilde B=\tilde\theta(G/K)$.  Claims 1 and 3 give that $\tilde\theta(gx)=\delta(g)\tilde\theta(x)$, for all $g\in\Gamma$ and almost every $x\in G/K$.

Next, let $x\in G/K$ and $y\in\Gamma x$ such that $\tilde\theta(x)=\tilde\theta(y)$. If  $y=gx$, for some $g\in\Gamma$, then we get that $\tilde\theta(y)=\delta(g)\tilde\theta(x)=\delta(g)\tilde\theta(y)$. By freeness, we derive that $\delta(g)=e$. Since $\delta$ is injective we deduce that $g=e$ and hence $x=y$. This argument implies
 that there exists  $\alpha\in [\mathcal R(\Lambda\curvearrowright H/L)]$ such that $\tilde\theta=\alpha\circ\theta$. 
 In particular, it follows that $\tilde\theta(\Gamma x)=\Lambda\tilde\theta(x)\cap\tilde B$, for $m_{G/K}$-almost every $x\in G/K$.
 
 Now, since $\delta(\Gamma)<\Lambda$ has finite index and $\Lambda<H$ is dense, it follows that $H_0=\overline{\delta(\Gamma)}<H$ is a finite index closed subgroup. Since $H$ is connected, we derive that $H_0=H$, hence $\delta(\Gamma)<H$ is dense. Therefore, the left translation actions of $\delta(\Gamma)$ on $H$ and on $H/L$ are ergodic. Since $\tilde B\subset H$ is a $\delta(\Gamma)$-invariant set of positive measure, we derive that  $\tilde B=K/L$ and hence $B=K/L$, almost everywhere. 
 
 By combining the last two paragraphs, we deduce that $\tilde\theta(\Gamma x)=\Lambda\tilde\theta(x)$, for almost every $x\in G/K$. Since $\tilde\theta(\Gamma x)=\delta(\Gamma)\tilde\theta(x)$, for almost every $x\in G/K$,  we conclude that $\delta(\Gamma)=\Lambda$.
 
 Altogether, we get that $\tilde\theta:(G/K,m_{G/K})\rightarrow (H/L,m_{H/L})$ is an isomorphism of probability spaces such that $\tilde\theta\Gamma{\tilde\theta}^{-1}=\Lambda$.
Since $\cap_{g\in G}gKg^{-1}=\{e\}$ and $\cap_{h\in H}hLh^{-1}=\{e\}$, the closure of $\Gamma$ in Aut$(G/K,m_{G/K})$ is equal to $G$, and the closure of $\Lambda$ in Aut$(H/L,m_{H/L})$ is equal to $H$. We therefore deduce that $\tilde\theta  G{\tilde\theta}^{-1}=H$. Thus, $\delta$ extends to a continuous isomorphism $\delta:G\rightarrow H$ so that for all $g\in G$ we have  $\tilde\theta(gx)=\delta(g)\tilde\theta(x)$, for almost every $x\in G/K$. 

By Fubini's theorem, there is $x_0K\in G/K$ such that $\tilde\theta(gx_0K)=\delta(g)\tilde\theta(x_0K)$, for almost every $g\in G$. Let $y\in H$ such that $\delta(x_0)^{-1}\tilde\theta(x_0K)=yL$. 
Then $\tilde\theta(gK)=\delta(g)yL$, for almost every $g\in G$. 
Finally, by using this identity, it is easy to see that $\delta(K)=yLy^{-1}$. 
\hfill$\blacksquare$

Specializing Theorem \ref{conjugacy} to the case when $K=L=\{e\}$, we obtain the following analogue of the first part of Theorem \ref{mainthm} for certain translation compact actions.

\begin{corollary}\label{OE}
 Let $\Gamma$ and $\Lambda$ be two countable icc groups together with dense embeddings $\tau:\Gamma\hookrightarrow G$ and  $\sigma:\Lambda\hookrightarrow H$ into connected compact metrizable groups $G$ and $H$. 
Assume that $\pi_1(G)$ and $\pi_1(H)$ are finite.
Assume that the left translation action $\Gamma\curvearrowright (G,m_G)$ has spectral gap.

Then the actions $\Gamma\curvearrowright (G,m_G)$ and $\Lambda\curvearrowright (H,m_H)$ are stably orbit equivalent if and only if  there exists a continuous isomorphism $\delta:G\rightarrow H$ such that $\delta(\Gamma)=\Lambda$.
\end{corollary}

{\it Proof.}
Assume that the actions $\Gamma\curvearrowright (G,m_G)$ and $\Lambda\curvearrowright (H,m_H)$ are stably orbit equivalent. Since the action $\Gamma\curvearrowright (G,m_G)$ is strongly ergodic, it follows that $\Lambda\curvearrowright (H,m_H)$ is also strongly ergodic and hence has spectral gap  by Proposition \ref{tau}. The conclusion then follows from Theorem \ref{conjugacy}.
\hfill$\blacksquare$

\section {Proofs of Corollaries \ref{corB} and \ref{corC}}\label{sect4}

The aim of this section is to prove Corollaries \ref{corB} and \ref{corC}.
We start by introducing some notation.

\begin{notation}\label{notation}
 Let $S$ be a set of primes. 
\begin{itemize}
\item  We define the profinite groups $$G_S=\prod_{p\in S}SL_2(\mathbb F_p),\;\;\;\; H_S=\prod_{p\in S}GL_2(\mathbb F_p)\;\;\;\;\text{and}\;\;\;\;\;K_S=\prod_{p\in S}SL_2(\mathbb Z_p).$$

\item For every $p$, we denote by $\pi_p: SL_2(\mathbb Z)\rightarrow  SL_2(\mathbb F_p)$ the reduction modulo $p$ and by \\ $\tau_p:SL_2(\mathbb Z)\rightarrow SL_2(\mathbb Z_p)$ the natural embedding. 

\item Whenever $S$ is infinite, we view $SL_2(\mathbb Z)$ as a subgroup of $G_S$ via the diagonal embedding $\pi_S=(\pi_p)_{p\in S}$.

\item We  view $SL_2(\mathbb Z)$ as a subgroup of $K_S$ via the diagonal embedding $\tau_S=(\tau_p)_{p\in S}$.

\item For a subgroup $\Gamma<SL_2(\mathbb Z)$, we denote by $G_{\Gamma,S}$ and $K_{\Gamma,S}$ the closure of $\Gamma$ in $G_S$ and $K_S$, respectively.

\end{itemize}
\end{notation}

Next, let us record the following fact that we will use several times in the next three sections:

\begin{fact}\label{strong} 
If $\Gamma<SL_2(\mathbb Z)$ is a non-amenable subgroup, then 
\begin{itemize}
\item $G_{\Gamma,S}<G_S$ and $K_{\Gamma,S}<K_S$ are open subgroups (see Example \ref{profinite}), and 
\item the left translation actions $\Gamma\curvearrowright G_{\Gamma,S}$ and $\Gamma\curvearrowright K_{\Gamma,S}$ have spectral gap (see Example \ref{selberg}).
\end{itemize}
\end{fact}

We continue with several elementary rigidity results about  homomorphisms between the groups defined above. We start by collecting some properties of $SL_2(\mathbb F_p)$ that we will use repeatedly.

\begin{proposition}\label{sl_2} Let $p$ be a prime and denote by $I\in SL_2(\mathbb F_p)$ the identity matrix. Then  
\begin{enumerate}
\item $|SL_2(\mathbb F_p)|=p(p^2-1)$.
\item If $p\geqslant 5$, then the only proper normal subgroups of $SL_2(\mathbb F_p)$ are $\{I\}$ and $\{\pm I\}$. 
\item Any automorphism of $SL_2(\mathbb F_p)$ is given by conjugation with an element of $GL_2(\mathbb F_p)$.
\item Any proper subgroup $L<SL_2(\mathbb F_p)$  with $|L|>60$ is  $2$-step solvable. In particular, for every $a,b,c,d\in L$ we have that $[[a,b],[c,d]]=I$. 
\end{enumerate}
\end{proposition}

The first three facts are well-known. The last fact is a consequence of Dickson's classification of subgroups of $SL_2(\mathbb F_p)$ (see \cite [Theorem 7 and Proposition 3]{BG05} and the references therein). 

Below, whenever $S\subset T$, we consider the natural embeddings $G_S<G_T$ and $K_S<K_T$.

\begin{lemma}\label{lemma1} Let $S$ and $T$ be two  sets of primes $\geqslant 7$. 

If $\delta:G_S\rightarrow G_T$ is an injective homomorphism, then $S\subset T$ and we can find $g\in H_S$ such that $\delta(x)=gxg^{-1}$, for all $x\in G_S$. In particular, $\delta(G_S)=G_S$.
\end{lemma}

{\it Proof.}  The lemma is an immediate consequence of the following:

\noindent {\bf Claim.} If $p\geqslant 7$ and $q\geqslant 5$ are primes and $\rho:SL_2(\mathbb F_p)\rightarrow SL_2(\mathbb F_q)$ is a non-trivial homomorphism, then $p=q$ and we can find $h\in GL_2(\mathbb F_p)$ such that $\rho(x)=hxh^{-1}$, for all $x\in SL_2(\mathbb F_p)$.

{\it Proof of the Claim}. Since $\ker(\rho)$ is a proper normal subgroup of $SL_2(\mathbb F_p)$, we have that either $\ker(\rho)=\{I\}$ or $\ker(\rho)=\{\pm I\}$. Thus, $\rho(SL_2(\mathbb F_p))$ is isomorphic to either $SL_2(\mathbb F_p)$ or $PSL_2(\mathbb F_p)$. As $p\geqslant 5$, neither of these groups is solvable (they are perfect groups, by Proposition \ref{sl_2}, part (2)). Hence, $\rho(SL_2(\mathbb F_p))$ is not solvable. Moreover,  it has cardinality either $p(p^2-1)$ or $\frac{p(p^2-1)}{2}$.

 Since $\frac{p(p^2-1)}{2}>60$, by using Proposition \ref{sl_2}, part (4), we deduce that $\rho(SL_2(\mathbb F_p))=SL_2(\mathbb F_q)$. Since the equation $q(q^2-1)=\frac{p(p^2-1)}{2}$ has no solutions $p\geqslant 7$, $q\geqslant 5$, we conclude that $p=q$ and that $\rho$ is an automorphism of $SL_2(\mathbb F_p)$.
By Proposition \ref{sl_2}, part (3), we are done.
\hfill$\blacksquare$

\begin{corollary} \label{cor2} Let $S$ and $T$ be two infinite sets of primes. Let $G<G_S$ be an open subgroup  and $\Gamma<SL_2(\mathbb Z)$ be an infinite subgroup. 

If there exists an injective  continuous homomorphism $\delta:G\rightarrow G_T$ such that $\delta(\Gamma\cap G)\subset SL_2(\mathbb Z)$, then $|S\Delta T|<\infty$.

\end{corollary}

{\it Proof.}
Let $T_0=T\setminus\{2,3,5\}$. Since $G_{T_0}<G_T$ is an open subgroup and $\delta$ is continuous, we get that $\delta^{-1}(G_{T_0})<G$ is an open subgroup. 
Thus, $\delta^{-1}(G_{T_0})<G_S$ is an open subgroup and therefore we can find a set $S_0\subset S\setminus\{2,3,5\}$ such that $G_{S_0}\subset \delta^{-1}(G_{T_0})$ and $S\setminus S_0$ is finite. 
Altogether, we have that $\delta_{|G_{S_0}}:G_{S_0}\rightarrow G_{T_0}$ is an injective homomorphism.

By  Lemma \ref{lemma1} we get that $S_0\subset T_0$ and that $\delta(G_{S_0})=G_{S_0}$. We claim that $T_0\setminus S_0$ is finite.
Assume by contradiction that $T_0\setminus S_0$ is infinite. Then we have that $SL_2(\mathbb Z)\cap G_{S_0}=\{e\}$ (where we view both $SL_2(\mathbb Z)$ and $G_{S_0}$ as subgroups of $G_{T_0})$. On the other hand,  since $\delta$ is injective, $\delta(\Gamma\cap G_{S_0})$ is an infinite subgroup of $SL_2(\mathbb Z)\cap \delta(G_{S_0})=SL_2(\mathbb Z)\cap G_{S_0}$.  This gives a contradiction.

Finally, the claim implies that $|S_0\Delta T_0|<\infty$ and hence that $|S\Delta T|<\infty$.
\hfill$\blacksquare$

For the rest of this section, for $m\geqslant 1$, we denote $SL_2(m\mathbb Z)=\ker(SL_2(\mathbb Z)\rightarrow SL_2(\mathbb Z/m\mathbb Z))$.

\begin{lemma}\label{lema2} Let $S$ be an infinite set of primes and $m,n\geqslant 1$. Assume that $p\nmid m$, for every $p\in S$.  
Assume that there exists $g\in H_S=\prod_{p\in S}GL_2(\mathbb F_p)$ such that  $g\pi_S(SL_2(m\mathbb Z))g^{-1}=\pi_S(SL_2(n\mathbb Z))$. 

Then  $m=n$ and we can find $k\in GL_2(\mathbb Z)$ and $l$ in the center of $H_S$ such that $g=\pi_S(k)l$.
\end{lemma} 

{\it Proof.} Write $g=(g_p)_{p\in S}$.    
 Denote by $A$ the set of $x\in\mathbb M_2(\mathbb Z)$ for which there exists $y\in\mathbb M_2(\mathbb Z)$ such that $g_p\pi_p(x)g_p^{-1}=\pi_p(y)$, for all $p\in S$.
 Since $S$ is infinite, such an element $y$ must be unique and we denote it by $\phi(x)$. Moreover, $A$ is a subring of $\mathbb M_2(\mathbb Z)$  and the map $\phi:A\rightarrow \mathbb M_2(\mathbb Z)$ is an injective ring homomorphism.

By the hypothesis we have that $A$ contains the ring generated by $SL_2(m\mathbb Z)$ inside $\mathbb M_2(\mathbb Z)$. Since the latter contains $m^2\mathbb M_2(\mathbb Z)$, we deduce that $m^2\mathbb M_2(\mathbb Z)\subset A$.

Next, we find a matrix $k\in\mathbb M_2(\mathbb Z)$ with non-zero determinant such that $\phi(x)k=kx$, for all $x\in A$. 
We view every $v\in\mathbb Z^2$ as a $2\times 1$ matrix over $\mathbb Z$ and denote by $v^t$ its transpose. Let  $v\in\mathbb Z^2$ and $w\in\mathbb (m^2\mathbb Z)^2$ be non-zero vectors. Since $vw^t\in m^2\mathbb M_2(\mathbb Z)\setminus\{0\}$ and  $\phi$ is injective we have that $\phi(vw^t)\not=0$. Let $z\in\mathbb Z^2$ such that $\phi(vw^t)z\not=0$. 

We define $k:\mathbb Z^2\rightarrow\mathbb Z^2$ by letting $k(\omega)=\phi(\omega w^t)z$.  Fix  $\omega\in\mathbb Z^2$ and $x\in A$. Since $\phi$ is multiplicative on $A$ we have that $\phi(x\omega w^t)=\phi(x)\phi(\omega w^t)$.
Thus,   $k(x\omega)=\phi(x\omega w^t)z=\phi(x)\phi(\omega w^t)z=\phi(x)k(\omega).$ Hence, if we view $k$ as an element of $\mathbb M_2(\mathbb Z)$, then $kx=\phi(x)k$, for all $x\in A$. This implies that $\ker(k)$ is an $A$-invariant subgroup of $\mathbb Z^2$. Since $A\supset m^2\mathbb M_2(\mathbb Z)$, we get that $\ker(k)$ is equal to either $\{0\}$ of $\mathbb Z^2$. Since $k(v)\not=0$,  it follows that $k$ is injective, proving our claim.

Now, write $k=\bigl(\begin{smallmatrix}a&b\\ c&d\end{smallmatrix} \bigr)$, where $a,b,c,d$ are integers. After replacing $k$ with $\frac{1}{N}k$, for some $N\geqslant 1$, we may assume that the greatest common divisor of $a,b,c,d$ is 1. We claim that $k\in GL_2(\mathbb Z)$.
Note that if $x\in SL_2(m\mathbb Z)$, then the hypothesis gives that $\phi(x)\in SL_2(n\mathbb Z)$. Since we also have that $x\in A$, we get that $kxk^{-1}=\phi(x)\in SL_2(n\mathbb Z)$. By applying this relation to $x\in\{ \bigl(\begin{smallmatrix}1&m\\ 0&1\end{smallmatrix} \bigr), \bigl(\begin{smallmatrix}1&0\\ m&1\end{smallmatrix}\bigr)\}$, it follows that $m k\bigl(\begin{smallmatrix}0&1\\ 0&0\end{smallmatrix} \bigr)k^{-1},m k\bigl(\begin{smallmatrix}0&0\\1&0\end{smallmatrix} \bigr)k^{-1}\in n\mathbb M_2(\mathbb Z)$. 

These relations imply that $n\det(k)$ divides $ma^2,mb^2,mc^2$ and $md^2$. Thus,  $n\det(k)\mid m$ and so $n\mid m$. By symmetry, we also get that $m\mid n$. Hence $m=n$ and $\det(k)\in\{\pm 1\}$, thus $k\in GL_2(\mathbb Z)$.

Finally, let $p\in S$. Then for every $x\in A$  we have $$g_p\pi_p(x)g_p^{-1}\pi_p(k)=\pi_p(\phi(x)k)=\pi_p(kx)=\pi_p(k)\pi_p(x).$$ This shows that $g_p^{-1}\pi_p(k)$ commutes with $\pi_p(A)$. By using again the fact that  $A\supset m^2\mathbb M_2(\mathbb Z)$ and  $p\nmid m$, we get that  $g_p^{-1}\pi_p(k)$ is a multiple of the identity. Therefore, we can find $l_p\in\mathbb Z/p\mathbb Z\setminus\{0\}$ such that $\pi_p(k)l_p=g_p I$. Since $l=(l_p)_{p\in S}$ belongs to the center of $G$, the conclusion follows. \hfill$\blacksquare$

\begin{lemma} \label{lemma3} Let $p$ be a prime and $S$ be a set of primes such that $p\not\in S$. 
\begin{enumerate}

\item\emph{\cite{GG88}} If $K<K_S$ is an open subgroup and $\delta:K\rightarrow SL_2(\mathbb Z_p)$ is a continuous homomorphism, then $\delta(K)$ is finite.
\item If $L<SL_2(\mathbb Z_p)$ is an open subgroup, then there is no injective continuous homomorphism $\delta:L\rightarrow K_S$.
\end{enumerate}
\end{lemma}

{\it Proof}. We denote by $\rho_q:SL_2(\mathbb Z_q)\rightarrow SL_2(\mathbb F_q)$, for $q\geqslant 3$, and by  $\rho_2:SL_2(\mathbb Z_2)\rightarrow SL_2(\mathbb Z/4\mathbb Z)$ the obvious surjective homomorphisms.  Then  $\Gamma_q:=\ker(\rho_q)$ is a pro-$q$ normal open subgroup of $SL_2(\mathbb Z_q)$, for any prime $q$. (Recall that a profinite group $G$ is {\it pro-$q$} if  the index $[G:G_0]$ is a power of $q$, for any open subgroup $G_0<G$). 

Moreover, as is well-known, $\Gamma_q$ is torsion free. Thus, if $\Delta$ is a finite subgroup of $SL_2(\mathbb Z_q)$, then ${\rho_q}_{|\Delta}$ is injective and therefore $|\Delta|$ divides $ [SL_2(\mathbb Z_q):\Gamma_q]$.

(1)  This part follows from the proof of \cite[Lemma A.6.]{GG88}. For completeness, we include the argument from \cite{GG88}. 
We first treat the case when $S$ has one element. Thus, let $\delta:K\rightarrow SL_2(\mathbb Z_p)$ be a continuous homomorphism, where $K$ is an open subgroup of $SL_2(\mathbb Z_q)$, for some prime $q\not=p$. Since we have $[\delta(K):\delta(K\cap\Gamma_q)]\leqslant [K:K\cap\Gamma_q]<\infty$, we may assume that $K<\Gamma_q$.

 Fix $n\geqslant 2$ and define $\Gamma_{p,n}=\ker(SL_2(\mathbb Z_p)\rightarrow SL_2(\mathbb Z/p^n\mathbb Z))$. Then $\delta^{-1}(\Gamma_{p,n})$ is an open subgroup of $\delta^{-1}(\Gamma_p)$. Let $m=[\delta^{-1}(\Gamma_p):\delta^{-1}(\Gamma_{p,n})]$. Since $\delta^{-1}(\Gamma_p)$ is a pro-$q$ group (being an open subgroup of $K$), we can write $m=q^{\alpha}$, for some $\alpha\geqslant 0$. On the other hand, since $m|[\Gamma_p:\Gamma_{p,n}]$ and $\Gamma_p$ is pro-$p$,  we have that $m=p^{\beta}$, for some $\beta\geqslant 0$. Since $p\not=q$, we must have that $m=1$. 

In other words, $\delta^{-1}(\Gamma_p)=\delta^{-1}(\Gamma_{n,p})$, for all $n\geqslant 2$. Since $\cap_{n\geqslant 2}\Gamma_{p,n}=\{e\}$, we get that $\delta^{-1}(\Gamma_p)=\{e\}$. This implies that $\delta(K)$ is finite.

Now, in general, let $K<K_S$ be an open subgroup and $\delta:K\rightarrow SL_2(\mathbb Z_p)$ be a continuous homomorphism. Then $K$ contains an open subgroup $L<K_S$ of the form $L=\prod_{q\in S}L_q$, such that $L_q<SL_2(\mathbb Z_q)$ is an open subgroup, for all $q\in S$, and the set $\{q\in S|\; L_q\not=SL_2(\mathbb Z_q)\}$ is finite. For every subset $F\subset S$, we denote $L_F=\prod_{q\in F}L_q<K_S$.

The first part of the proof gives that $\delta(L_q)$ is a finite subgroup of $SL_2(\mathbb Z_p)$, for all $q\in S$.
  Since $\{\delta(L_q)\}_{q\in S}$ are mutually commuting subgroups of $SL_2(\mathbb Z_p)$, we  get that $\delta(L_F)$ is a finite group, whenever $F\subset S$ is finite. This implies that $|\delta(L_F)|\leqslant N:=[SL_2(\mathbb Z_p):\Gamma_p]$. Thus, we can choose $F$ such that $|\delta(L_F)|\geqslant |\delta(L_{F'})|$, for any other finite subset  $F'$ of $S$. In particular, if $F'\supset F$ then $\delta(L_{F'})=\delta(L_{F})$. Since $\delta$ is  continuous, we deduce that $\delta(L)=\delta(L_F)$, hence $\delta(L)$ is finite. Since $L<K$ has finite index, we conclude that $\delta(K)$ is finite as well.

(2) Assume by contradiction that there is an injective continuous homomorphism $\delta:L\rightarrow K_S$, where $L<SL_2(\mathbb Z_p)$ is an open subgroup and $p\not\in S$. We may take $L<\Gamma_p$ so that $L$ is pro-$p$. 

We denote $\rho=(\rho_q)_{q\in S}:K_S\rightarrow\prod_{q\in S}SL_2(\mathbb Z_q)/\Gamma_q$ and claim that  $\rho\circ\delta$ is injective. 
To see this, for $q\in S$, we denote by $\sigma_q:K_S\rightarrow SL_2(\mathbb Z_q)$ the quotient homomorphism. Then $\sigma_q\circ\delta:L\rightarrow SL_2(\mathbb Z_q)$ is a continuous homomorphism. Since $p\not\in S$,  part (1) yields that $\sigma_q(\delta(L))$ is finite. Since $\Gamma_q$ is torsion free we derive that $\sigma_q(\delta(L))\cap\Gamma_q=\{e\}$, for all $q\in S$.

Now,  assume that $\rho(\delta(x))=e$, for some $x\in L$. Then, given $q\in S$, we have $\rho_q(\sigma_q(\delta(x)))=e$ or, equivalently, $\sigma_q(\delta(x))\in\Gamma_q$. Using  the last paragraph, we get that $\sigma_q(\delta(x))=e$, hence $\delta(x)=e$. Since $\delta$ is injective, we deduce that $x=e$. Hence, $\rho\circ\delta$ is also injective. 

Next, we note that $\rho_q\circ\sigma_q\circ\delta:L\rightarrow SL_2(\mathbb Z_q)/\Gamma_q$ is not onto, for any $q\in S$. This is because $L$ is  pro-$p$  while the cardinality of $SL_2(\mathbb Z_q)/\Gamma_q$ is equal to $q(q^2-1)$, if $q\geqslant 3$, and to $40$, if $q=2$.

To finish the proof we use an idea from \cite{Ga02}.
Let $L_q=\rho_q(\sigma_q(\delta(L)))$. Then $L_q$ is a proper subgroup of $SL_2(\mathbb Z_q)/\Gamma_q\cong SL_2(\mathbb Z_q)$, when $q\geqslant 3$, and of $SL_2(\mathbb Z_2)/\Gamma_2\cong SL_2(\mathbb Z/4\mathbb Z),$ when $q=2$. By Proposition \ref{sl_2}, part (4), we have that either $|L_q|\leqslant 60$ or $[[a,b],[c,d]]=e$, for all $a,b,c,d\in L_q$.
In either case we deduce that $[[a,b],[c,d]]^{60!}=e$, for all $a,b,c,d\in L_q$ and for every $q\in S$. Since $\rho(\delta(L))\subset\prod_{q\in S} L_q$ and $\rho\circ\delta$ is injective, we conclude that  $[[a,b],[c,d]]^{60!}=e$, for all $a,b,c,d\in L.$ This is however impossible since $L$ contains a non-abelian free group (e.g. $L\cap SL_2(\mathbb Z)$).
\hfill$\blacksquare$

\begin{corollary} \label{cor3} Let $S$ and $T$ be two sets of primes. Let $K<K_S$ be an open subgroup and $\Gamma<SL_2(\mathbb Z)$ be an infinite subgroup. 

If there exists an injective continuous homomorphism $\delta:K\rightarrow K_T$ such that $\delta(\Gamma\cap K)\subset SL_2(\mathbb Z)$, then $S=T$.

\end{corollary}
{\it Proof.} Firstly, Lemma \ref{lemma3} (2) implies that $S\subset T$. Secondly, to show that $T\subset S$, assume by contradiction that we can find $p\in T\setminus S$ and let $\rho:K_T\rightarrow SL_2(\mathbb Z_p)$ be the quotient homomorphism. By Lemma \ref{lemma3} (1) we get that $\rho(\delta(K))$ is finite. On the other hand, by the hypothesis, $\Lambda=\delta(\Gamma\cap K)$ is an infinite subgroup of $SL_2(\mathbb Z)$. Since the restriction of $\rho$ to $SL_2(\mathbb Z)$ is injective, we deduce that $\rho(\Lambda)$ is infinite, which gives a contradiction.
\hfill$\blacksquare$
\vskip 0.05in

We are now ready to prove corollaries \ref{corB} and \ref{corC}.

{\bf Proof of Corollary \ref{corB}.}  By Fact \ref{strong}, $G_{\Gamma,S}<G_S$ is an open subgroup and the action $\Gamma\curvearrowright G_{\Gamma,S}$ has spectral gap, for any infinite set of primes $S$ and every non-amenable subgroup $\Gamma<SL_2(\mathbb Z)$.

Let $S$ and $T$ be two infinite sets of primes. Below, we prove assertions (1) and (2) separately.

(1) Assume that either $SL_2(\mathbb Z)\curvearrowright G_S$ is stably orbit equivalent to $SL_2(\mathbb Z)\curvearrowright G_T$ , or $\mathcal R(\Gamma\curvearrowright G_{\Gamma,S})$ is Borel reducible to $\mathcal R(\Lambda\curvearrowright G_{\Lambda,T})$. In either case, by applying Theorem \ref{mainthm}, we can find an open subgroup $G_0<G_S$, a closed subgroup $G_1<G_T$, and a continuous isomorphism $\delta:G_0\rightarrow G_1$ such that $\delta(\pi_S(SL_2(\mathbb Z))\cap G_0)=\pi_T(SL_2(\mathbb Z))\cap G_1$.

Let $S_1=S\setminus\{2,3,5\}$ and $T_0=T\setminus\{2,3,5\}$. Then $G_{S_1}\cap \delta^{-1}(G_1\cap G_{T_0})$ is an open subgroup of $G_0$ and thus of $G_S$. 
Hence, we can find a subset $S_0\subset S_1$ such that $S\setminus S_0$ is finite, $G_{S_0}\subset G_0$, and $\delta(G_{S_0})\subset G_{T_0}.$
Lemma \ref{lemma1} implies that $S_0\subset T_0$, $\delta(G_{S_0})=G_{S_0}$ and that there is $g\in H_{S_0}$ such that $\delta(x)=gxg^{-1}$, for all $x\in G_{S_0}$. Moreover, we have that $\delta(\pi_S(SL_2(\mathbb Z))\cap G_{S_0}))=\pi_T(SL_2(\mathbb Z))\cap G_{S_0}$.

 Next, we denote by $m$ the product of the primes in $S\setminus S_0$. Then $\pi_S(SL_2(\mathbb Z))\cap G_{S_0}=\pi_{S_0}(SL_2(m\mathbb Z))$. Since $\pi_T(SL_2(\mathbb Z))\cap G_{S_0}$ is infinite (hence, non-trivial) it follows that $T\setminus S_0$ is finite. Moreover, if $n$ is the product of the primes in $T\setminus S_0$, then $\pi_T(SL_2(\mathbb Z))\cap G_{S_0}=\pi_{S_0}(SL_2(n\mathbb Z))$.
Altogether, we get that $g\pi_{S_0}(SL_2(m\mathbb Z))g^{-1}=\pi_{S_0}(SL_2(n\mathbb Z))$. 

Lemma \ref{lema2} implies that $m=n$ and thus $S\setminus S_0=T\setminus S_0$. This gives that $S=T$, as claimed.

(2) Let $\Gamma,\Lambda<SL_2(\mathbb Z)$ be two non-amenable subgroups. Assume that either $\Gamma\curvearrowright G_{\Gamma,S}$ is stably orbit equivalent to $\Lambda\curvearrowright G_{\Lambda,T}$, or $\mathcal R(\Gamma\curvearrowright G_{\Gamma,S})$ is Borel reducible to $\mathcal R(\Lambda\curvearrowright G_{\Lambda,T})$. In either case, Theorem \ref{mainthm} implies that we can find an open subgroup $G_0<G_{\Gamma,S}$, a closed subgroup $G_1<G_{\Lambda,T}$, and a continuous isomorphism $\delta:G_0\rightarrow G_1$ such that $\delta(\Gamma\cap G_0)=\Lambda\cap G_1$.

Since $G_{\Gamma,S}<G_S$ is open,  $G_0<G_S$ is open. By Corollary \ref{cor2}, we conclude that $|S\Delta T|<\infty$. \hfill$\blacksquare$

 {\bf Proof of Corollary \ref{corC}.} By Fact \ref{strong}, $K_{\Gamma,S}<K_S$ is an open subgroup and the action $\Gamma\curvearrowright K_{\Gamma,S}$ has spectral gap, for any  nonempty set of primes $S$ and every non-amenable subgroup $\Gamma<SL_2(\mathbb Z)$.

Let $\Gamma,\Lambda<SL_2(\mathbb Z)$ be non-amenable subgroups and $S,T$ be nonempty set of primes.
Assume that either $\Gamma\curvearrowright K_{\Gamma,S}$ is stably orbit equivalent to $\Lambda\curvearrowright K_{\Lambda,T}$, or $\mathcal R(\Gamma\curvearrowright K_{\Gamma,S})$ is Borel reducible to $\mathcal R(\Lambda\curvearrowright K_{\Lambda,T})$. Theorem \ref{mainthm} then implies that we can find an open subgroup $K_0<K_{\Gamma,S}$, a closed subgroup $K_1<K_{\Lambda,T}$, and a continuous isomorphism $\delta:K_0\rightarrow K_1$ such that $\delta(\Gamma\cap K_0)=\Lambda\cap K_1$.

Since $K_{\Gamma,S}<K_S$ is open, $K_0<K_S$ is open, and Corollary \ref{cor3} implies that $S=T$. \hfill$\blacksquare$

We end this section by proving the following strengthening of part of Corollary \ref{corC}.

\begin{corollary}\label{nohom} Let $S$, $T$ be nonempty sets of primes, and $\Gamma,\Lambda<SL_2(\mathbb Z)$ non-amenable subgroups.

 If  $T\not\subset S$, then any homomorphism from $\mathcal R(\Gamma\curvearrowright K_{\Gamma,S})$ to $\mathcal R(\Lambda\curvearrowright K_{\Lambda,T})$ is trivial.
 \end{corollary}

{\it Proof.}
Assume that there is a non-trivial homomorphism from $\mathcal R(\Gamma\curvearrowright K_{\Gamma,S})$ to $\mathcal R(\Lambda\curvearrowright K_{\Lambda,T})$. By Corollary \ref{bborel} we can find an open subgroup $K_0<K_{\Gamma,S}$ and a continuous homomorphism $\delta:K_0\rightarrow K_T$ such that $\delta(\Gamma\cap K_0)\subset\Lambda$ and $\delta(K_0)\not\subset\Lambda$.
Also, note that Fact \ref{strong} implies that $K_0<K_S$ is an open subgroup.

Let $p\in T\setminus S$ and $\rho:K_T\rightarrow SL_2(\mathbb Z_p)$ be the quotient homomorphism. By Lemma \ref{lemma3} (1) we have that $\rho(\delta(K_0))$ is finite. Thus, we can find an open subgroup $K_1<K_0$ so that $\rho(\delta(K_1))=\{e\}$. Since the restriction of $\rho$ to $SL_2(\mathbb Z)$ is injective and $\delta(\Gamma\cap K_1)\subset SL_2(\mathbb Z)$, we deduce that $\delta(\Gamma\cap K_1)=\{e\}$. Since $\Gamma\cap K_1<K_1$ is dense, it follows that $\delta(K_1)=\{e\}$. Since $\Gamma\cap K_0<K_0$ is dense, we also have that $(\Gamma\cap K_0) K_1=K_0$. From this we conclude that $\delta(K_0)=\delta(\Gamma\cap K_0)\subset\Lambda$. This gives a contradiction and finishes the proof of the corollary.
\hfill$\blacksquare$

\section{Proof of Corollary \ref{corD}} \label{padic}
This section is devoted to the proof of Corollary \ref{corD}. 
In fact, we will establish a stronger result from which Corollary \ref{corD} will follow easily. Before stating this result, let us recall some notation.

For a prime $p$, let $PG(1,\mathbb Q_p)=\mathbb Q_p\cup\{\infty\}$ be the projective line over the field $\mathbb Q_p$ of $p$-adic numbers. Consider the Borel action of $GL_2(\mathbb Q_p)$ on $PG(1,\mathbb Q_p)$ by linear fractional transformations:

$$\begin{pmatrix} a&b\\c&d\end{pmatrix}\cdot x=\frac{ax+b}{cx+d}.$$

\begin{theorem}\label{Qp} Let $p\not=q$ be primes, and $\Gamma<SL_2(\mathbb Z)$, $\Lambda<GL_2(\mathbb Q_q)$ be countable subgroups.
Assume that $\Gamma$ is non-amenable and denote by $G$ the closure of $\Gamma$ in $SL_2(\mathbb Z_p)$.

Then  any homomorphism from $\mathcal R(\Gamma\curvearrowright G)$ to $\mathcal R(\Lambda\curvearrowright PG(1,\mathbb Q_q))$ is trivial with respect to $m_G$.
\end{theorem}

{\it Proof.} 
Since $\Gamma<SL_2(\mathbb Z)$ is non-amenable, Fact \ref{strong} gives that $G<SL_2(\mathbb Z_p)$ is an open subgroup and the action $\Gamma\curvearrowright (G,m_G)$ has spectral gap.

We denote $PGL_2(\mathbb Q_q)=GL_2(\mathbb Q_q)/Z_q$, where $Z_q=\{zI|z\in\mathbb Q_q^*\}$ is the center of $GL_2(\mathbb Q_q)$. Let $\rho:GL_2(\mathbb Q_q)\rightarrow PGL_2(\mathbb Q_q)$ be the quotient homomorphism. Notice $Z_q$ acts trivially on $PG(1,\mathbb Q_q)$ and consider the resulting action $PGL_2(\mathbb Q_q)\curvearrowright PG(1,\mathbb Q_q)$.

Moreover, if we let $\Omega:=\{(x_1,x_2,x_3)\in PG(1,\mathbb Q_q)^3|\;x_i\not=x_j,\;\text{for all}\;1\leqslant i<j\leqslant 3\}$, then we have

{\bf Claim 1.}  The action $PGL_2(\mathbb Q_q)\curvearrowright\Omega$ is transitive and free.

{\it Proof of Claim 1}. This claim is well-known, but we include a proof for completeness.  

Let $x_1,x_2,x_3\in PG(1,\mathbb Q_q)$ be distinct points. Consider the matrix $g=\begin{pmatrix} x_2-x_3&-x_1(x_2-x_3)\\x_2-x_1&-x_3(x_2-x_1)\end{pmatrix}$. Then $g\in GL_2(\mathbb Q_q)$ and the function $x\rightarrow g\cdot x=\frac{(x-x_1)(x_2-x_3)}{(x-x_3)(x_2-x_1)}$ maps $x_1\mapsto 0,x_2\mapsto 1$ and $x_3\mapsto\infty$. This shows that the action $PGL_2(\mathbb Q_q)\curvearrowright\Omega$ is transitive.

 If $g=\begin{pmatrix}a&b\\c&d\end{pmatrix}\in GL_2(\mathbb Q_q)$ and the equation $g\cdot x=\frac{ax+b}{cx+d}=x$ has at least three distinct solutions $x\in PG(1,\mathbb Q_q)$, then  $a=d$ and $b=c=0$, hence $g\in Z_q$. This shows that the action $PGL_2(\mathbb Q_q)\curvearrowright\Omega$ is also free. 
\hfill$\square$

Now, since $Z_q$ acts trivially on $PG(1,\mathbb Q_q)$,  it follows that $\mathcal R(\Lambda\curvearrowright PG(1,\mathbb Q_q))$ is Borel isomorphic to $\mathcal R(\rho(\Lambda)\curvearrowright PG(1,\mathbb Q_q)$.  Let $\theta:G\rightarrow PG(1,\mathbb Q_q)$ be a homomorphism from $\mathcal R(\Gamma\curvearrowright G)$ to $\mathcal R(\rho(\Lambda)\curvearrowright PG(1,\mathbb Q_q))$. Assume  by contradiction that $\theta$ is not trivial.

{\bf Claim 2.} $(\theta(x_1),\theta(x_2),\theta(x_3))\in\Omega$, for almost every $(x_1,x_2,x_3)\in G^3$.

{\it Proof of Claim 2.} If the claim is false, then the set $\{(x_1,x_2)\in G^2|\theta(x_1)=\theta(x_2)\}$ has positive measure. 
This would then imply that there exists $y\in PG(1,\mathbb Q_q)$ such that the set $\{x\in G|\theta(x)=y\}$ has positive measure. Since $\Gamma<G$ is dense, we would further get that $\theta(x)\in\Lambda y$, for almost every $x\in G$. This contradicts our assumption that $\theta$ is not trivial. \hfill$\square$

Since the action $\Gamma\curvearrowright (G,m_G)$ has spectral gap and any transitive action is smooth,  Claims 1 and 2 imply that the hypothesis of Theorem \ref{homrig2}, part (1), is verified. 
Thus, we can find an open subgroup $G_0<G$, a continuous homomorphism $\delta:G_0\rightarrow PGL_2(\mathbb Q_q)$, a Borel map $\phi:G_0\rightarrow \rho(\Lambda)$, and $y\in PG(1,\mathbb Q_q)$  such that $\delta(\Gamma\cap G_0)\subset \rho(\Lambda)$ and $\theta(x)=\phi(x)\delta(x)y$, for almost every $x\in G_0$.

{\bf Claim 3.} We can find an open subgroup $G_1<G_0$ such that $\delta(G_1)=\{e\}$.

{\it Proof of Claim 3.} Let $K_q=I+q\mathbb M_2(\mathbb Z_q)=\{\bigl(\begin{smallmatrix}1+a&b\\c&1+d\end{smallmatrix}\bigr)|\;a,b,c,d\in q\mathbb Z_q\}$. Then $K_q<GL_2(\mathbb Q_q)$ is an open compact subgroup, hence $\rho(K_q)<PGL_2(\mathbb Q_q)$ is an open subgroup. Since $\delta$ is continuous and $G_0<SL_2(\mathbb Z_p)$ is an open subgroup, we can find an open subgroup $G_1<G_0$ such that $\delta(G_1)\subset\rho(K_q)$. Moreover, we may assume that $G_1$ is a pro-$p$ group (see the proof of Lemma \ref{lemma3}).

Now, let $K_{q,n}=I+q^n\mathbb M_2(\mathbb Z_q)$. Then $K_{q,n}$ is an open subgroup of $K_q$, for every $n\geqslant 1$, and $\{K_{q,n}\}_{n\geqslant 1}$ is a basis of open neighborhoods of the identity in $K_q$. 
Since $K_q/K_{q,n}\cong I+q\mathbb M_2(\mathbb Z/q^n\mathbb Z)$ and $| I+q\mathbb M_2(\mathbb Z/q^n\mathbb Z)|=q^{4(n-1)}$, we deduce that $K_q$ is a pro-$q$ group. This implies that $\rho(K_q)$, being a continuous image of $K_q$, is also a pro-$q$ group.

Since $G_1$ is  pro-$p$, $K_q$ is pro-$q$ and $p\not=q$, any continuous homomorphism from $G_1$ to $K_q$ is trivial (see e.g. the proof of Lemma \ref{lemma3}). This proves the claim.\hfill$\square$

Finally, Claim 3 implies that $\theta(x)\in\Lambda y$, for almost every $x\in G_1$. Since $\Gamma<G$ is dense and $G_1<G$ is open, we have that $\Gamma G_1=G$. From this we derive that $\theta(x)\in\Lambda y$, for almost every $x\in G$, and hence that $\theta$ is trivial.
\hfill$\blacksquare$

{\bf Proof of Corollary \ref{corD}.} Let $\Gamma<GL_2(\mathbb Q_p)$ and $\Lambda<GL_2(\mathbb Q_q)$ be countable subgroups, for some primes $p\not=q$, such that $\Gamma_0:=\Gamma\cap SL_2(\mathbb Z)$ is non-amenable. Assume by contradiction that there exists a Borel reduction $\theta:PG(1,\mathbb Q_p)\rightarrow PG(1,\mathbb Q_q)$ from $\mathcal R(\Gamma\curvearrowright PG(1,\mathbb Q_p))$ to $\mathcal R(\Lambda\curvearrowright PG(1,\mathbb Q_q))$.

Next, note that the action $SL_2(\mathbb Z_p)\curvearrowright PG(1,\mathbb Q_p)$ is Borel isomorphic to the left multiplication action $SL_2(\mathbb Z_p)\curvearrowright SL_2(\mathbb Z_p)/K_p$, where $K_p<SL_2(\mathbb Z_p)$ denotes the subgroup of lower triangular matrices. Indeed, the former action  is transitive (see e.g. \cite[Lemma 6.1]{Th01}) and $K_p$ is the stabilizer of $0\in PG(1,\mathbb Q_p)$. This implies that there exists a unique $SL_2(\mathbb Z_p)$-invariant Borel probability measure on $PG(1,\mathbb Q_p)$, which we denote by $\mu_p$.

Let $\pi:SL_2(\mathbb Z_p)\rightarrow SL_2(\mathbb Z_p)/K_p$ be the quotient map, and define $\Theta:=\theta\circ\pi:SL_2(\mathbb Z_p)\rightarrow PG(1,\mathbb Q_q)$. Denote also by $G$ the closure of $\Gamma_0$ in $SL_2(\mathbb Z_p)$. Since $\Gamma_0$ is non-amenable, Fact \ref{strong} implies that $G<SL_2(\mathbb Z_p)$ is an open subgroup. Let $k=[SL_2(\mathbb Z_p):G]$ and $g_1,...,g_k\in SL_2(\mathbb Z_p)$ such that $SL_2(\mathbb Z_p)=\sqcup_{i=1}^kGg_i$. For $1\leqslant  i\leqslant k$, let
 $\Theta_i:G\rightarrow PG(1,\mathbb Q_p)$ be given by $\Theta_i(x)=\Theta(xg_i)$. 
 
 Then $\Theta_i$ is a homomorphism from $\mathcal R(\Gamma_0\curvearrowright G))$ to $\mathcal R(\Lambda\curvearrowright PG(1,\mathbb Q_q))$. By Theorem \ref{Qp} we conclude that $\Theta_i$ is trivial, i.e. there is $y_i\in PG(1,\mathbb Q_q)$ such that $\Theta_i(x)\in\Lambda y_i$, for $m_G$-almost every $x\in G$. From this it follows that $\theta(x)\in\cup_{i=1}^k\Lambda y_i$, for $\mu_p$-almost every $x\in PG(1,\mathbb Q_p)$.  This contradicts the fact that $\theta$ is countable-to-one. \hfill$\blacksquare$

\section{Proofs of Corollaries \ref{corE} and \ref{corF} } 
In this section we will use Corollary \ref{cor1} to deduce Corollaries \ref{corE} and \ref{corF}.
We start by proving Corollary \ref{corE}.
As in the beginning of Section \ref{sect4}, we denote $G_S=\prod_{p\in S}SL_2(\mathbb F_p)$ and $H_S=\prod_{p\in S}GL_2(\mathbb F_p)$, for any infinite set of primes $S$.
Let $\pi_p:\mathbb M_2(\mathbb Z)\rightarrow \mathbb M_2(\mathbb F_p)$ be the reduction modulo $p$.

Let $\alpha=\bigl(\begin{smallmatrix}1&0\\ 0&-1\end{smallmatrix} \bigr)\in GL_2(\mathbb Z).$ We denote still by $\alpha$ the element $(\pi_p(\alpha))_{p\in S}$ of $H_S$.
Then $\alpha^2=I$ and $\alpha\in H_S$ normalizes $G_S$. 
Below we consider the semidirect product $G_S\rtimes\mathbb Z/2\mathbb Z$, where $\mathbb Z/2\mathbb Z=\{I,\alpha\}$ and $\alpha$ acts on $G_S$ by conjugation.
Note also that $\alpha$ commutes with the center $Z=\{\pm I\}$ of $G_S$.

We begin by proving the first part of Corollary \ref{corE}.
\begin{corollary} Let $S$ be an infinite set of primes and denote $\mathcal R=\mathcal R(SL_2(\mathbb Z)\curvearrowright G_S)$.

Then Out$(\mathcal R)\cong (G_S/Z)\rtimes\mathbb Z/2\mathbb Z.$
\end{corollary}

{\it Proof.} We define $\sigma_h:G_S\rightarrow G_S$, for $h\in G_S$, and $\sigma_{\beta}:G_S\rightarrow G_S$, for $\beta\in\{I,\alpha\}$, by letting

 $$\sigma_h(x)=xh\;\;\;\text{and}\;\;\;\sigma_{\beta}(x)=\beta x\beta^{-1},\;\;\;\text{for all $x\in G$.}\;\;\;$$

As $\sigma_h$ commutes with the action $SL_2(\mathbb Z)\curvearrowright G_S$ and $\beta$ normalizes $SL_2(\mathbb Z)$, we have $\sigma_h,\sigma_{\beta}\in$ Aut$(\mathcal R)$. Moreover, since $\sigma_{\beta}\sigma_{h}\sigma_{\beta}^{-1}=\sigma_{\beta h\beta^{-1}}$, we have  a homomorphism $\sigma:G_S\rtimes\mathbb Z/2\mathbb Z\rightarrow$ Aut$(\mathcal R)$. We denote  $\rho=\varepsilon\circ\sigma:G_S\rtimes\mathbb Z/2\mathbb Z\rightarrow$ Out$(\mathcal R)$, where  $\varepsilon:$ Aut$(\mathcal R)\rightarrow$ Out$(\mathcal R)$  is the quotient homomorphism. 

To prove the corollary, we only have to argue that $\ker(\rho)=Z$ and that $\rho$ is surjective. Firstly, let $h\in G_S$ and $\beta\in\{I,\alpha\}$ such that $\sigma_{(h,\beta)}=\sigma_{h}\sigma_{\beta}$ belongs to $[\mathcal R]$. Thus, we can find $\gamma\in SL_2(\mathbb Z)$ and a Borel set $A\subset G_S$ of positive (Haar) measure such that \begin{equation}\label{beta}\beta x\beta^{-1}h=\gamma x,\;\;\;\text{for all $x\in A$}.\end{equation}

This implies that $\gamma^{-1}\beta$ commutes with $xy^{-1}$, for all $x,y\in A$. Since $AA^{-1}$ contains an open neighborhood of the identity in $G_S$ we deduce that $\gamma^{-1}\beta$ commutes with $G_T$, for some subset $T\subset S$ with $S\setminus T$ finite. Hence, for all $p\in T$, we have that $\pi_p(\gamma^{-1}\beta)\in\{\pm I\}$.
Since $T$ is infinite and $\gamma^{-1}\beta\in GL_2(\mathbb Z)$, it follows easily that $\gamma^{-1}\beta= {\pm I}$. Thus, we must have that $\beta=I$ and $\gamma=\pm I$. By using equation \ref{beta}, we derive that $h=\pm I$, showing that $(h,\beta)\in Z$.

To show that $\rho$ is surjective, let $\theta\in$ Aut$(\mathcal R)$.
 Then by Corollary \ref{cor1}, after replacing $\theta$ with $\theta\circ\tau$, for some $\tau\in [\mathcal R]$, we can find open subgroups $G_0,G_1<G_S$, an isomorphism $\delta:G_0\rightarrow G_1$ and $h\in G_S$ such that  $\delta(SL_2(\mathbb Z)\cap G_0)=SL_2(\mathbb Z)\cap G_1$ and $\theta(x)=\delta(x)h$, for all $x\in G_1$.

Let $S_0\subset S$ be a subset such that $S_0\subset S\setminus\{2,3,7\}$, $S\setminus S_0$ is finite, $G_{S_0}\subset G_0$ and $\delta(G_{S_0})\subset G_1\cap G_{S_0}$. 
By Lemma \ref{lemma1} we get that $\delta(G_{S_0})=G_{S_0}$ and that there is $g\in H_{S_0}$ such that $\delta(x)=gxg^{-1}$, for all $x\in G_{S_0}$. Since $\delta(SL_2(\mathbb Z)\cap G_{S_0})=SL_2(\mathbb Z)\cap G_{S_0}$ by applying Lemma \ref{lema2} we deduce that there is $k\in GL_2(\mathbb Z)$ such that $\delta(x)=kxk^{-1}$, for all $x\in G_{S_0}$. 

Let $\beta\in\{I,\alpha\}$ such that $k=\gamma\beta$, for some $\gamma\in SL_2(\mathbb Z)$. Then for almost every $x\in G_{S_0}$ we have that $\theta(x)=\delta(x)h=\gamma\beta x\beta^{-1}\gamma^{-1}h=\gamma\sigma_{(\gamma^{-1}h,\beta)}(x)$. Thus,  the set $A$ of $x\in G_S$ such that $\theta(x)\in\Gamma\sigma_{\gamma^{-1}h,\beta}(x)$ contains $G_{S_0}$. Since $A$  is invariant under $\mathcal R$, and $\mathcal R$ is ergodic, we deduce that $A=G_S$, almost everywhere.

 This implies that $\varepsilon(\theta)=\varepsilon(\sigma_{(\gamma^{-1}h,\beta}))=\rho_{(\gamma^{-1}h,\beta)}$, which proves that $\rho$ is surjective.
\hfill$\blacksquare$

We continue by establishing the second part of Corollary \ref{corE}.

\begin{corollary} Let $S$ be an infinite set of primes and $\Gamma<SL_2(\mathbb Z)$ be a non-amenable subgroup. Let $G_{\Gamma,S}$ denote the closure of $\Gamma$ in $G_S$.

Then the equivalence relation $\mathcal R=\mathcal R(\Gamma\curvearrowright G_{\Gamma,S})$ and the II$_1$ factor $M=L^{\infty}(G_{\Gamma,S})\rtimes\Gamma$ have trivial fundamental groups, i.e. $\mathcal F(\mathcal R)=\mathcal F(M)=\{1\}$.
\end{corollary}
{\it Proof.} Fact \ref{strong} gives that $G_{\Gamma,S}<G_S$ is an open subgroup and the action $\Gamma\curvearrowright G_{\Gamma,S}$ has spectral gap. Let $t\in\mathcal F(\mathcal R)$. Applying Corollary \ref{cor1} provides  open subgroups $G_0,G_1<G_{\Gamma,S}$ such that $t=\frac{[G_{\Gamma,S}:G_0]}{[G_{\Gamma,S}:G_1]}$ and a continuous isomorphism $\delta:G_0\rightarrow G_1$. 

Since $G_{\Gamma,S}<G_S$ is open, we get that $G_0<G_S$ is open. Thus, we may find a subset $S_0$ of $S$ such that $S_0\subset S\setminus\{2,3,5\}$, $S\setminus S_0$ is finite, $G_{S_0}\subset G_0$ and $\delta(G_{S_0})\subset G_1\cap G_{S_0}$. 
But then applying Lemma \ref{lemma1} to the injective homomorphism $\delta_{|G_{S_0}}:G_{S_0}\rightarrow G_{S_0}$ gives that $\delta(G_{S_0})=G_{S_0}$. Hence, we get that $$t=\frac{[G_{\Gamma,S}:G_0]}{[G_{\Gamma,S}:G_1]}=\frac{[G_1:G_{S_0}]}{[G_0:G_{S_0}]}=\frac{[\delta(G_0):\delta(G_{S_0})]}{[G_0:G_{S_0}]}=1.$$

Finally,  by a result of N. Ozawa and S. Popa \cite[Corollary 3]{OP07},  $M$ has a unique Cartan subalgebra, up to unitary conjugacy. This implies that $\mathcal F(M)=\mathcal F(\mathcal R)=\{1\}$.
\hfill$\blacksquare$

\begin{remark}
Note that in the above proof one can use \cite[Theorem 1.1 and Remark 4.1]{Io11} instead of \cite[Corollary 3]{OP07} to conclude that $\mathcal F(M)=\mathcal F(\mathcal R)$.
\end{remark}

We end this section by proving Corollary \ref{corF}.

{\bf Proof of Corollary \ref{corF}}. We prove separately assertions (1) and (2).

(1) Let $p$ be a prime and denote $\Gamma=SL_2(\mathbb Z), G=SL_2(\mathbb Z_p)$, $\mathcal  R=\mathcal R(\Gamma\curvearrowright G)$.
Our goal is to show that $Out(\mathcal R)$ is isomorphic to a $(\mathbb Z/2\mathbb Z)^2$-extension of $PSL_2(\mathbb Q_p)$.
To this end, let $K=SL_2(\mathbb Q_p)$. 
Then $K$ is a unimodular locally compact group, i.e. it admits a Haar measure $m_K$ which is invariant under both left and right multiplication with elements of $K$. Moreover, we have that: 

{\bf Claim 1.} $m_K$ is invariant under the conjugation action of $GL_2(\mathbb Q_p)$ on $K$. 

{\it Proof of Claim 1}. To see this, recall that if we parametrize (a co-null subset of) $SL_2(\mathbb Q_p)$ by $\{\bigl(\begin{smallmatrix}x&y\\z&\frac{yz+1}{x}\end{smallmatrix}\bigr)|\;x\in\mathbb Q_p\setminus\{0\},y,z\in\mathbb Q_p\}$, then, up to a multiplicative factor, $m_K$ is given by the differential form $\frac{1}{x}dx\wedge dy\wedge dz$. 

Now, let $a\in\mathbb Q_p\setminus\{0\}$ and put $\zeta=\bigl(\begin{smallmatrix}a&0\\0&1\end{smallmatrix}\bigr)$. Since $\zeta\bigl(\begin{smallmatrix}x&y\\z&t\end{smallmatrix}\bigr)\zeta^{-1}=\bigl(\begin{smallmatrix}x&ay\\\frac{z}{a}&t\end{smallmatrix}\bigr)$, we deduce that $m_K$ is invariant under conjugation with $\zeta$. Since $m_K$ is also invariant under conjugation with elements from $SL_2(\mathbb Q_p)$ and every $\eta\in GL_2(\mathbb Q_p)$ can be written as $\eta=\bigl(\begin{smallmatrix}\det(\eta)&0\\0&1\end{smallmatrix}\bigr)\eta'$, with $\eta'\in SL_2(\mathbb Q_p)$, the claim follows. \hfill$\square$

Next, we let $$H=\{h\in GL_2(\mathbb Q_p)|\det(h)=\pm p^n, \;\text{for some $n\in\mathbb Z$}\}\;\;\;\;\;\text{and}\;\;\;\;\;\Lambda=H\cap GL_2(\mathbb Z[\frac{1}{p}]).$$ Following S. Gefter \cite[Remark 2.8]{Ge96} and A. Furman \cite[Proof of Theorem 1.6]{Fu03}, we will define a homomorphism $\rho:H\rightarrow Out(\mathcal R)$. Fix $h\in H$.
We claim that there are $\lambda\in\Lambda$ and $g\in G$ such that $h=g^{-1}\lambda$. Note that $\lambda_0=\bigl(\begin{smallmatrix}\det(h)&0\\0&1\end{smallmatrix}\bigr)\in\Lambda$ and $h\lambda_0^{-1}\in K=SL_2(\mathbb Q_p)$. Since $SL_2(\mathbb Z[\frac{1}{p}]$ and $G$ are dense and, respectively, open in $K$, we can find $\lambda_1\in SL_2(\mathbb Z[\frac{1}{p}])$ and $g\in G$ such that $\lambda_0^{-1}h=g^{-1}\lambda_1$. Since $\lambda_1\lambda_0\in\Lambda$, the claim is proven.

Let $G'=\lambda^{-1} G\lambda\cap G$ and define $\sigma_h:G'\rightarrow G$ by letting $\sigma_h(x)=\lambda xh^{-1}=(\lambda x\lambda^{-1})g$.   Then $G'<G$ is an open subgroup. Moreover, Claim 1 implies that $m_G(\sigma_h(A))=m_G(A)$, for any Borel subset $A\subset G'$. Furthermore, since $\Lambda\cap G=\Gamma$, if $x,y\in G$, then $\Gamma x=\Gamma y$ if and only if $\Lambda x=\Lambda y$. Also, if $x,y\in G'$ then $\Lambda x=\Lambda y$ if and only if $\Lambda\sigma_h(x)=\Lambda\sigma_h(y)$. Altogether, we get that if $x,y\in G'$, then $$\Gamma x=\Gamma y\;\Longleftrightarrow\; \Gamma\sigma_h(x)=\Gamma\sigma_h(y)$$

This implies that $\sigma_h$ extends to an automorphism of $\mathcal R$, which we still denote by $\sigma_h$.

We define $\rho_h=\varepsilon(\sigma_h)$, where $\varepsilon:Aut(\mathcal R)\rightarrow Out(\mathcal R)$ is the quotient homomorphism.  Note that $\rho_h$ only depends on $h$, and not on the choices made in its definition. 
Furthermore, it is easy to see that $\rho:H\rightarrow Out(\mathcal R)$ is a homomorphism.

 The rest of the proof is divided between two claims.

{\bf Claim 2.} $\ker(\rho)=Z:=\{\pm p^nI|n\in\mathbb Z\}\subset H$. 

{\it Proof of Claim 2.} Let $h\in \ker(\rho)$.  Let $\lambda\in\Lambda$ and $g\in G$  such that $h=g^{-1}\lambda$. Since $\sigma_h\in [\mathcal R]$, we can find $\gamma\in\Gamma$ such that $A=\{x\in G|\gamma x=\lambda xh^{-1}\}$ has positive measure.
Notice that $\lambda^{-1}\gamma$ commutes with $xy^{-1}$, for all $x,y\in A$. Since $A$ has positive measure, $AA^{-1}$ contains an open subgroup of $G$. We get that $\lambda^{-1}\gamma\in\Lambda$ commutes with $\bigl(\begin{smallmatrix}1&p^n\\ 0&1\end{smallmatrix}\bigr)$ and $\bigl(\begin{smallmatrix}1&0\\ p^n&1\end{smallmatrix}\bigr)$, for some $n\geqslant 1$.  
From this it follows that $\lambda^{-1}\gamma=cI$, for some $c\in\mathbb Z[\frac{1}{p}]$. Since $\det(\lambda^{-1}\gamma)\in\{\pm p^n|n\in\mathbb Z\}$, we get that $c=\pm p^m$, for some $m\in\mathbb Z$. 
Since there is $x\in G$ such that $\gamma x=\lambda xh^{-1}$, we finally get that $h=c^{-1}I=\pm p^{-m}I\in Z$.  \hfill$\square$

{\bf Claim 3.} $\rho$ is surjective.

{\it Proof of Claim 3.} Let $\theta\in Aut(\mathcal R)$. We will prove that $\varepsilon(\theta)=\rho_{h}$, for some $h\in H$.

 Since the action $\Gamma\curvearrowright (G,m_G)$ has spectral gap, Corollary \ref{cor1} implies that after composing $\theta$ with an element from $[\mathcal R]$, we can find open subgroups $G_0,G_1<G$, a continuous isomorphism $\delta:G_0\rightarrow G_1$ and $g\in G$ such that $\delta(\Gamma\cap G_0)=\Gamma\cap G_1$ and $\theta(x)=\delta(x)g^{-1}$, for almost every $x\in G_0$.

Next, by a result of R. Pink, since $G_0,G_1<H=SL_2(\mathbb Q_p)$ are compact open subgroups, $\delta$ extends to a continuous automorphism of $H$ (see \cite[Corollary 0.3]{Pi98}). Since the field $\mathbb Q_p$ has no non-trivial automorphisms, we conclude that there exists $\lambda\in GL_2(\mathbb Q_p)$ such that $\delta(x)=\lambda x\lambda^{-1}$, for all $x\in G_0$.

Since $\delta(\Gamma\cap G_0)=\Gamma\cap G_1$, we get   that $\lambda(\Gamma\cap G_0)\lambda^{-1}\subset\Gamma$.
It follows that $\lambda SL_2(p^n\mathbb Z)\lambda^{-1}\subset\Gamma$, for some $n\geqslant 1$.
Since the subring of $\mathbb M_2(\mathbb Z)$ generated by $SL_2(m\mathbb Z)$ contains $m^2M_2(\mathbb Z)$,  for every $m\in\mathbb Z$, we deduce that 
$p^{2n}\lambda\mathbb M_2(\mathbb Z)\lambda^{-1}\subset \mathbb M_2(\mathbb Z)$.
 If we write $\lambda=\bigl(\begin{smallmatrix}a&b\\ c&d\end{smallmatrix}\bigr)$,  where $a,b,c,d\in\mathbb Q_p$, then we get that $xy\in p^{-2n}\det(\lambda)\mathbb Z$, for all $x,y\in\{a,b,c,d\}$. Thus, after replacing $\lambda$ with $k\lambda$, for some $k\in\mathbb Q_p\setminus\{0\}$, we may assume that $a,b,c,d\in\mathbb Z$, i.e. $\lambda\in GL_2(\mathbb Q_p)\cap\mathbb M_2(\mathbb Z)$.

Moreover, we may assume that the greatest common divisor of $a,b,c,d$ is equal to $1$. Since $a^2,b^2,c^2,d^2\in p^{-2n}\det(\lambda)\mathbb Z$, we get that $\det(\lambda)\mid p^{2n}$, hence $\det(\lambda)=\pm p^m$, for some $m\in\mathbb N$. This shows that $\lambda\in\Lambda$. Since $\theta(x)=\lambda x\lambda^{-1}g$, for almost every $x\in G_0$,  it follows that $\varepsilon(\theta)=\rho_h$, where $h=g^{-1}\lambda\in H$.
\hfill$\square$

Claim 2 and Claim 3 imply that $Out(\mathcal R)\cong H/Z$. Let $\pi:H\rightarrow H/Z$ be the quotient homomorphism. Then $\pi(SL_2(\mathbb Q_p))\cong PSL_2(\mathbb Q_p)$ and   $(H/Z)/PSL_2(\mathbb Q_p)=\{\bigl(\begin{smallmatrix}1&0\\0&x\end{smallmatrix}\bigr)PSL_2(\mathbb Q_p)|\;x=\pm 1,\pm p\}.$ It is now easy to see that $(H/Z)/PSL_2(\mathbb Q_p)\cong (\mathbb Z/2\mathbb Z)^2$.

\vskip 0.05in

(2) Let $\Gamma<SL_2(\mathbb Z)$ be a non-amenable subgroup. Denote by $L$ the closure of $\Gamma$ in $SL_2(\mathbb Z_p)$, and let $\mathcal S=\mathcal R(\Gamma\curvearrowright L)$,
$M=L^{\infty}(L)\rtimes\Gamma$. Our goal is to show that $\mathcal F(\mathcal S)=\mathcal F(M)=\{1\}$. 

 By \cite[Corollary 3]{OP07}  we have that $M$ has a unique Cartan subalgebra, up to unitary conjugacy. Thus, we have that
 $\mathcal F(M)=\mathcal F(\mathcal S)$ and hence it suffices to show that $\mathcal F(\mathcal S)=\{1\}$.

Towards this goal, let $t\in\mathcal F(\mathcal S)$. Note that by Fact \ref{strong},  $L<SL_2(\mathbb Z_p)$ is an open subgroup and the action $\Gamma\curvearrowright (L,m_L)$ has spectral gap. By applying Corollary \ref{cor1} we can find open subgroups $L_0,L_1<L$ such that $t=\frac{[L:L_0]}{[L:L_1]}$ and a continuous homomorphism $\delta:L_0\rightarrow L_1$.

Since $L_0,L_1<SL_2(\mathbb Q_p)$ are compact open subgroups, \cite [Corollary 0.3]{Pi98} implies that there exists $\lambda\in GL_2(\mathbb Q_p)$ such that $\delta(x)=\lambda x\lambda^{-1}$, for all $x\in L_0$.
As in the proof of part (1) we denote by $m_K$  the Haar measure of $K=SL_2(\mathbb Q_p)$. Since $m_K$ is invariant under the conjugation action of $GL_2(\mathbb Q_p)$ by Claim 1, we deduce that $m_K(L_0)=m_K(\lambda L_0\lambda^{-1})=m_K(L_1)$. 

Since $L<K$ is an open subgroup, we have that $m_K(L)>0$. Since $L$ is a compact group, the uniqueness of the Haar measure of $L$ implies that we can find a constant $c>0$ such that $m_K(A)=c\; m_L(A)$, for any Borel subset $A\subset L$.
Thus, we deduce that $m_L(L_0)=m_L(L_1)$ or, equivalently, $[L:L_0]^{-1}=[L:L_1]^{-1}$. This shows that $t=1$ and therefore that $\mathcal F(\mathcal S)=\{1\}$.
\hfill$\blacksquare$

\section{Proofs of Corollaries \ref{corG} and \ref{corH}}\label{tree}

In this section we use Theorem \ref{conjugacy} to calculate the outer automorphism groups of  
equivalence relations arising from the natural  actions of rather general countable subgroups of $\Gamma<SO(n+1)$ on $S^n$ and $P^n(\mathbb R)$. In particular, we derive Corollaries \ref{corG} and \ref{corH}, leading to examples of treeable equivalence relations with trivial outer automorphism group.

We start by fixing some notation:

\begin{notation}
Let $n\geqslant 2$.
\begin{itemize}
\item We denote by $\lambda_n$ the Lebesgue probability measure of the $n$-dimensional sphere, $S^n$, and consider the pmp action $SO(n+1)\curvearrowright (S^n,\lambda_n)$. 

\item We denote by $\mu_n$ the probability measure on the $n$-dimensional real projective space $P^n(\mathbb R)$ obtained by pushing forward $\lambda_n$ through the quotient map $S^n\rightarrow P^n(\mathbb R)$: $\xi\rightarrow [\xi]$.

\item We consider the action $\mathbb Z/2\mathbb Z\curvearrowright (S^n,\lambda_n)$ given by  the involution $T(x)=-x$. 
Then $(S^n,\lambda_n)/(\mathbb Z/2\mathbb Z)=(P^n(\mathbb R),\mu_n)$. 
 Moreover, since the  actions of $\mathbb Z/2\mathbb Z$ and $SO(n+1)$ on $S^n$ commute, we have a pmp action $SO(n+1)\curvearrowright (P^n(\mathbb R),\mu_n)$.
\end{itemize}
\end{notation}

\begin{theorem}\label{OUT} Let  $\Gamma<G:=SO(n+1)$ be a countable icc dense subgroup, for some $n\geqslant 2$.  Assume that the left translation action $\Gamma\curvearrowright (G,m_G)$ has spectral gap. 
\begin{enumerate}
\item Let $\mathcal R=\mathcal R(\Gamma\curvearrowright (G,m_G))$. Then we have:
\begin{enumerate}
\item $\mathcal F(\mathcal R)=\{1\}$. 
\item If $n$ is even, then  $Out(\mathcal R)\cong N_G(\Gamma)/\Gamma\times G$.
\item If $n$ is odd, then $Out(\mathcal R)\cong (N _{G}(\Gamma)\times G)/\tilde\Gamma$, where $\tilde\Gamma=\{(\alpha\gamma,\alpha I)|\gamma\in\Gamma,\alpha=\pm 1\}$ and $I\in G$ is the identity matrix.
\end{enumerate}

\item Let $\mathcal S=\mathcal R(\Gamma\curvearrowright (S^n,\lambda_n))$. Then we have:

\begin{enumerate}
\item $\mathcal F(\mathcal S)=\{1\}.$
\item If $n$ is even, then Out$(\mathcal S)\cong N_G(\Gamma)/\Gamma\times(\mathbb Z/2\mathbb Z)$.
\item If $n$ is odd, then $Out(\mathcal S)\cong (N_G(\Gamma)\times\tilde K)/K_0$, where $K=\{\begin{pmatrix} 1 \;\;0\\0 \;\;a \end{pmatrix}|\;a\in SO(n)\}$, $\tilde K=\{\begin{pmatrix} \alpha \;\;0\\0 \;\;\alpha a \end{pmatrix}|\;a\in SO(n),\alpha=\pm 1\}$, and   $K_0=\{(\alpha\gamma,\alpha k)|\alpha=\pm 1,\gamma\in\Gamma,k\in K\}$.

\end{enumerate}

\item Let $\mathcal T=\mathcal R(\Gamma\curvearrowright (P^n(\mathbb R),\mu_n))$ and assume that $-I\not\in\Gamma$. Then we have:

\begin{enumerate}
\item $\mathcal F(\mathcal T)=\{1\}$. 
\item If $n$ is even, then Out$(\mathcal T)\cong N_G(\Gamma)/\Gamma$.
\item If $n$ is odd, then Out$(\mathcal T)\cong N_G(\Gamma)/\bar{\Gamma}$, where $\bar{\Gamma}=\{\alpha\gamma|\alpha=\pm 1,\gamma\in\Gamma\}$.
\end{enumerate}

\end{enumerate}

\end{theorem}

\begin{remark} 
J. Bourgain and A. Gamburd proved that if $\Gamma<H=SU(2)$ is a dense subgroup  generated by finitely many matrices $\{g_1,...,g_k\}$ having algebraic entries, then the representation $\pi:\Gamma\rightarrow \mathcal U(L^2(H,m_H))$ (equivalently, the action $\Gamma\curvearrowright (H,m_H)$) has spectral gap \cite{BG06}. More recently, they proved that this result holds for $H=SU(n)$, whenever $n\geqslant 2$ \cite{BG11}.
In particular, it follows that if $\Gamma<G=SO(3)$ is any dense subgroup which is generated by finitely many matrices with algebraic entries, then the action $\Gamma\curvearrowright (G,m_G)$ has spectral gap.
This provides a large family of actions to which the above theorem applies.

\end{remark}

{\bf Proof.} Let us first record a fact that we will use repeatedly: $(\star$) $\Gamma\curvearrowright (G,m_G)$ has spectral gap, $G$ is connected, $\pi_1(G)\cong\mathbb Z/2\mathbb Z$
 is finite, and $\Gamma$ is icc, hence has no non-trivial normal subgroups.
 \vskip 0.05in
By using $(\star)$ and applying Theorem \ref{conjugacy} to $H=G$, $K=\{e\}$, $L=\{e\}$, we get that $\mathcal F(\mathcal R)=\{1\}$. Moreover, Theorem \ref{conjugacy} shows that if $\theta\in$ Aut$(\mathcal R)$, then after composing $\theta$ with an element of $[\mathcal R$], we can find an automorphism $\delta:G\rightarrow G$ with $\delta(\Gamma)=\Gamma$, and $y\in G$, such that $\theta(g)=\delta(g)y$, for almost every $g\in G$. Since $G$ has no outer automorphisms,  $\delta(g)=zgz^{-1}$, for some $z\in N_G(\Gamma)$.
Therefore,  if $\theta:N_G(\Gamma)\times G\rightarrow$ Aut$(\mathcal R)$ denotes the homomorphism defined by $\theta_{z,w}(g)=zgw$, and  $\varepsilon:$ Aut$(\mathcal R)\rightarrow$ Out$(\mathcal R)$ denotes the natural quotient homomorphism, then $\varepsilon\circ\theta$ is surjective.

Let $z\in N_G(\Gamma)$ and $w\in G$ such that $\theta_{z,w}\in [\mathcal R]$. Let $\gamma\in\Gamma$ such that the set $A=\{g\in G|zgw=\gamma g\}$ has positive measure. Then $zgh^{-1}z^{-1}=\gamma gh^{-1}\gamma^{-1}$, for all $g,h\in A$. Hence $\gamma^{-1}z$ commutes with $AA^{-1}$ and further with $G_0:=\cup_{n\geqslant 1}(AA^{-1})^n$. Since $A$ has positive measure, $AA^{-1}$ contains a neighborhood of $I\in G$.  Thus, $G_0$ is an open subgroup of $G$. Since $G$ is connected,  we get that $G_0=G$ and therefore $\gamma^{-1}z$ must be in the center of $G$.

If $n$ is even, then $Z(G)=\{I\}$ and hence $z=\gamma$. Since there is $g\in G$ such that $zgw=\gamma g$, we  further get that $w=1$. This shows that the kernel of $\varepsilon\circ\theta$ is equal to $\Gamma\times\{I\}$, hence $Out(\mathcal R)\cong N_G(\Gamma)/\Gamma\times G$. If $n$ is odd, then $Z(G)=\{\pm I\}$. It follows that $z=\alpha\gamma$ and $w=\alpha$, for $\alpha\in \{\pm I\}$. This shows that the kernel of $\varepsilon\circ\theta$ is equal to $\tilde\Gamma$,  hence $Out(\mathcal R)\cong (N _{G}(\Gamma)\times G)/\tilde\Gamma$. 

\vskip 0.05in
This finishes the proof of part (1) of Theorem \ref{OUT}.
We now turn to the proofs of parts (2) and (3).

We begin by giving  an alternative description of the actions of $G=SO(n+1)$ on $S^n$ and $P^n(\mathbb R)$.
Let $\xi=(1,0,...,0)\in S^n$. Since the action $G\curvearrowright S^n$ is transitive and Stab$_G(\xi)=K$, it follows that
the action $G\curvearrowright (S^n,\lambda_n)$ is isomorphic to the action $G\curvearrowright (G/K,m_{G/K})$. 

Also, notice that $\tilde K=\{g\in G|\;g\xi=\pm\xi\}$ and $K<\tilde K$ is an open normal subgroup of index two. Since the action $G\curvearrowright P^n(\mathbb R)$ is transitive and Stab$_G([\xi])=\tilde K$, it follows that
the action $G\curvearrowright (P^n(\mathbb R),\mu_n)$ is isomorphic to the action $G\curvearrowright (G/\tilde K,m_{G/\tilde K})$. 

We are now ready to calculate $\mathcal F(\mathcal S)$ and Out$(\mathcal S)$. Note that $(\star)$ ensures that conditions (1)-(3) of Theorem \ref{conjugacy} 
are satisfied for $\Lambda=\Gamma$, $H=G$ and $L=K$. Since condition (4) is also satisfied by Claim 2, we are in position to apply Theorem \ref{conjugacy}. Thus, we deduce that
$\mathcal F(\mathcal S)=\{1\}$. Moreover, we get that if $\theta\in$ Aut$(\mathcal S)$, then after composing $\theta$ with an element from $[\mathcal S]$, we can find an automorphism $\delta:G\rightarrow G$ and $y\in G$ such that $\delta(\Gamma)=\Gamma$, $\delta(K)=yKy^{-1}$, and $\theta(gK)=\delta(g)yK$, for almost every $g\in G$.

Since $G$ has no outer automorphisms, we can find $z\in G$ such that $\delta(g)=zgz^{-1}$, for all $g\in G$.
Thus, if we put $w=z^{-1}y$, then $z\in N_G(\Gamma)$, $w\in N_G(K)$, and $\theta(gK)=zgwK$, for almost every $g\in G$.
 Since $K$ stabilizes $\xi\in S^n$ and $w$ normalizes $K$, it follows that $K$ stabilizes $w\xi\in S^n$. This easily implies that $w\xi=\pm\xi$, i.e. $w\in\tilde K$.

Consider the well-defined homomorphism $\theta:N_G(\Gamma)\times\tilde K\rightarrow$ Aut$(\mathcal S)$ given by $\theta_{z,w}(gK)=zgwK$. Then the above shows that $\varepsilon\circ\theta:N_G(\Gamma)\times\tilde K\rightarrow$ Out$(\mathcal S)$ is surjective, where $\varepsilon:$ Aut$(\mathcal S)\rightarrow$ Out$(\mathcal S)$ denotes the quotient homomorphism.

Let $z\in N_G(\Gamma)$ and $w\in\tilde K$ such that $\theta_{z,w}\in [\mathcal S]$. Then we can find $\gamma\in\Gamma$ such that the set $A=\{gK\in G/K|zgwK=\gamma gK\}$ has positive measure.
Let $\alpha\in\{\pm 1\}$ such that $w\xi=\alpha\xi$. Then for every $gK\in A$ we have that $\alpha z(g\xi)=\gamma(g\xi)$. Since $m_{G/K}(A)>0$, Claim 1 implies that $\alpha z=\gamma$. Since there exists $gK\in G/K$ such that $zgwK=\gamma gK$, we derive that  $wK=\alpha K$ and hence $\alpha I\in G$. 

If $n$ is even, then $-I\not\in G$, which forces $\alpha=1$. Hence, $w\in K$ and $z=\gamma\in\Gamma$. This gives that the kernel of $\varepsilon\circ\theta$ is equal to $\Gamma\times K$ and thus Out$(\mathcal S)\cong (N_G(\Gamma)\times\tilde K)/(\Gamma\times K)\cong N_G(\Gamma)/\Gamma\times(\mathbb Z/2\mathbb Z)$. If $n$ is odd, then $-I\in G$ and we get that the kernel of $\varepsilon\circ\theta$ is equal to $K_0$. This finishes the proof of part (2) of Theorem \ref{OUT}.

Finally, we compute $\mathcal F(\mathcal T)$ and Out$(\mathcal T)$. Recall that we may identify $\mathcal T=\mathcal R(\Gamma\curvearrowright G/\tilde K)$. 
By $(\star)$, conditions (1)-(3) of Theorem \ref{conjugacy} are satisfied for $H=G$, $K$ and $L$ both equal to $\tilde K$, and $\Lambda=\Gamma$. 
Since $-I\not\in\Gamma$,  Claim 2 shows that condition (4) of Theorem \ref{conjugacy} is also verified. 

By applying Theorem \ref{conjugacy} we get that $\mathcal F(\mathcal T)=\{1\}$.
Moreover,  it follows that if $\theta\in$ Aut$(\mathcal S)$, then after composing $\theta$ with an element from $[\mathcal S]$, we can find $z\in N_G(\Gamma)$ and $w\in N_G(\tilde K)$ such that $\theta(g\tilde K)=zgw\tilde K$, for almost every $g\in G$. Since $\tilde K$ stabilizes $[\xi]\in P^n(\mathbb R)$ and $w$ normalizes $\tilde K$, we have that $\tilde K$ stabilizes $[w\xi]\in P^n(\mathbb R)$.  This immediately implies that $[w\xi]=[\xi]$, hence $w\in\tilde K$. As a consequence, we have that $\theta(g\tilde K)=zg\tilde K$, for almost every $g\in G$.

Define  $\theta:N_G(\Gamma)\rightarrow$ Aut$(\mathcal T)$ by letting $\theta_z(g\tilde K)=zg\tilde K$. If $\varepsilon:$ Aut$(\mathcal T)\rightarrow$ Out$(\mathcal T)$ denotes the quotient homomorphism, then the above implies that $\varepsilon\circ\theta$ is surjective.

 Let $z\in N_G(\Gamma)$ such that $\theta_z\in [\mathcal T]$. Then  $A=\{g\tilde K\in G/\tilde K|zg\tilde K=\gamma g\tilde K\}$ has positive measure, for some $\gamma\in\Gamma$. If $g\tilde K\in A$, then $z(g\xi)=\pm\gamma(g\xi)$. By using Claim 1, we get that $z=\alpha\gamma$, for $\alpha\in\{\pm 1\}$. 

If $n$ is even, then $-I\not\in G$, so we must have $\alpha=1$, hence $z=\gamma\in\Gamma$. This shows that the kernel of $\theta$ is equal to $\Gamma$ and therefore Out$(\mathcal T)\cong N_G(\Gamma)/\Gamma$. 
If $n$ is odd, then $-I\in G$ and we get that $\ker(\varepsilon\circ\theta)\subset\bar{\Gamma}$. Since $\bar{\Gamma}$ is  contained in the kernel of $\varepsilon\circ\theta$, we deduce that Out$(\mathcal T)\cong N_G(\Gamma)/\bar{\Gamma}$. 
\hfill$\blacksquare$

\subsection{Proof of Corollary \ref{corG}} 
Let $\Gamma<G$ be a subgroup which contains matrices $\{g_1,...,g_k\}$ with algebraic entries such that the subgroup $\Gamma_0$ of $\Gamma$ generated by $\{g_1,...,g_k\}$ is dense in $G$. Let us briefly explain how \cite{BG06} implies that the action $\Gamma\curvearrowright (G,m_G)$ has spectral gap.

Recall that there is a surjective continuous homomorphism $\Phi:H=SU(2)\rightarrow G=SO(3)$ whose kernel is $\{\pm I\}$ (see e.g. \cite[Section 1.6.1]{Ha03}). More precisely,  $\Phi$ is given by 
$$\Phi(\;(\begin{matrix} x & y\\-\bar{y} & \bar{x}\end{matrix})\;)=\begin{pmatrix} \Re(x^2-y^2) & \Im(x^2+y^2) & -2\;\Re(xy)\\ -\Im(x^2-y^2) & \Re(x^2+y^2) &  2\;\Im(xy)\\ 2\;\Re(x\bar{y}) & 
2\;\Im(x\bar{y}) & |x|^2-|y|^2\end{pmatrix}\;\;\text{for all $x,y\in\mathbb C$ with $|x|^2+|y|^2=1$}.$$

Let $h_1,...,h_k\in H$ such that $\Phi(h_1)=g_1,...,\Phi(h_k)=g_k$, and denote by $\Lambda<H$ the subgroup generated by $\{h_1,...,h_k\}$. Since $\Gamma$ is dense in $G$, the closure $H_0=\bar{\Lambda}<H$ satisfies $\Phi(H_0)=G$. Since $\Phi$ is a $2$-$1$ map, we get that $[H:H_0]\leqslant 2$. Since $H$ is connected, we must have $H_0=H$, hence $\Lambda$ is dense in $H$. Moreover, since the entries of $g_1,...,g_k$ are algebraic, the above formula for $\Phi$ implies that the entries of $h_1,...,h_k$ are also algebraic.

By a result of J. Bourgain and A. Gamburd \cite[Theorem 1]{BG06}, since $\Lambda<H$ is a finitely generated, dense subgroup, which is generated by matrices $\{h_1,...,h_k\}$ with algebraic entries, the action $\Lambda\curvearrowright (H,m_{H})$ has spectral gap. This readily implies that the action $\Gamma\curvearrowright (G,m_G)$ has spectral gap.

Corollary \ref{corG} now follows directly from Theorem \ref{OUT}.
\hfill$\blacksquare$

We continue  by proving the following more general form of Corollary \ref{corH}.

\begin{corollary}\label{cor}
Let $p,q\geqslant 3$ be natural numbers. Denote by $\alpha_p$ the rotation about the $x$-axis by angle $\frac{2\pi}{p}$ and by $\beta_q$ the rotation about 
the $z$-axis by angle $\frac{2\pi}{q}$, i.e.

$$\alpha_p=\begin{pmatrix} 1 & 0 & 0\\ 0 & \cos{\frac{2\pi}{p}} & -\sin{\frac{2\pi}{p}}\\ 0 & \sin{\frac{2\pi}{p}} & \cos{\frac{2\pi}{p}}\end{pmatrix}\;\;\;\;\text{and}\;\;\;\;\;
 \beta_q=\begin{pmatrix}  \cos{\frac{2\pi}{q}} & -\sin{\frac{2\pi}{q}} & 0 \\ \sin{\frac{2\pi}{q}} & \cos{\frac{2\pi}{q}} & 0 \\ 0 & 0 & 1\end{pmatrix}.$$
 
Denote by $\Gamma=\Gamma(p,q)$ the subgroup of $G=SO(3)$ generated by $\alpha_p$ and $\beta_q$.

Then the following hold: 
 \begin{enumerate}
 \item If $p$, $q$ are odd and $p\not=q$, then  $\Gamma\cong (\mathbb Z/p\mathbb Z)*(\mathbb Z/q\mathbb Z)$,    $Out(\mathcal R(\Gamma\curvearrowright G))\cong (\mathbb Z/2\mathbb Z)^{3}\times G$
  $Out(\mathcal R(\Gamma\curvearrowright S^2))\cong (\mathbb Z/2\mathbb Z)^{3}$
   and $Out(\mathcal R(\Gamma\curvearrowright P^2(\mathbb R))\cong (\mathbb Z/2\mathbb Z)^2$.
  
  \item If $p\geqslant 4$ is even and $q\geqslant 3$ is odd, then  $\Gamma\cong (\mathbb Z/p\mathbb Z)*_{\mathbb Z/2\mathbb Z}D_q$, $Out(\mathcal R(\Gamma\curvearrowright G))\cong (\mathbb Z/2\mathbb Z)^{2}\times G$,
  $Out(\mathcal R(\Gamma\curvearrowright S^2))\cong (\mathbb Z/2\mathbb Z)^2$ and $Out(\mathcal R(\Gamma\curvearrowright P^2(\mathbb R))\cong \mathbb Z/2\mathbb Z$.
  
  \item If $p\geqslant 4$ is even, $q=2s$, $s\geqslant 3$ odd, and $p\not=q$, then   $\Gamma\cong D_p*_{D_2}D_q$ and we have that $Out(\mathcal R(\Gamma\curvearrowright G))\cong (\mathbb Z/2\mathbb Z)\times G$,
  $Out(\mathcal R(\Gamma\curvearrowright S^2))\cong \mathbb Z/2\mathbb Z$ and $Out(\mathcal R(\Gamma\curvearrowright P^2(\mathbb R))\cong\{e\}$.
\end{enumerate}
\end{corollary}

Here, $D_n$ denotes the dihedral group with $2n$ elements.
More precisely, $D_n$ is the semi-direct product $D_n=(\mathbb Z/n\mathbb Z)\rtimes \mathbb Z/2\mathbb Z$ associated with the order two automorphism $x\mapsto -x$ of $\mathbb Z/n\mathbb Z$.

The isomorphisms between $\Gamma$ and the corresponding amalgamated free product groups are due to C. Radin and L. Sadun \cite[Corollary 2]{RS98}. Notice that they  imply that $\Gamma$ is icc. 

{\it Proof.} 
Note that  $\alpha_p$, $\beta_q$ have algebraic entries, and $\Gamma<G=SO(3)$ is dense (since it contains a copy of $\mathbb F_2$). By applying Corollary \ref{corG} we deduce that $Out(\mathcal R(\Gamma\curvearrowright G))\cong (N_G(\Gamma)/\Gamma)\times G$, $Out(\mathcal R(\Gamma\curvearrowright S^2))\cong N_G(\Gamma)/\Gamma\times(\mathbb Z/2\mathbb Z)$ and
 $Out(\mathcal R(\Gamma\curvearrowright P^2(\mathbb R))\cong N_G(\Gamma)/\Gamma.$

In the rest of the proof, we calculate $N_G(\Gamma)$ is each of the three cases. Towards this, we let $g\in N_G(\Gamma)$ and denote by $\rho$ the automorphism of $\Gamma$ given by $\rho(x)=gxg^{-1}$. Also, we define $$a=\begin{pmatrix}  1 & 0 & 0 \\ 0 & -1 & 0 \\ 0 & 0 & -1 \end{pmatrix}\;\;\;\;\text{and}\;\;\;\;b=\begin{pmatrix} -1 & 0 & 0 \\ 0 & -1 & 0 \\ 0 & 0 & 1\end{pmatrix}.$$

(1) Assume that $p$, $q$ are  odd and  $p\not=q$. Denote  $\Gamma_1=\langle\alpha_p\rangle\cong\mathbb Z/p\mathbb Z$ and $\Gamma_2=\langle\beta_q\rangle\cong \mathbb Z/q\mathbb Z$. Then $\Gamma=\Gamma_1*\Gamma_2$ by \cite[Corollary 2]{RS98}. The Kurosh subgroup theorem implies that we can find $h_1,h_2\in\Gamma$ and $i_1,i_2\in\{1,2\}$ such that $\rho(\Gamma_1)\subset h_1\Gamma_{i_1}h_1^{-1}$  and 
$\rho(\Gamma_2)\subset h_2\Gamma_{i_2}h_2^{-1}$. Since $\Gamma_i$ and $h\Gamma_ih^{-1}$ cannot generate $\Gamma$, for any $h\in\Gamma$ and $i\in\{1,2\}$,  we conclude that $i_1\not=i_2$. Since $p\not=q$, we must have that $i_1=1$ and $i_2=2$.
  It follows that  $\rho(\Gamma_1)=h_1\Gamma_1h_1^{-1}$ and $\rho(\Gamma_2)=h_2\Gamma_2h_2^{-1}$. 
  
  Since $h_1\Gamma_1h_1^{-1}$ and $h_2\Gamma_2h_2^{-1}$ generate $\Gamma$, we can find $k_1\in\Gamma_1,k_2\in \Gamma_2$ such that $h_2^{-1}h_1=k_2^{-1}k_1$. Denote $l=h_1k_1^{-1}=h_2k_2^{-1}\in\Gamma$.  Then $\rho(\Gamma_1)=l\Gamma_1l^{-1}$ and $\rho(\Gamma_2)=l\Gamma_2l^{-1}$.
  
   Therefore, we have that $l^{-1}g\in N_G(\Gamma_1)\cap N_G(\Gamma_2)$. 
  It is easy to see that \begin{equation}\label{norma}N_G(\Gamma_1)=\{\begin{pmatrix} \det{A} & 0 \\ 0 & A\end{pmatrix}\;|\;A\in O(2)\}\;\;\;\text{and}\;\;\;N_G(\Gamma_1)=\{\begin{pmatrix} A & 0 \\ 0 & \det{A}\end{pmatrix}\;|\; A\in O(2)\}.\end{equation}

If we let $D<G$ be the subgroup consisting of diagonal matrices, then  equation \ref{norma} implies that $D=N_G(\Gamma_1)\cap N_G(\Gamma_2)$. Hence $l^{-1}g\in D$, showing that $g\in\langle\Gamma,D\rangle$. In conclusion,  $N_G(\Gamma)=\langle\Gamma,D\rangle$.

Finally, notice that $D=\{I, a, b, ab\}$,  $a=\alpha_{2p}^p$, $b=\beta_{2q}^q, \alpha_p=\alpha_{2p}^2$, $\beta_q=\beta_{2q}^2$. Therefore, we get that $\langle\Gamma,D\rangle=\Gamma(2p,2q)$. Moreover, \cite[Corollary 2]{RS98} gives that $\Gamma(2p,2q)=\langle\alpha_p,b\rangle*_{\langle a,b\rangle}\langle\beta_q,a\rangle$. Since $\Gamma=\Gamma(p,q)=\langle\alpha_p\rangle*\langle\beta_q\rangle$, it follows that $N_G(\Gamma)/\Gamma=\Gamma(2p,2q)/\Gamma(p,q)\cong(\mathbb Z/2\mathbb Z)\times (\mathbb Z/2\mathbb Z)$.

(2) Assume that $p\geqslant 4$ is even and $q\geqslant 3$ is odd. Note that $\alpha_p^{\frac{p}{2}}=a$ normalizes the cyclic group $\langle\beta_q\rangle$. Denote $\Gamma_1=\langle\alpha_p\rangle\cong\mathbb Z/p\mathbb Z$, $\Gamma_2=\langle\beta_q,a\rangle\cong D_q$ and $\Lambda=\langle a\rangle\cong \mathbb Z/2\mathbb Z$. Then by \cite[Corollary 2]{RS98} we have that $\Gamma=\Gamma_1*_{\Lambda}\Gamma_2$. 

Since  $\Gamma_1\not\cong\Gamma_2$, by reasoning as in the proof of (1) we can find $l\in\Gamma$ such that  $\rho(\Gamma_1)=l\Gamma_1l^{-1}$ and $\rho(\Gamma_2)=l\Gamma_2l^{-1}$. Hence, $l^{-1}g\in N_G(\Gamma_1)\cap N_G(\Gamma_2)$. It is easy to see that  $N_G(\Gamma_1)\cap N_G(\Gamma_2)=D$, which implies that $g\in\langle\Gamma, D\rangle$ and hence $N_G(\Gamma)=\langle\Gamma,D\rangle$. 

Since $D=\{I, a, b, ab\}$, $a=\alpha_p^{\frac{p}{2}}\in\Gamma$, $b=\beta_{2q}^q$ and $\beta_q=\beta_{2q}^2$, we deduce that $\langle\Gamma,D\rangle=\Gamma(p,2q)$. By \cite[Corollary 2]{RS98} we have that $\Gamma(p,2q)=\langle\alpha_p,b\rangle*_{\langle a,b\rangle}\langle\beta_{2q},a\rangle$. Since $\Gamma(p,q)=\langle\alpha_p\rangle*_{\langle a\rangle}\langle \beta_{2q}^2,a\rangle$, we finally conclude that $N_G(\Gamma)/\Gamma=\Gamma(p,2q)/\Gamma(p,q)\cong\mathbb Z/2\mathbb Z.$

(3) Assume that $p\geqslant 4$ is even, $q=2s$, $s$ odd, and $p\not=q$. Note that $\alpha_p^{\frac{p}{2}}=a$ normalizes the cyclic group $\langle\beta_q\rangle$ and $\beta_q^{\frac{q}{2}}$ normalizes the cyclic group $\langle\alpha_p\rangle$. If we denote $\Gamma_1=\langle\alpha_p,b\rangle\cong D_p$,  $\Gamma_2=\langle\beta_q,a\rangle\cong D_q$ and $\Lambda=\langle a,b\rangle\cong (\mathbb Z/2\mathbb Z)^2$, then by \cite[Corollary 2]{RS98} we have that $\Gamma=\Gamma_1*_{\Lambda}\Gamma_2$.

Since  $\Gamma_1\not\cong\Gamma_2$, by reasoning as in the proof of (1) we can find $l\in\Gamma$ such that  $\rho(\Gamma_1)=l\Gamma_1l^{-1}$ and $\rho(\Gamma_2)=l\Gamma_2l^{-1}$. Hence, $l^{-1}g\in N_G(\Gamma_1)\cap N_G(\Gamma_2)$. It is easy to see that  $N_G(\Gamma_1)\cap N_G(\Gamma_2)=D$, which implies that $g\in\langle\Gamma, D\rangle$ and hence $N_G(\Gamma)=\langle\Gamma,D\rangle$. Since $D\subset\Gamma$, we get that $N_G(\Gamma)=\Gamma$.
\hfill$\blacksquare$

Finally, let us note the the actions $\Gamma(p,q)\curvearrowright G$ are neither stably orbit equivalent nor Borel reducible to each other, for varying values of $(p,q)$. More generally, we have:

\begin{corollary}
Let $\Gamma,\Lambda$ be countable icc dense subgroups of $G=SO(3)$. Assume that $\Gamma$ contains matrices $g_1,...,g_k$  which have algebraic entries and generate a dense subgroup of $G$. Then we have:

\begin{enumerate}
\item The actions $\Gamma\curvearrowright (G,m_G)$ and $\Lambda\curvearrowright (G,m_G)$ are stably orbit equivalent if and only if there exists $g\in G$ such that $g\Gamma g^{-1}=\Lambda$.
\item $\mathcal R(\Gamma\curvearrowright G)\leqslant_{B}\mathcal R(\Lambda\curvearrowright G)$ if and only if there exists $g\in G$ such that $g\Gamma g^{-1}=\Lambda$.
\item There exists a non-trivial homomorphism from $\mathcal R(\Gamma\curvearrowright G)$ to $\mathcal R(\Lambda\curvearrowright G)$ if and only if there exists $g\in G$ such that $g\Gamma g^{-1}\subset\Lambda$.

\end{enumerate}

Moreover, let $p,q,r,s\geqslant 3$ be integers such that $(p,q)\not=(r,s), (p,q)\not=(s,r)$ and $(p,q),(r,s)\not\in (4\mathbb Z)^2$.
Then the actions $\Gamma(p,q)\curvearrowright (G,m_G)$ and $\Gamma(r,s)\curvearrowright (G,m_G)$ are not stably orbit equivalent, and the equivalence relations  $\mathcal R(\Gamma(p,q)\curvearrowright G)$ and $\mathcal R(\Gamma(r,s)\curvearrowright G)$ are not comparable with respect to Borel reducibility.

\end{corollary}

{\it Proof.} As in the proof of Corollary \ref{corG}, the main result of \cite{BG06} implies that the action $\Gamma\curvearrowright (G,m_G)$ has spectral gap.
We claim that if  $\delta:G\rightarrow H/\Delta$ is a non-trivial continuous homomorphism, where  $H<G$ is a subgroup and $\Delta<Z(H)$ is a subgroup, then $H=G$, $\Delta=\{e\}$ and there exists $g\in G$ such that $\delta(x)=gxg^{-1}$, for all $x\in G$. This is because $G$ is a simple group, any proper subgroup $H\lneq G$ has dimension strictly smaller than $G$, and $G$ has no outer automorphisms.

By applying Corollary \ref{OE} and Corollary \ref{Bor}, the assertions (1)-(3) follow immediately.

To see the moreover assertion, note that the assumptions made and \cite[Corollary 2]{RS98} imply that $\Gamma(p,q)$, $\Gamma(r,s)$ are non-isomorphic icc groups. By applying (1) and (2) we get the conclusion. \hfill$\blacksquare$

\begin{remark}\label{Swi}
In \cite{Sw94}, S. \'{S}wierczkowski exhibits an interesting family of embeddings of $\mathbb F_2$ into $G=SO(3)$. Let $a,b$ be integers such that $b>0$ and $|a|\leqslant b$. Put $c=b^2-a^2$ and define $\Gamma$ to be the subgroup of $G$ generated by the following rotation matrices: $$A=\begin{pmatrix} \frac{a}{b}&-\frac{\sqrt{c}}{b}&0\\\frac{\sqrt{c}}{b}&\frac{a}{b}&0\\0&0&1 \end{pmatrix}\;\;\;\text{and}\;\;\; B=\begin{pmatrix} 1&0&0\\0&\frac{a}{b}&-\frac{\sqrt{c}}{b}\\0&\frac{\sqrt{c}}{b}&\frac{a}{b}\end{pmatrix}.$$

The main result of \cite{Sw94} asserts that if $\frac{a}{b}\notin\{0,\pm\frac{1}{2},\pm 1\}$, then $\Gamma=\langle A\rangle*\langle B\rangle\cong\mathbb F_2$.
Thus,  if $\frac{a}{b}\notin\{0,\pm\frac{1}{2},\pm 1\}$, then since the entries of $A$ and $B$ are algebraic, Corollary \ref{corG} applies to $\Gamma$ and, more generally, to any non-cyclic subgroup $\Gamma_0<\Gamma$. Therefore, the calculation of the outer automorphism groups of the equivalence relations associated to the actions of $\Gamma$  on $G, S^2, P^2(\mathbb R)$ reduces to the calculation of $N_G(\Gamma)$. We were, however, unable to compute $N_G(\Gamma)$ or, more generally, $N_G(\Gamma_0)$.
\end{remark}

\end{document}